\documentclass[mathpazo]{cicp}

\usepackage[dvipsnames]{xcolor}
\usepackage{subcaption}
\usepackage{graphicx}
\usepackage{multirow}

\newcommand{\cU}{\mathcal{U}}
\newcommand{\cD}{\mathcal{D}}
\newcommand{\cY}{\mathcal{Y}}
\newcommand{\cK}{\mathcal{K}}
\newcommand{\cM}{\mathcal{M}}
\newcommand{\cN}{\mathcal{N}}
\newcommand{\cF}{\mathcal{F}}

\newcommand{\cB}{\mathcal{B}}
\newcommand{\U}{\mathbb{U}}
\newcommand{\V}{\mathbb{V}}
\newcommand{\W}{\mathbb{W}}

\newcommand{\NN}{\mathcal{NN}}
\newcommand{\eRNN}{Exp-ResNN}
\newcommand{\RNN}{Gl-ResNN}
\newcommand{\ba}{{\bf a}}
\newcommand{\bc}{{\bf c}}

\newcommand{\bfell}{{\boldsymbol \ell}}
\newcommand{\bPhi}{{\boldsymbol \Phi}}
\newcommand{\bPsi}{{\boldsymbol \Psi}}

\newcommand{\bo}{{\bf o}}

\newcommand{\bw}{{\bf w}}

\newcommand{\bx}{{\bf x}}
\newcommand{\by}{{\bf y}}
\newcommand{\bz}{{\bf z}}

\newcommand{\R}{\mathbb{R}}
\newcommand{\cE}{\mathbb{E}}

\newcommand{\bzero}{{\bf 0}}

\newcommand{\mE}{\mathcal{E}}
\newcommand{\mw}{\textcolor{black}}
\newcommand{\mwrevise}{\textcolor{black}}

\newcommand{\be}{\begin{equation}}
\newcommand{\ee}{\end{equation}}
\newcommand{\Wp}{{\widetilde{\W}_h^\perp}}
\newcommand{\e}{\varepsilon}

\renewcommand\div{\mathop{\rm div}}
\newcommand\dist{\mathop{\rm dist}}

\def\argmin{\mathop{\rm argmin}}

\newtheorem{rem}{Remark}

\begin{document}
\title{Nonlinear Reduced DNN Models for State Estimation}


\author[Dahmen W et.~al.]{Wolfgang Dahmen\affil{1},
      Min Wang\affil{2}\comma\corrauth, and Zhu Wang\affil{1}}
\address{\affilnum{1}\ Department of Mathematics, University of South Carolina, Columbia, SC 29208. \\
          \affilnum{2}\ Department of Mathematics, Duke University, Durham, NC 27708.}
\emails{{\tt dahmen@math.sc.edu} (W.~Dahmen), {\tt wangmin@math.duke.edu} (M.~Wang),
         {\tt wangzhu@math.sc.edu} (Z.~Wang)}

\begin{abstract}
We propose in this paper a data driven state estimation scheme for generating nonlinear reduced models for parametric families of PDEs, directly providing {\em data-to-state maps}, represented in terms of Deep Neural Networks. A major constituent is a {\em sensor-induced} decomposition
of a model-compliant Hilbert space warranting approximation in problem relevant metrics. It plays a similar role as in   a Parametric Background Data Weak framework for state estimators based on Reduced Basis concepts. Extensive numerical
tests shed light on several optimization strategies that are to improve robustness
and performance of such estimators.

\end{abstract}


\maketitle

\section{Introduction}
\label{sec1}

Understanding complex ``physical systems'' solely through observational data
is an attractive but unrealistic objective if one insists on certifiable
accuracy quantification. This,  in turn, is an essential precondition for 
prediction capability. In fact, unlike application scenarios where
an abundance of data are available, data acquisition for   ``Physics Informed Learning
Task'' 
typically relies on sophisticated sensor technology and is often expensive 
or even harmful.
Therefore, a central task is to develop efficient ways for fusing the information
provided by data with {\em background information} provided by
physical laws governing the observed states of interest, typically represented
by partial differential equations (PDEs). In principle, this falls into the framework
of ``Physics Informed Neural Networks'' (PINN), however, with some noteworthy distinctions
explained next.

The central objective of this note is to explore a machine learning approach
to state estimation in the above sense. Our contributions concern two major aspects:

(i) In contrast to typical PINN formulations, we   employ  loss functions that 
are {\em equivalent} to the error of the estimator  in a norm that is imposed by the continuous 
mathematical 
model. More precisely, this norm corresponds to a {\em stable variational formulation}
of the PDE family. In other words, the generalization error for this loss function
measures the accuracy of the estimator in a problem intrinsic norm without
imposing any additional regularity properties.

(ii) When employing estimators, represented
  as Deep Neural Networks (DNNs),  one has to accept a significant and unavoidable uncertainty
  about optimization success. Due to (i), one can at least measure the achieved accuracy at
  any stage of the optimization. We therefore take this fact as a starting point for
  a systematic {\em computational exploration} of a simple  optimization strategy that seems to
  be particularly natural in combination with ResNet architectures. 
  
  Regarding (i), the proposed approach is, in principle,  applicable to a much wider scope
  of problems than discussed below. Last but not least, in order to facilitate comparisons 
  with other recovery schemes, specifically with methods that are based on {\em Reduced Basis}
  concepts, the numerical experiments focus on elliptic families of PDEs with {\em parameter dependent} diffusion fields. However, for this problem class we discuss in detail two rather different scenarios,
  namely diffusion coefficients with an {\em affine} parameter dependence, as well as {\em log-normal}
  parameter dependence. It is well known that the first scenario offers favorable conditions 
  for Reduced Basis methods which have been well studied for this type of models and
  can therefore serve for comparisons. While in this case nonlinear schemes using neural networks
  do not seem to offer decisive advantages in terms of achievable certifiable estimation accuracy nor computational efficiency
  we see an advantage of the DNN approach in the second scenario because it seems that they can 
  be better adapted to the challenges of this problem class.
  
  It should be noted though that the present approach shares some conceptual constituents with so called One-Space methods or PBDW (Parametric Background Data Weak) methods (see \cite{CDDFMN,BCDDPW2,MPPY}). We therefore briefly recollect some related basic ideas in Section 
  \ref{ssec:1-space}. An important element is to represent the sensor functionals as elements of  the trial space
 $\U$ for the underlying PDE.  The $\U$-orthogonal projection to their span, termed ``measurement space'',
 provides a natural ``zero-order approximation'' to the observed state.
 To obtain an improved more accurate reconstruction, the data need to be ``lifted''
 to the complement space. We view the construction of such a ``lifting map'' as ``learning'' the
 expected ``label'' associated with a given observation, see Section \ref{sec:2}.
 This, in turn, is based on first projecting ``synthetic data'' in terms of parameter snapshots, to the $\U$-orthogonal complement of the measurement space.
 We then extract via SVD from these projected data a sufficiently
 accurate ``effective'' complement space that captures corresponding components of the solution manifold
 with high accuracy in $\U$. The lifting map is then expressed in terms of the coefficients
 of a $\U$-orthogonal basis of the effective complement space which, in turn, are represented
 by a neural network. The fact that we learn the coefficients of a $\U$-orthogonal basis 
 allows us to control the accuracy of the estimator in the problem-relevant norm.
 Thus, the proposed method is based on combining  POD based model reduction 
 with neural network regression. 
 In this regard, it shares similar concepts with the recent work in \cite{BHKS}
 and \cite{DDP}.  However, a major distinction is that in \cite{BHKS,DDP} the whole parameter-to-solution map is constructed while  the present work focuses on generating directly a {\em data-to-state} map through regression in 
 {\em sensor coordinates} (see Section \ref{ssec:2.3}). 
 
 Concerning (ii), successfully training a neural network remains a serious issue. Neither
 can one guarantee to actually exploit the expressive power of a given network architecture, nor
 does it seem possible to estimate the required computational cost. 
Therefore, it is more important to measure a given optimization outcome in a problem-relevant metric for making it possible to enable the assessment of the estimation quality obtained in the end. 
 In the second part of the paper we computationally
 explore a natural training strategy that, as will show, renders optimization more reliable, stable
 and robust with regard to varying algorithmic parameters like depth, width or learning rates.
 To counter gradient decay when increasing depth, we opt for ResNet architectures. 
 Moreover, we systematically compare stochastic gradient descent applied to all trainable parameters, termed ``plain \mwrevise{\RNN} training'' to an {\em expansion strategy} that starts
 with a shallow network and successively adds further blocks in combination with a {\em blockwise optimization}. That is, at every stage   we optimize only the trainables in a single block while freezing the remaining
 parameters. In that sense we do not fix any network architecture beforehand but expand it dynamically while monitoring loss-decay, see
 Subsection \ref{ssec:step3}. Section \ref{sec:num-res} is then devoted to extensive 
 numerical studies comparing plain \mwrevise{\RNN} and the expansion strategy for both application scenarios.
   
\section{Problem Formulation and Conceptual Preview}\label{sec:2}

\subsection{Parametric PDE models}\label{ssec:1.1}

For a wide scope
of applications, the underlying governing laws can be represented by a family of partial differential
equations (PDEs)
\be
\label{parfam}
\cF(u,\by;f)=0,\quad \by\in \cY,
\ee
with data $f$ and coefficients  depending 
on parameters $\by$ ranging over some compact parameter domain $\cY\subset \R^{d_y}$.
Focusing on linear problems, what matters in the present context is to   identify first
a {\em stable variational formulation} for \eqref{parfam}, i.e., to identify
a(n infinite-dimensional) {\em trial-space} $\U$ and a {\em test-space} $\V$,
such that
\be
\label{general_weak}
\cF(u,\by;f)(v)=0 = f(v)- (\cB_\by u)(v),\quad v\in \V,
\ee 
is well-posed, meaning that the bilinear form
\be
\label{b}
b(u,v;\by) = (\cB_\by u)(v),\quad u\in \U,\, v\in \V,
\ee
satisfies the conditions of the Babuska-Necas Theorem (continuity-, inf-sup-,
surjectivity-condition), \cite{LBB}. 

When $\by$ traverses $\cY$ the respective solutions $u(\by)$ mark the viable states of interest and form what is often   often called the {\em solution manifold}
\be
\label{sm}
\cM := \{ u(\by): \cF(u(\by),\by;f)=0,\quad \by\in\cY\}\subset \U,
\ee
which is in most relevant cases a compact subset of $\U$. The central task
considered in this work is to recover, from a moderate number of measurements, a state $u\in \U$ under the prior $u\in\cM$.

Stable formulations are available for a wide scope of problems to which the following
approach then applies. To be specific, we focus as a guding example second order
elliptic problems 
\be
\label{ellip}
\cF(u,\by;f)= f- \cB_\by u:=  f+\div (a(\cdot;\by)\nabla u),\quad \mbox{in }\, \Omega\subset \R^{d_x},
\quad u|_{\partial\Omega}=0,
\ee
where $f$ is a forcing term.
Examples are Darcy's equation for the pressure in ground-water flow or electron impedance tomography.
Both involve   second order elliptic equations as core models and 
{\em parameter dependent}   diffusion coefficients that  describe permeability or conductivity, respectively, where in the latter case the model is complemented
by Robin-type boundary conditions.
 A large parameter dimension $d_y$ reflects substantial  model complexity.
In a probabilistic framework\mw{,} the parametric representations of the coefficients
could  {arise}, for instance, from Karhunen-Lo\`{e}ve expansions
of a random field that   represent numerically  ``unresolvable'' features.  In this case the number of parameters is ideally
 {\em infinite} and parameter truncation causes additional model bias. 

For \eqref{ellip} a proper choice of trial- and test-spaces 
is $\U = \V= H^1_0(\Omega)$, provided that one has {\em uniform ellipticity},
namely there exist constants $0<r\le R<\infty$ such that 
\be
\label{ue}
r\le a(\cdot, \by)\le R,\quad \mbox{in }\, \Omega,\,\, \by \in \cY.
\ee
 For other problem types, such as 
convection dominated problems or time-space formulations for
parabolic problems, one may have to choose $\V$ differently from $\U$, see
e.g. \cite{DPW,DSW}.


Later in our numerical examples  we consider  
two types of diffusion coefficients: 

\begin{itemize}
\item[(S1)] {\em Piecewise constant coefficients:} Given $\{\Omega_j\}_{j\in I}$ a non-overlapping partition of $\Omega$, $\bigcup_{j\in I} \Omega_j = \Omega$, the diffusion coefficient in the model \eqref{ellip} is assumed to be of the following form: 
\begin{equation}
a(x; \by) = a_0(x) + \sum_{j\in I} y_j \chi_j(x),
 \label{eq:model_a}
\end{equation}
where $y_j$, the $j$-th component of $\by$, obeys a uniform distribution  and $\chi_j(x)$ is the characteristic function that takes the values $1$ in $\Omega_j$ and $0$ in its complement $\Omega\setminus\Omega_j$. This correponds to the affine parameter representation \eqref{eq:model_a}.
We assume for simplicity $a_0(x)= 1$ and $y_j \sim U[-1/2, 1/2]$ so that the diffusion coefficient is piecewise constant while the uniform ellipticity \eqref{ue} of the problem is guaranteed with a moderate condition of the variational formulation. 
\item[(S2)] {\em Lognormal case:} Suppose the diffusion coefficient has the following parameterized form: 
\begin{equation}
a(x; \by) = a_0(x) + a_1 e^{z(x, \by)},
 \label{eq:model_a_lognormal}
\end{equation}
where  $a_0$ is a continuous non-negative function on $\overline{\Omega}$, $a_1$ is a positive constant and $z(x, \by)$ is a zero-mean Gaussian random field 
\begin{equation*}
z(x, \by) = \sum_{j=1}^\infty \sqrt{\mu_j} \xi_j(x) \eta_j(\by), \quad x\in \Omega \text{ and } \by\in \cY.
\end{equation*}

Here, $\{\eta_j\}_{j\geq 1}$ form an orthogonal system over $\cY$ where the components of $\by$ are i.i.d. $\mathcal{N}(0,1)$, and the sequence $\{(\mu_j, \xi_j)\}_{j\geq 1}$ are real eigenpairs of the covariance integral operator  
\begin{equation*}
(Cv)(x) = \int_D c(x, x') v(x') dx'
 \label{eq:a_lognormal_ker}
\end{equation*}
associated with the Mat\'{e}rn model
\begin{equation*}
c(x,x') = \sigma^2 \frac{2^{1-\nu}}{\Gamma(\nu)} (2\sqrt{\nu}\, r)^{\nu} K_{\nu} (2\sqrt{\nu}\, r),
 \label{eq:matern}
\end{equation*}
in which  
$\sigma^2$ is the marginal variance, 
$\nu>\frac{1}{2}$ is the smoothness parameter of the random field, 
$\Gamma$ is the Gamma function, 
$r = \sqrt{\frac{(x_1-x_2')^2}{\lambda_{x_1}^2}+\frac{(x_2-x_2')^2}{\lambda_{x_2}^2}}$, 
$\lambda_{x_1}$ and $\lambda_{x_2}$ are the correlation lengths along $x_1$- and $x_2$-coordinates, 
and 
$K_{\nu}$ the modified Bessel function of second kind. 

At the discrete level, to generate realizations of the stationary Gaussian process over the grid points of $\Omega$, we use the circulant embedding approach in \cite{dietrich1997fast,kroese2015spatial}. 
Its main idea is to embed the covariance matrix to a block circulant matrix so that FFT can be applied for a fast evaluation. 
\end{itemize}

The reason for considering these two scenarios lies in the following principal distinctions: 
In scenario (S1) the choice of $\cY$ ensures that the diffusion parameter always satisfies \eqref{ue} so that 
\be\label{weak}
a(u,v;\by):= \int_\Omega a(\cdot;\by)\nabla u\cdot\nabla v \,dx = f(v),\quad v\in \U.
\ee
is stable over $\mathbb{U} = H_0^1(\Omega)=\V$  and possesses 
 for each $\by\in \cY$ a unique   weak solution $u=u(\by)\in \U$.
However, the condition of \eqref{weak} depends 
on $R/r$ and deteriorates when this quotient grows. 
In scenario (S2) this latter aspect is aggravated further 
because \eqref{ue} holds only with high probability. 
Moreover, in contrast to scenario (S1), the diffusion coefficient no longer is affine in $\by\in\cY$ which is known to poses challenges to the construction of certifiable
Reduced Bases, see e.g. \cite{RHP}. 


 
 \subsection{Sensors and Data}
 In addition to the model assumption that an observed state $u\in \U$ (nearly-)
 belongs to $\cM$ we wish to utilize   
   external information in terms of measurements or data.
Throughout this note we will assume that the number $m$ of measurements  
$\bo = (o_1,\ldots,o_m)^\top\in \R^m$ of an unknown state $u$ is fixed and of moderate size. We will always assume that the data 
are produced by {\em sensors} 
\be
\label{sensors}
o_i= \ell_i(u),\quad i=1,\ldots,m, \quad \bfell:= (\ell_1,\ldots,\ell_m)\in (\U')^m,
\ee
i.e.,   we assume, for simplicity, in what follows the $\ell_i$ to be {\em bounded linear functionals}.  

Estimating $u$ from
such data is ill-posed, already due to a possible severe under-sampling.
Under the assumption that $u$  (at least nearly) satisfies the PDE for {\em some} 
parameter $\by\in \cY$ one may also ask for such a parameter $\by$ that explains the
data best. Since the {\em parameter-to-solution map} $\by\to u(\by)$ is not 
necessarily injective, this latter parameter estimation problem is typically
even more severely ill-posed and nonlinear.

A common approach to regularizing both estimation tasks is 
{\em Bayesian inversion}. An alternative is to fix from the start
a presumably good enough discretization, say in terms of a large
finite element space $\U_h$ to then solve (a large scale optimization)
problem
\be
\label{opt}
\min_{\by\in\cY,u_h\in \U_h}\Big\{\|\bo - \bfell(u_h)\|^2+ \lambda\|\cF_h(u_h,\by;f)\|^2
+ \mu \|R_h(u_h)\|^2\Big\},
\ee
where $\|\cdot\|$ is just the Euclidean norm,  $\lambda$ and $\mu$ are weight parameters and the last summand represents a regularization term which is needed
since $d_y+{\rm dim}\,\U_h$ is typically much larger than the number of measurements $m$. Questions arising in this context are:
How to choose $\lambda, \mu$
and the regularization operator $R_h$; perhaps, more importantly, should one
measure deviation from measurements and the discrete residual - closeness to
the model - in the {\em same} metric?

Moreover, for each new data instance $\bo'$ one has
to solve the same large-scale problem again. Thus, employing a {\em reduced model} \cite{spacetime,YB}
for approximating the {\em parameter-to-solution map} $\by\mapsto u(\by)$ would
serve two purposes, namely mitigating under-sampling and speeding forward
simulations. 
Specifically, we apply similar concepts as in \cite{BCDDPW2,MPPY,CDDFMN,CDMN}
to base state estimation on {\em reduced modeling}, as described next.

\subsection{Sensor-Coordinates}\label{ssec:2.3}
Rather than tying a discretization directly to the estimation task, and hence
to a specific regularization, as in \eqref{opt}, reduced order modeling methods such as the Reduced Basis Methods (RBM) first prepare in an {\em offline phase} a reduced order model that requires the bulk of computation. It typically takes place in a ``truth-space''
$\U_h\subset \U$, usually a finite element space of sufficiently large dimension that is expected to comply with
envisaged estimation objectives (which could be even adjusted at a later stage). In particular, the reduced order model can be adapted to the solution manifold and the sensor system. As such it is only used indirectly in 
the recovery process and should be viewed as representing ``computing in $\U$''. 
A first important ingredient of the methods in \cite{BCDDPW2,MPPY,CDDFMN,CDMN}
is to {\em Riesz-lift} the functionals in $\bfell$ from $\U'$ to $\U$, thereby
subjecting them to the same metric as the states $u$. The obtained Riesz-representers $\phi_i\in \U$ of $\ell_i\in \U'$ then span an $m$-dimensional
subspace $\W\subset \U$, referred to as {\em mesurement space}. Thus, $\U$-orthogonal
projections $P_\W u$ of some $u\in \U$ to $\W$ encode the same data-information about
$u$ as $\bo=\bfell(u)$. This induces the decomposition
\be
\label{deco1}
\U= \W \bigoplus \W^\perp,\quad u = P_\W u + P_{\W^\perp} u,
\ee
so to speak representing any state $u$ in ``sensor coordinates'' $w=P_\W u
$, and ``labels'' $w^\perp= P_{\W^\perp} u$.

In these terms, recovering a state $u$ from its measurement $w=P_\W u$, means 
to approximate the label $P_{\W^\perp} u\in \W^\perp$. Thus, the state estimation
can be viewed as seeking a map 
\be
\label{map}
A : w\mapsto A(w),\quad A(w)= w +B(w),
\ee
where, in principle, $B:\W \to \W^\perp$ could be any map that hopefully  
exploits the fact that $w=P_\W u$ for some $u\in \cM$ in an effective way.

\subsection{Affine Recovery Map}
\label{ssec:1-space}
The methods in \cite{BCDDPW2,MPPY,CDDFMN} determine the {\em lifting-map}
$B: \W\to \W^\perp$ as a {\em linear} or {\em affine} map, termed as ``One-Space-Methods'', and in \cite{CDMN} as a {\em piecewise affine} map combining 
One-Space concepts with parameter domain decomposition and model selection.
As shown in \cite{CDDFMN}\mw{,} any affine map $B$ is characterized by an 
{\em affine} subspace
$\U_n = \bar u + \widetilde\U_n\subset \U$, 
where $\bar u$ is a suitable chosen offset state and  
$\widetilde\U_n$ is a linear space of dimension $n$, for which the estimator $A=A_{\U_n}$ satisfies 
\be\label{ustar}
u^*(w):=A_{\U_n}(w) =\argmin \{u\in w+\W^\perp : \|u- P_{\U_n}u\|_\U \} .
\ee
Moreover, $u^*(w)$ can be computed efficiently as a linear least-squares problem in $\U_n$ followed by a 
simple correction in $\W$. Uniqueness is ensured if
\be\label{mun}
\mu(\widetilde\U_n,\W) := 
\sup_{v\in \U_n}\frac{\|v\|_{\U}}{\|P_{\W}v\|_\U} <\infty,
\ee
which is the case if and only if $\U_n\cap \W^\perp = \{0\}$.
$\mu(\widetilde\U_n,\W)$ is actually computable as one over the smallest singular value of the cross-Gramian of an basis for $\U_n$ and $\W$. Hence, it has a geometric interpretation because it relates to the {\em angle} between the spaces $\widetilde\U_n$ and $\W$, tending to infinity when this angle approaches $\pi/2$.
This affects estimation accuracy directly since 
\be\label{err}
\sup_{u\in \cM}\|u^*(P_{\U_n}u) - u\|_\U 
\le \mu(\widetilde\U_n,\W)\e_n, \quad\mbox{when}\,\,   \dist\,(\cM,\U_n)_\U:= \sup_{u\in \cM}\inf_{z\in \U_n}\|u-z\|_\U\le \e_n.
\ee
This is actually best possible when using, as sole prior, the knowledge that $u\in \cM$
also belongs to the convex set $\cK({\U_n},\e_n):= \{u\in \U: \inf_{z\in \U_n}\|u-z\|_\U\le \e_n\}$, \cite{BCDDPW2}. Thus, overall estimation accuracy involves
a competition between the approximation property of the space $\U_n$ and
``visibility'' from $\W$. In fact, $\U_n$ can have at most dimension $n=m$
since otherwise $\mu(\widetilde\U_n,\W)=\infty$. Hence, the rigidity of an affine space severly limits estimation quality. In fact, it is shown in \cite{CDMN}
that restricting to such convex priors will generally fail to meet
natural estimation benchmarks which calls for employing {\em nonlinear
reduced models}. 

In this regard, the following comments will provide some orientation. One-space-methods appear to nevertheless perform very well for problems
where the solution manifold $\cM$ has rapidly decaying
 Kolmogorov $n$-widths
 \be\label{widths}
 d_n(\cM)_\U := \inf_{{\rm dim}\V_n\le n}\dist(\cM,\V_n)_\U.
 \ee
For elliptic problems \eqref{ellip} this is known to be the case, even in 
high parameter dimensional regimes, 
when \eqref{ue} holds and the diffusion coefficients depend {\em affinely}
on the parameters, as is the case in scenarion (S1), \eqref{eq:model_a}.
Approximations from judiciously chosen linear spaces can then be very effective.
Moreover, affine parameter-dependence is also instrumental in the {\em greedy}
construction methods of Reduced Bases for rate-optimal reduced linear models
which play a central role in one-space-methods, 
\cite{RHP,BCDDPW,DPW,CDDN,CDMN}. It is therefore interesting to see
how the nonlinear estimators proposed in the present paper compare with
affine recovery methods in scenarios where the latter are known to perform well.
For this reason we include scenario (S1) in our numerical tests.

That said, the effectivity of one-space-methods depends in a rather sensitive
way on the above favorable conditions, namely rapid decay of $n$-widths, uniform
ellipticity, and
affine parameter dependence. We therefore include scenario (S2) with log-normal
parameter dependence. In fact, affine dependence no longer holds, the behavior of $n$-widths is much less clear,
and the coefficient field may nearly degenerate, challenging the validity 
of \eqref{ue}. 

\subsection{A Regression Framework}\label{ssec:regr}

To avoid monitoring pointwise errors in high-dimensions
we opt now for a mean-square accuracy quantification which responds in a less
sensitive way to
high parameter dimensionality and blends naturally into a learning context.
A natural probabilistic model could be  based on a  probability measure
$\mu$  
 on $\U$ with support on $\cM$. $u= (w=P_\W u,w^\perp= P_{\W^\perp} u)$, 
 is then viewed as a random variable. The optimal estimator would then be
 the nonlinear map 
 \be
 \label{optest}
A^*(w) = w+ B^*(w),
\ee
where the conditional expectation $B^*(w) = \int_{w+\W^\perp} P_{\W^\perp}u\, d\mu(u|w)$  is 
the {\em regression function} minimizing 
\be
\label{obj}
\int_\U \|u - A(P_\W u)\|^2_\U d\mu(u) = \int_\W \int_{w+\W^\perp} \|w^\perp -
B(w)\|_\U^2d\mu(u|w) d\mu(w)
\ee 
over all mappings of the form $A(w)= w+ B(w)$, $B:\W\to \W^\perp$.

The central goal in what follows is to construct numerical estimators that
approximate $B^*$ well. Note that this approximation should take place in $\U$
for the estimator to respect the natural problem metrics.

\section{State Estimation Algorithm}\label{sec:3}
\subsection{Computational Setting}
We describe next in more details how to set up a learning problem in the problem compliant norm $\|\cdot\|_\U$.
We adhere to the sensor-induced 
decomposition \eqref{deco1}. Of course, the underlying Riesz-lifts of the measurement functionals $\bfell$ cannot be computed exactly but need to be approximated within controllable accuracy.
In the spirit of Reduced Basis methodology, we employ a sufficiently large ``truth space'' $\U_h\subset \U$, which we choose here as a conforming finite element 
 space. The scheme is based on the following steps.\\

\noindent 
{\bf (1) The Measurement Space:}
As a major part of the offline stage, one then solves
the $m$ Galerkin problems 
\be
\label{Rieszh}
(\tilde\phi_h^i,v_h)_{\U}=\ell_i(v_h), \quad \forall v_h\in \U_h ,\text{ and } i=1,\ldots, m,
\ee 
providing (approximate)
{\em Riesz representers} $\tilde\phi_h^i\in \U_h$ of the linear functionals $\ell_i$. 
Then define 
\be\label{Wh}
\W_h:= \verb|span| \{\phi_h^1, \ldots, \phi_h^m\}
\ee 
where the $\{\phi_h^1, \ldots, \phi_h^m\}$ result from orthonormalizing the lifted  functionals  $\{\tilde\phi_h^1, \ldots, \tilde\phi_h^m\}$. 
Hence, 
\be\label{orthproj}
P_{\W_h} u:= \sum_{i=1}^m (u,\phi_h^i)_\U\phi^i_h,
\ee
is the orthogonal projector from $\U_h$ onto 
$\W_h$ that encodes the information provided by
the sensor system as elements of the trial space $\U_h$.

Specifically, given an observation vector $\bfell(u)\in \R^m$ of some observed state $u\in \U$, we can determine $P_{\mathbb{W}_h}(u)$ as follows. Let $\tilde\bPhi$, $\bPhi \in \U_h^m$ denote the column vectors,
obtained by lining up the functions $\tilde \phi_h^i, \phi_h^i$, $i=1,\ldots,m$. Moreover let 
$\mathbf{C}\in \R^{m\times m}$ denote the
(lower triangular) matrix that realizes the change of bases $\mathbf{C}\tilde\Phi= \Phi$. One readily checks that 
\be\label{proja}
\mathbf{w}:=  \mathbf{C}\ell(u_h) =\big((u_h,\phi^1_h)_\U,\ldots,(u_h,\phi^m_h)_\U\big)^\top,
\ee
i.e., the $\U$-orthogonal projection of $u$ to the measurement space $\W_h$ is given by
\be\label{by}
P_{\W_h} u = \sum_{i=1}^m w_i\phi^i_h =: 
\bPhi^\top \mathbf{w}, 
\ee
providing the best approximation to an observed space from $\W_h$.
This can therefore be viewed
already as a ``zero-order'' reconstruction of
the observed space in $\W_h\subset\U$.

Note that the problems \eqref{Rieszh} are always {\em elliptic} Galerkin 
problems, regardless of the nature of \eqref{parfam}. The Riesz-lift (as a mapping from 
$\U'$ to $\U$) has condition number equal to one.
 Therefore, these preparatory offline calculations are stable and 
often a posteriori error estimates allow one assess the accuracy of Riesz-lifts, controlled 
by the choice of the truth space $\U_h\subset \U$. This should not be confused with the accuracy of $P_{\W_h}u$ with respect to $u$ which is primarily limited by the number and type of sensors.

So far, we have ignored noise of the data
$\bo=\bfell(u)$. Noise in the observation vector $\bo$ carries over to noise in the
coefficients $\bw$ from \eqref{proja} possibly inflated by the condition of the transformation matrix $\mathbf{C}$, see \cite{CDMN} for a more detailed discussion of this issue.

 
Our subsequent numerical experiments refer to the model problem \eqref{ellip} where the
sensors are given by the local averages
 of $u$ over subdomains $B_i\subset \Omega$
\be\label{elli}
\ell_i\left(u\left(x, \by\right)\right) = \frac{1}{|B_i|} \int_{B_i} u(x, \by)\, dx,
\ee
with local neighborhoods $B_i$, and where  the diffusion coefficients are
piecewise constants on a $4\times 4$ checkerboard partition of $\Omega$, see
\eqref{eq:model_a}. Thus the parametric dimension is $16$.

\begin{figure}[h!]
\centering
\begin{subfigure}{.45\textwidth}
  \centering
  \includegraphics[width=\linewidth]{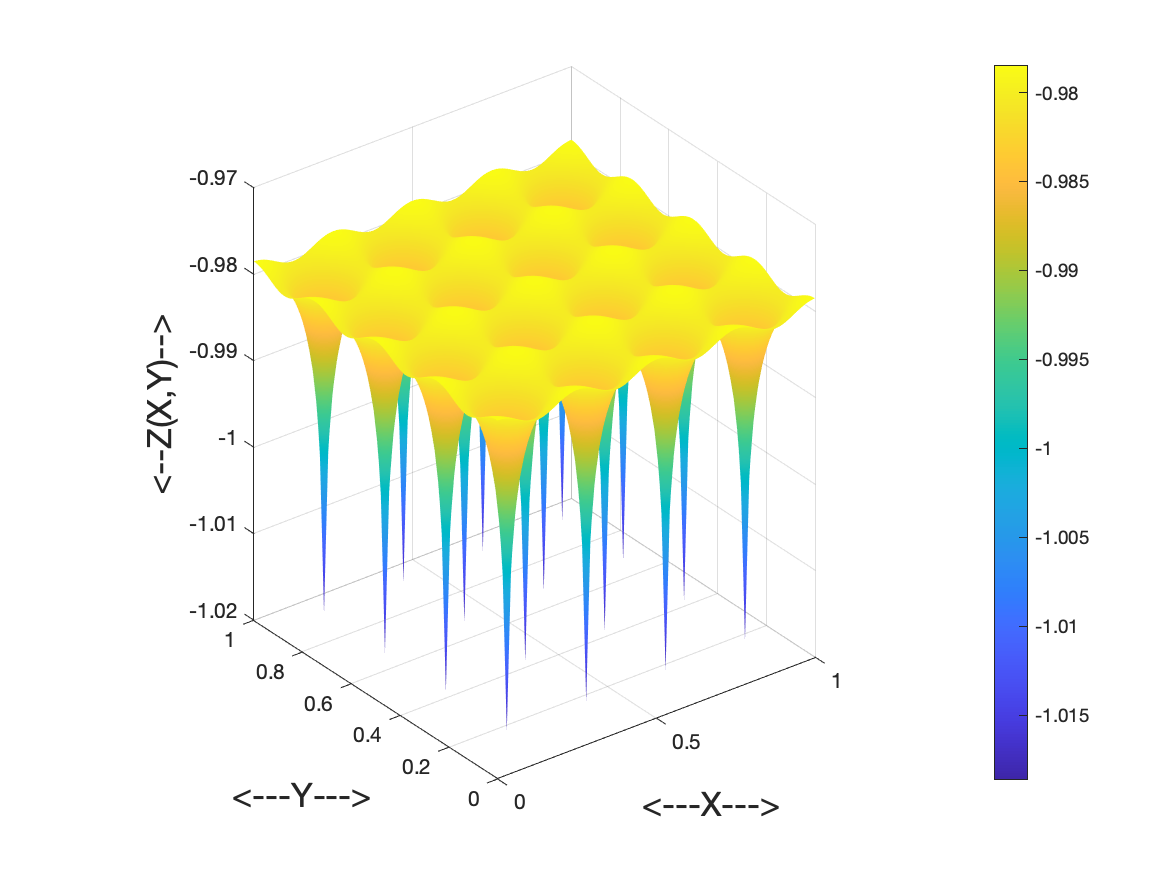}
   \caption{$\phi_h^1$}
\end{subfigure}\hfil
\begin{subfigure}{.45\textwidth}
  \centering
    \includegraphics[width=\textwidth]{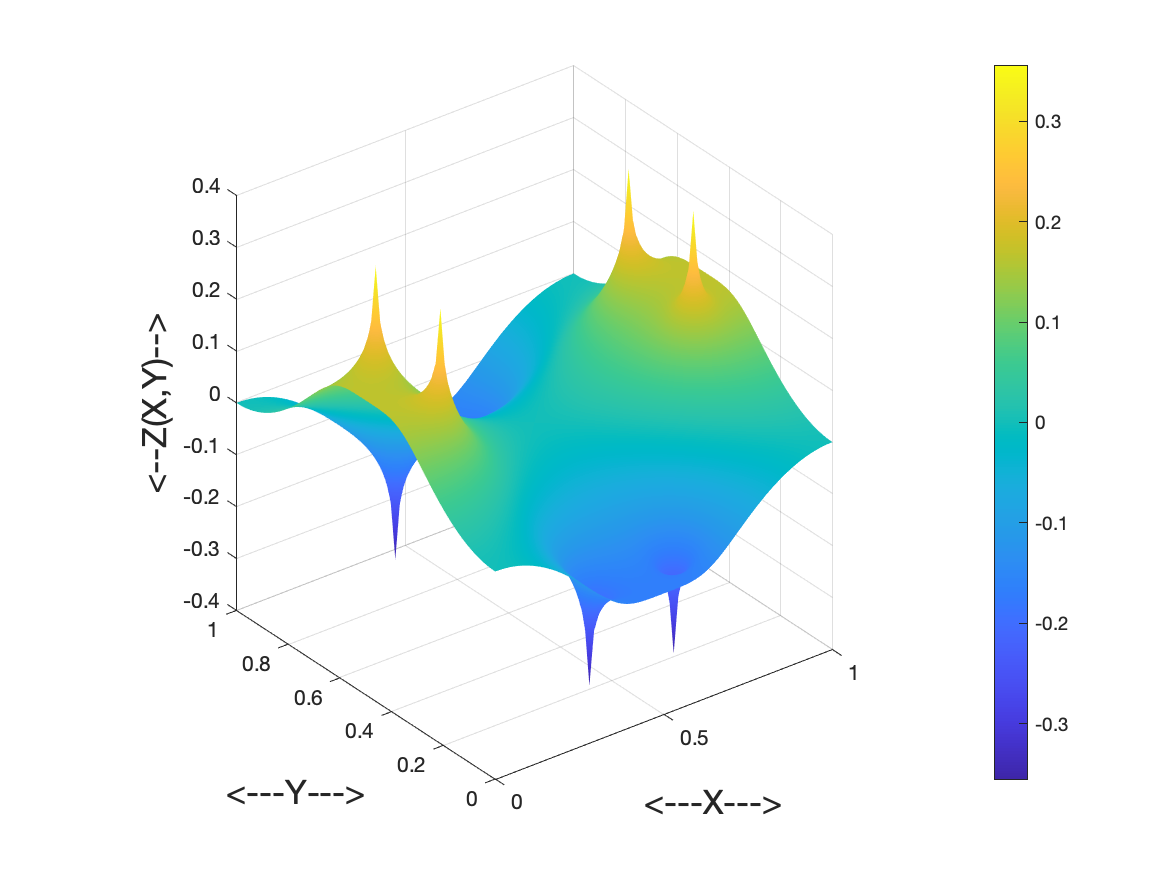}
    \caption{$\phi_h^6$}
\end{subfigure}
\begin{subfigure}{.45\textwidth}
  \centering
  \includegraphics[width=\linewidth]{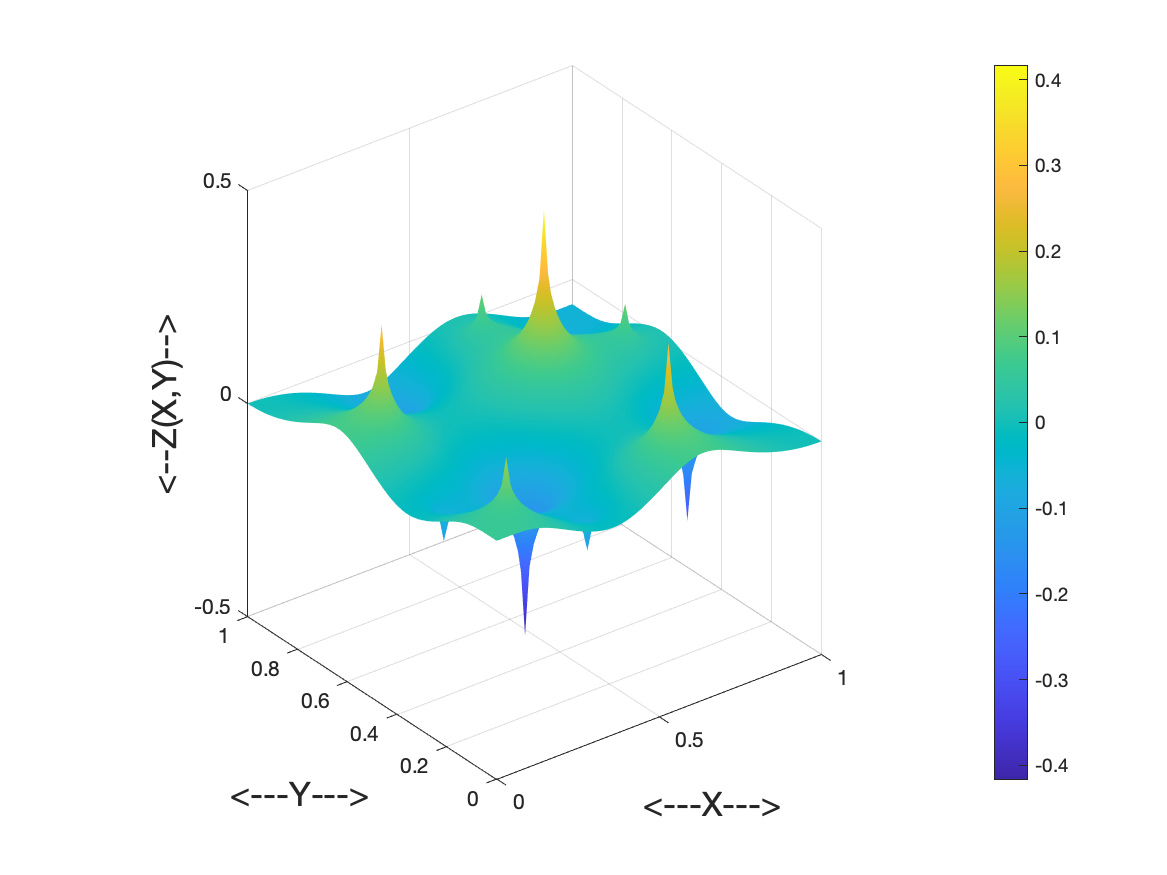}
   \caption{$\phi_h^{11}$}
\end{subfigure}\hfil
\begin{subfigure}{.45\textwidth}
  \centering
    \includegraphics[width=\textwidth]{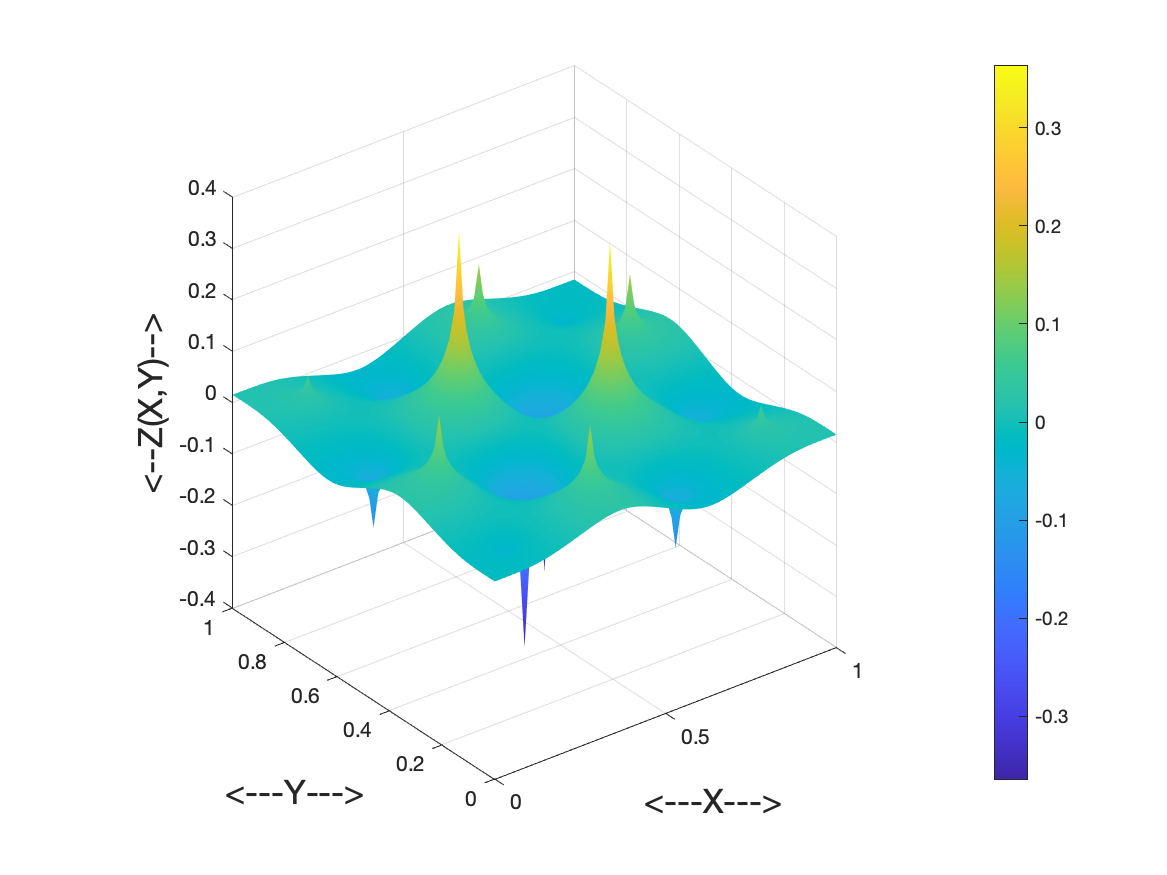}
     \caption{$\phi_h^{16}$}
\end{subfigure}
\caption{Example basis functions of the measurement space $\mathbb{W}_h$.}
\label{fig:sensorsReisz}
\end{figure}

Specifically, in subsequent experiments we realize those 
by taking   the average value of the FE solutions at the four points of a small square $B_i$ enclosing the sensor location. Denoting these four points by  $\{x_{p_j}\}_{j=1}^4$, we have
\be\label{four}
\ell_i(u_h(\by)) = \frac{1}{4} \sum_{j=1}^4 u_h (x_{p_j}, \by),
\ee
entering the right hand side of \eqref{Rieszh}.
{A numerical illustration of selected normalized basis, $\phi_h^1, \phi_h^6, \phi_h^{11}, \phi_h^{16}$, is given in Figure \ref{fig:sensorsReisz}, where the normalized basis vectors are obtained by the SVD of the matrix formed by the finite element coefficient vectors of the lifted functionals with respect to the inner product weighted by $(v_h, v_h)_{\U}$.} 

\medskip
\noindent
{\bf (2) Generating  Synthetic Data:}
We randomly pick  a set of parameter samples $\by^s$, for $s=1, \ldots, \widehat{N}$, and and employ a standard solver to compute the FE solutions $u_h^s=u_h(x; \by^s) \in \U_h$ as ``snapshots'' at the selected parameter values. These FE solutions provide the high-fidelity ``truth" data to be later used for training the estimator towards minimizing the regression objective \eqref{obj}. 
Corresponding synthetic measurements then take the form
\be
\label{train}
\bw^j= (w_{h,1}^j,\ldots,w_{h,m}^j)^\top,\quad w_h^j:= P_{\W_h} u_h^j
= \sum_{i=1}^m w_{h,i}^t \phi_h^i, \quad j=1,\ldots,\widehat N.
\ee
as detailed by \eqref{proja} and \eqref{by}.

In addition we need the ``training labels'' providing complement information
\be
\label{compl}
z_h^s= (I-P_{\mathbb{W}_h})u_h^s,\quad s=1,\ldots, \widehat{N}.
\ee


\medskip
\noindent
{\bf (3) Approximate $\W_h^\perp$:} To extract the essential information provided 
by $z_h^s$, 
we perform next
a singular value decomposition to the resulting point-cloud of finite-element coefficient vectors
${\bf z_h^s}\in \R^{N_h}$. Suppose we envisage an overall  estimation target tolerance $\eta>0$.  We then retain only
those $k$ left singular vectors corresponding to singular values larger than or equal to a value 
$\tilde\eta$ which is typically less than $\eta$, for the following  reason. 
First, incidentally, the SVD and the decay of singular values indicate whether the complement information of $\cM$ can be adequately captured by a linear space of acceptable size within some target tolerance.
While the $H^1$-norms of the $z^s_h$ are uniformly bounded, the coefficients $\bz^s_h$ in their respective
finite element representation convey accuracy only in $L_2$, and so does the truncation of singular values. Strictly speaking, employing standard inverse inequalities, one should take 
$\tilde\eta \le h\eta$, where $h$ is the mesh-size in $\U_h$.

We next $\U$-orthonormalize the finite element functions corresponding to the retained left singular  
vectors, arriving at a $\U$-orthonormal basis $\bPsi= \{\psi_h^1, \ldots, \psi_h^k\}\subset\W_h^\perp$.
We then  define 
\be
\label{Wperp}
\Wp := \verb|span|\{\psi_h^1, \ldots, \psi_h^k\}
\ee
as our {\em effective complement space} that accommodates the training labels $c_h^s$ given by
\be
\label{labels}
c^{{s}}_h(x):= P_{\mathbb{W}_h^\perp} u_h^{{s}}= \sum_{i=1}^k c_i^{{s}} \psi_h^i =\bPsi^\top\bc^s,\quad c_i^{{s}}:= (u_h^{{s}},\psi_h^{{s}})_\U,\quad {{s}}=1,\ldots,\widehat N.
\ee
Note that $c_h^s$ is essentially a compression of $z_h^s$.

\begin{figure}[h!]
\begin{subfigure}{.45\textwidth}
  \centering
  \includegraphics[width=\textwidth]{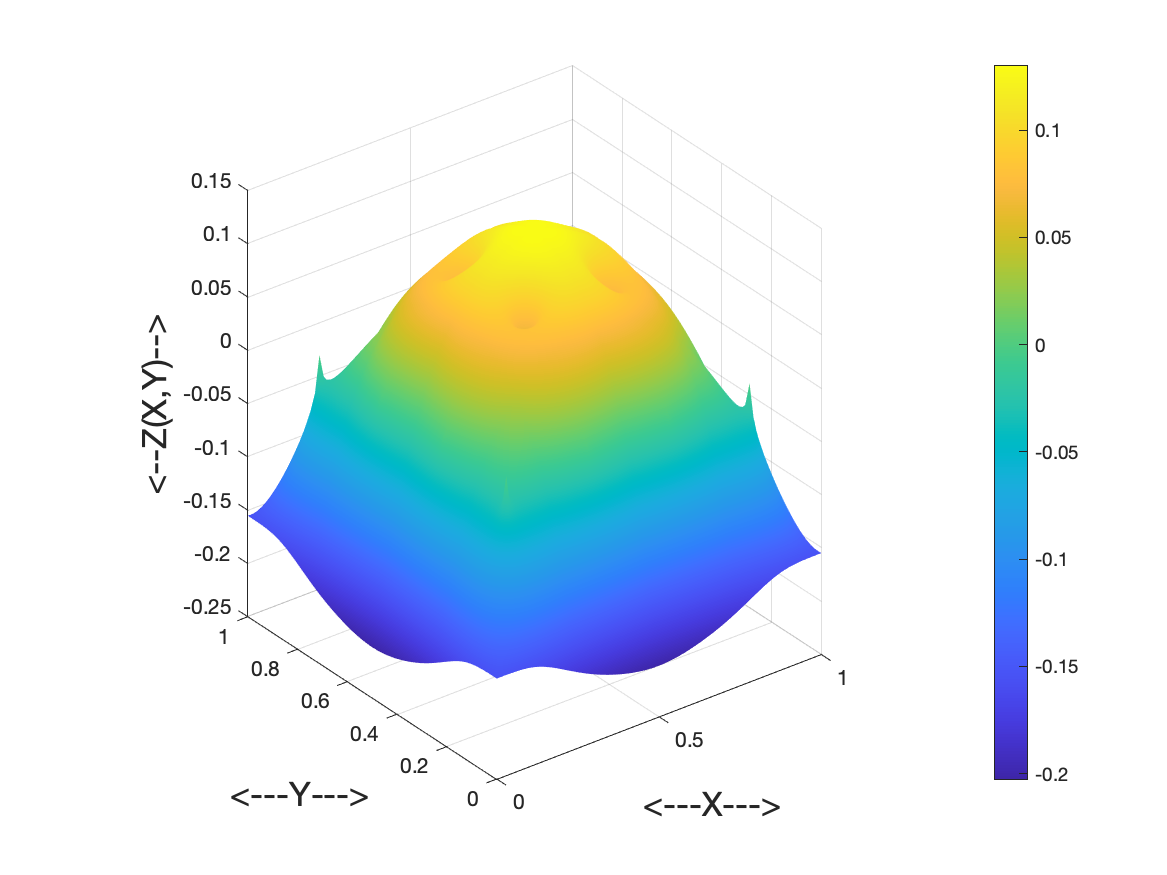}
  \caption{$\psi_h^1$}
\end{subfigure}\hfil
\begin{subfigure}{.45\textwidth}
  \centering
    \includegraphics[width=\textwidth]{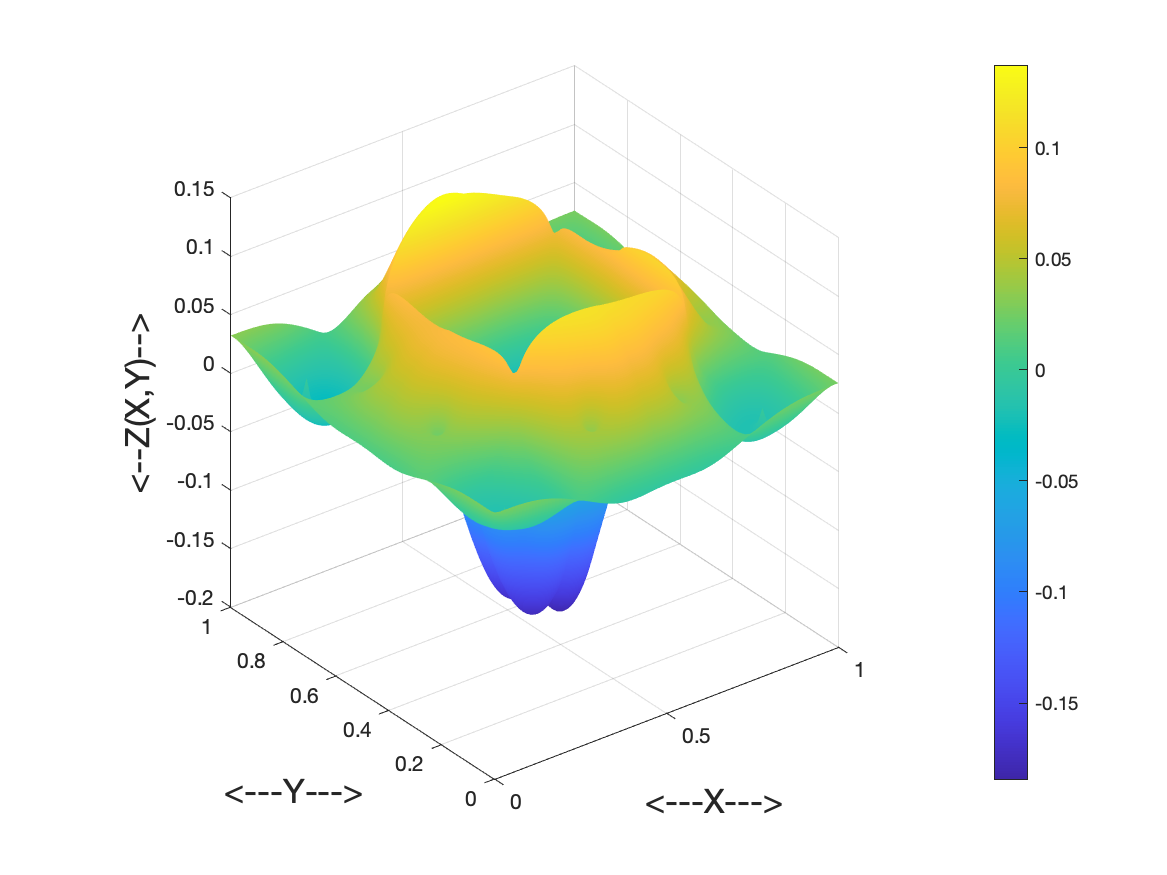}
      \caption{$\psi_h^6$}
\end{subfigure}
\begin{subfigure}{.45\textwidth}
  \centering
  \includegraphics[width=\textwidth]{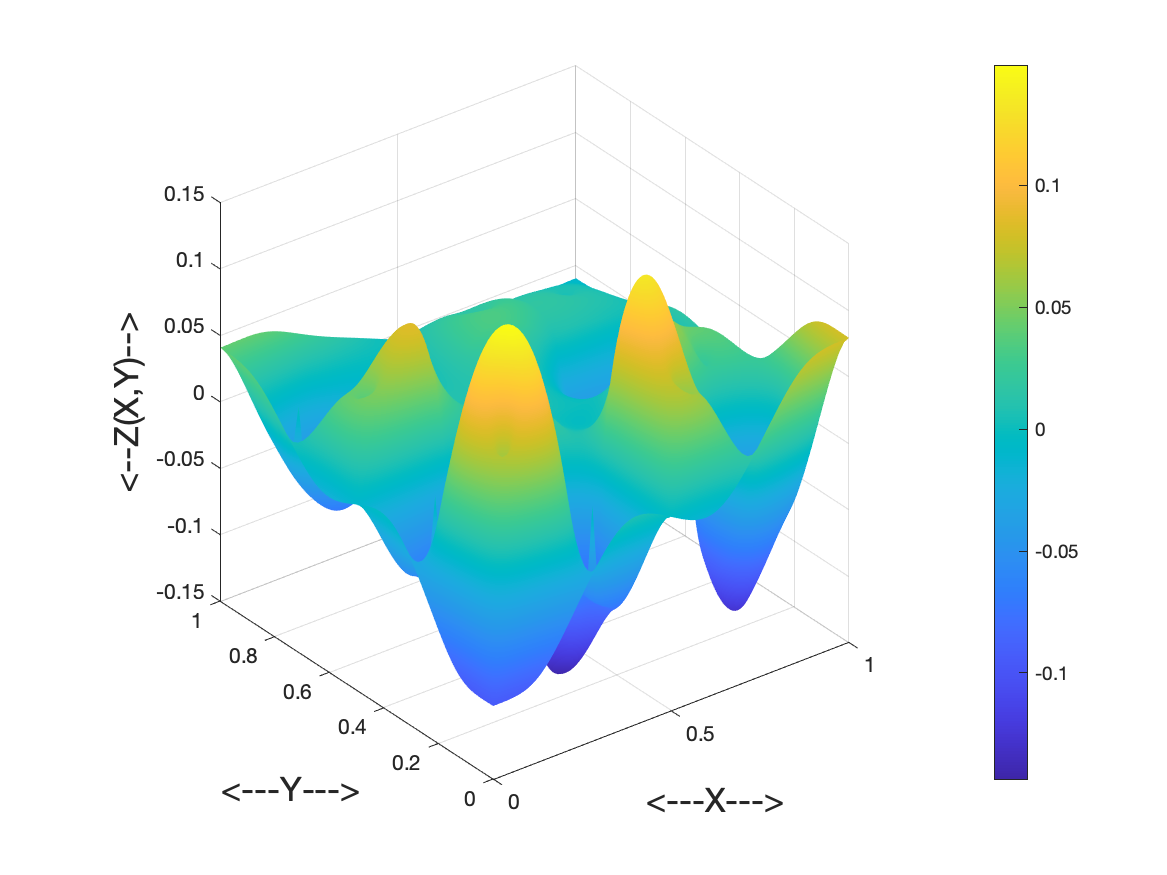}
    \caption{$\psi_h^{11}$}
\end{subfigure}\hfil
\begin{subfigure}{.45\textwidth}
  \centering
    \includegraphics[width=\textwidth]{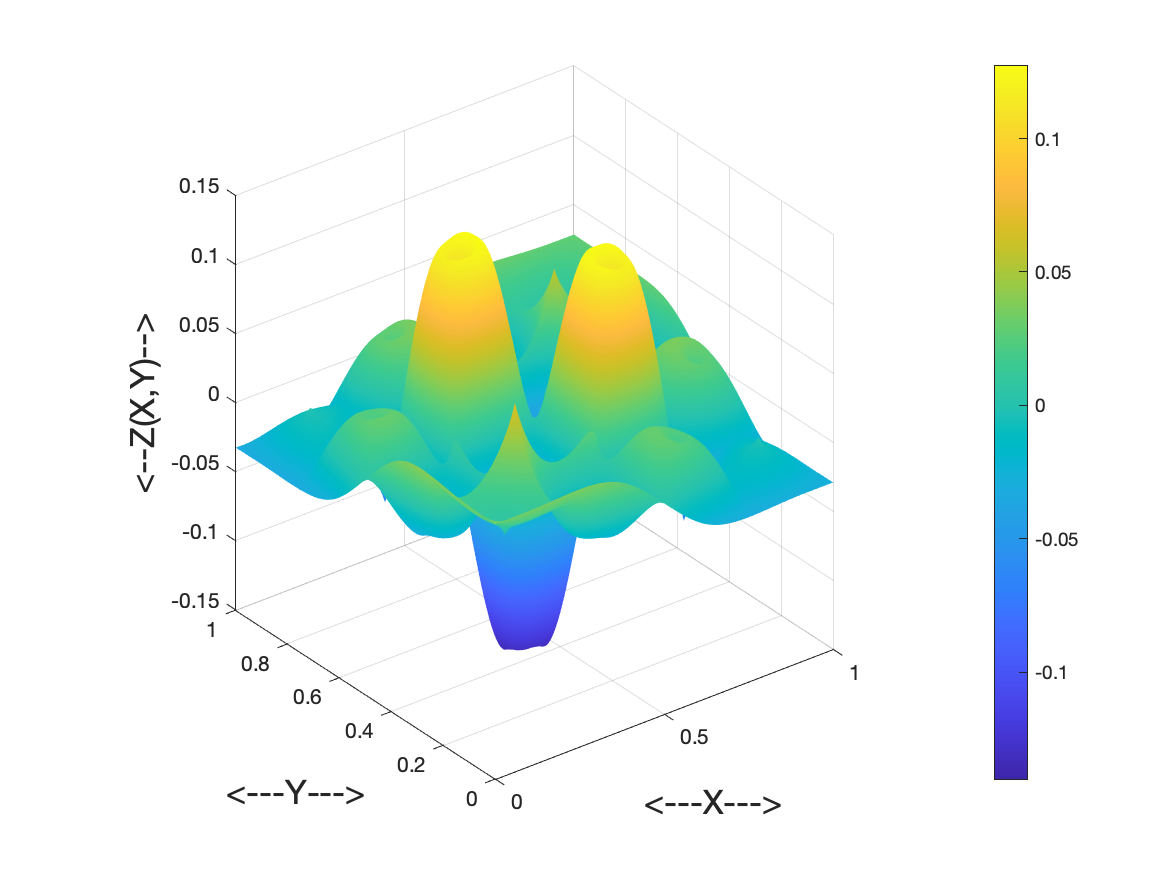}
      \caption{$\psi_h^{16}$}
\end{subfigure}
\caption{Example basis functions of the complement space $\Wp$.}
\label{fig:psi}
\end{figure}

In brief, after projecting the snapshot data into the complement of the measurement space,  we determine a set of orthogonal basis $\psi_h^{j}$ that optimally approximates the data in the sense that
$$
  \min_{\{\psi_h^1,\ldots, \psi_h^k\}} \frac{1}{\widehat{N}} \sum_{s=1}^{\widehat{N}}
  \left\| z_h^s - 
  \sum_{j=1}^k \left( z_h^s, \psi_h^j \right)_{\mathbb{U}} \, \psi_h^j 
  \right\|_{\U}^2 
$$
subject to the conditions $(\psi_h^j, \psi_h^i)_{\U} = \delta_{i,j}, \ 1 \leq i, j \leq k$, where $\delta_{i,j}$ is the Kronecker delta function. A numerical illustration of the normalized basis of $\Wp$, $\psi_h^1, \psi_h^6, \psi_h^{11}, \psi_h^{16}$, is shown in Figure \ref{fig:psi}.

\begin{rem}
\label{rem:accurate}
\end{rem}
Approximation in $\U$ is now realized by training
coefficient vectors in the Euclidean norm of $\R^k$ because
\be\label{norms}
\Big\|\sum_{j=1}^k a_j \psi_h^j\Big\|^2_\U = \sum_{j=1}^k a_j^2 =: \|\ba\|_2^2,\quad \ba \in \R^k.
\ee 
In other words, the estimator respects the intrinsic problem metrics which is 
a major difference from common approaches involving neural networks.

\medskip

As argued earlier, we wish to construct a nonlinear map that
recovers from data $w\in \W_h$ a state $w+ B(w)\in w+\Wp$  for an appropriate $B: \W_h\to \Wp$.
More specifically, for any $w_h=\sum_{i=1}^m w_i \phi_h^i=:{\bw^\top}\bPhi\in \W_h$, 
the envisaged mapping $B$ has the form 
\be\label{B}
B(w_h)= \bw_\perp(\bw)^\top \bPsi,
\ee
where 
\be\label{wperp}
\bw_\perp(\bw) = (w^1_{\perp,h}(\bw),\ldots, w^k_{\perp,h}(\bw))^\top\in \R^k, 
\ee
and
\be\label{Ncoeff}
\bw_{\perp}(\bw)=
\NN(\bw)
\ee
will be represented as a neural network with input data $\bw\in \R^m$ and output dimension $k$.
 \\


\noindent  {\bf  (4) Loss Function and Training the Neural Network (NN):} 
We randomly select a subset $\cU_{train}=\{u_h^1, \ldots, u_h^{N_{train}}\}$ from the ``truth" FE data $\cU_{sample}= \{u^{{s}}_h: {s}=1,\ldots,\widehat{N}\}$, generated in step (1), for training purpose and leave the rest $\cU_{ {ghost}}:=\cU_{sample}\setminus\cU_{train}$ for testing. The {associated} coefficient vectors  $\{\bw^1, \ldots, \bw^{N_{train}}\}$ and $\{\bc^1, \ldots, \bc^{N_{train}}\}$ are used as the training data for minimizing 
the natural empirical {loss} analogous to \eqref{obj}  
\be\label{loss}
\verb|Loss| =  R(\Theta):= \frac{1}{N_{train}} \sum_{t=1}^{N_{train}} \| \bc^t- \NN(\bw^t;\Theta) \|_2^2,
\ee
where we assume in what follows that the $\NN$ depends on a collection $\Theta$ of {\em trainable} or {\em hyper-parameters}, over which $\verb|Loss|$ is to be minimized.

Once the training of $\NN$ is completed, for any new measurement  $\bfell\left(u\left(\cdot, \by\right)\right)$
with coefficient vector $\bw\in \R^m$, defined by \eqref{proja}, \eqref{by}, the observed state $u(\cdot, \by)$    is approximated by 
 	$$
 	u \approx  
\bPhi_h^\top \bw + \bPsi_h^\top \NN(\bw).
$$

A natural question would be now to analyze the performance of the estimator
obtained by minimizing the loss. This could be approached by employing
standard machine learning concepts like Rademacher complexity in conjunction
with more specific assumptions on the network structure and the solution manifold. We postpone this task
to forthcoming work and address in the remainder of this paper instead a more elementary issue, namely 
 ``optimization success'' which is most essential for a potential merit of the
proposed schemes.
%

\subsection{A Closer Look at Step {(4)}}\label{ssec:step3}

\begin{figure}[!h]
\centering
\includegraphics[width=1\textwidth]{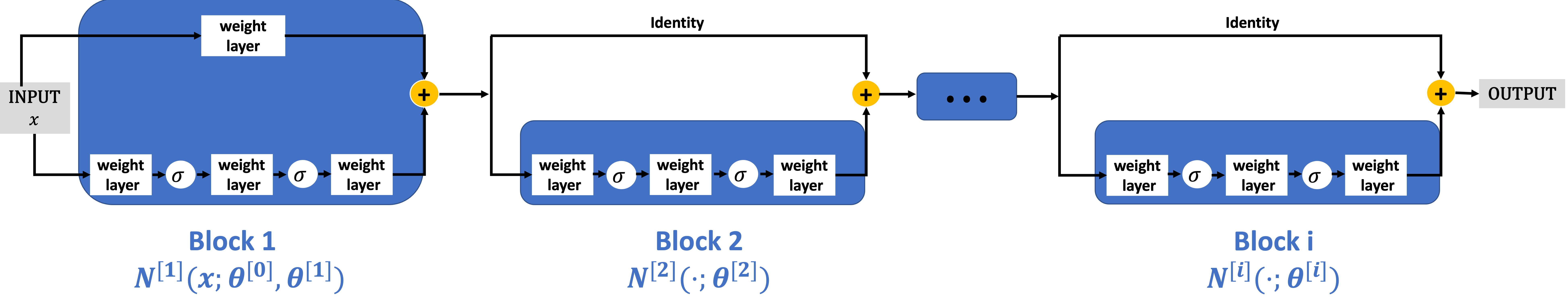}
\caption{ResNet with $i$ blocks}
\label{resnet_structure}
\end{figure}

It is a priori unclear which specific network architecture and which budget
of trainable parameters is appropriate for the given recovery problem.
Even if such structural knowledge were available, it is not clear
whether the expressive potential of a given network class can be exhausted
by the available standard optimization tools largely relying on stochastic gradient descent (SGD) concepts. Aside from formulating a model-compliant regression problem, our second primary goal is to explore numerical strategies
that, on the one hand, render optimization stable and  less sensitive on algorithmic settings whose most favorable choice is  usually not known in 
practice. On the other hand, we wish to incorporate and test some simple
mechanisms to {\em adapt} the network architecture, inspired by classical 
``nested iteration'' concepts in numerical analysis. Corresponding simple
mechanisms can be summarized as follows: (a) {A} ResNet architecture with its skip-connections is known to mitigate ``gradient damping'' which impedes the adjustment
of trainable parameters in lower blocks when using deep networks. (b) Each block
in a \mwrevise{ResNet} can be viewed as a perturbation of the identity and may therefore be
expected to support a stable incremental accuracy upgrade. (c) Instead of
``plain training'' of an \mwrevise{ResNet}, where \mwrevise{a stochastic gradient descent} is applied to all trainable parameters
simultaneously, we study an iterative \mwrevise{training} strategy that optimizes only single
blocks at a time while freezing trainable parameters in all other blocks.
Part of the underlying rationale is that training a shallow network is more reliable and efficient than training a deep network. We discuss these issues next in more detail. 
\vspace*{-4mm} 
\paragraph{Neural Network Setting:}
Specifically, we consider the following {\em Residual Neural Network} (\mwrevise{ResNet}) structure with $i$ blocks:
\begin{equation}\label{architecture}
\NN^{[i]}(\bx;\theta^{[0]},\theta^{[1]},\cdots,\theta^{[i]})  := \cN^{[i]}\big( \NN^{[i-1]}(\bx;\theta^{[0]},\theta^{[1]},\cdots,\theta^{[i-1]});\theta^{[i]}\big).
\end{equation}
with each block defined by 
\begin{equation*}
\begin{split}
\cN^{[1]}(\bx;\theta^{[0]},\theta^{[1]}) &:= W_3^{[1]}\sigma\big(W_2^{[1]}\sigma(W_1^{[1]}\bx+b_1^{[1]})+b_2^{[1]} \big)+ W^{[0]}\bx,\\
\cN^{[i]}(\bx, \theta^{[i]})  &:= W_{3}^{[i]}\sigma\big(W_2^{[i]}\sigma(W_1^{[i]}\bx+b_1^{[i]})+b_2^{[i]} \big)+ \bx, \qquad i =2,3,\cdots
\end{split}
  \end{equation*}
where ${\theta^{[i]}} = \{W_1^{[i]}, W_2^{[i]},W_3^{[i]}, b_1^{[i]}, b _2^{[i]}\}$ for $i =1,2, \cdots$ and $\theta^{[0]} = \{W^{[0]}\}$,  $\sigma(\cdot)$ is some pointwise nonlinear function\mwrevise{. In the numerical experiments of Section \ref{sec:num-res}, we specifically take the activation function $\sigma$ to be $ \tanh$. See} Figure \ref{resnet_structure} for an descriptive diagram of the \mwrevise{ResNet} structure. 

\paragraph{An Expansion Strategy:} As indicated earlier, we will not approximate the recovery map using a fixed network architecture. We rather start with training a single hidden layer until saturation. We then proceed approximating residual data with another 
shallow network, and iterate this process which results in a ResNet structure. The basic rationale is similar to ``full multigrid'' or ``nested iteration'' in numerical analysis.

To describe the procedure in more detail, let
$$
\cD := \{\bw^i: i=1,\ldots, N^{train}\}\times \{\bc^i: i=1,\ldots,N^{train}\}
$$
denote the set of training data,  as described before. In view of \eqref{norms}, 
we employ simple mean-squared loss functions. In these terms we may rewrite \eqref{loss} as
\begin{equation}
\mathcal{L}_\cD [f]: =  \frac{1}{\#\cD}\sum_{(\bw,\bc)\in \cD} \| \bc- f({\bw}) \|_2^2.
\end{equation}

In these terms, the procedure can be described as follows:
At the initial step $f$ is taken as a shallow network $\NN^{[1]}(\cdot; \theta^{[0]}, \theta^{[1]})$
and solve first  the optimization problem:  
\begin{equation*}
\begin{split}
\textbf{OP-1 :}\qquad &\min_{(\theta^{[0]}, \theta^{[1]})\in \Theta^{[1]}}  \qquad\mathcal{L}_\cD[\NN^{[1]}(\cdot; \theta^{[0]}, \theta^{[1]})],
\end{split}
\end{equation*}
where $\Theta^{[1]}$ denotes the budget of hyper-parameters used at this first stage, encoding
in particular, the widths. \textbf{OP-1} is treated \mwrevise{with SGD-based optimizers} until the loss stagnates or reaches a (local) minimum. The resulting hyper-parameters are denoted 
by $(\theta_*^{[0]}, \theta_*^{[1]})$.

In a next step, we expand the neural network by introducing a new block and solve
\begin{equation*}
    \begin{split}
\textbf{OP-2 :}\qquad &\min_{(\theta^{[0]}, \theta^{[1]}, \theta^{[2]})\in \Theta^{[2]}}  \quad\mathcal{L}_\cD[\NN^{[2]}(\cdot;\theta^{[0]}, \theta^{[1]},\theta^{[2]})],
\end{split}
\end{equation*}
over an extended global budget of hyper-parameters $\Theta^{[2]}$.
\noindent Moreover,   \textbf{OP-2} itself is solved by sweeping over sub-problems 
as follows. Fixing the values of the trained parameters $\theta^{[1]}= \theta^{[1]}_*$ and 
$\theta^{[0]}=\theta^{[0]}_*$ we apply gradient descent first only over the newly added trainable parameters $\theta^{[2]}$. That is\mw{,} we solve, 
\begin{equation*}
   \textbf{OP-2$^*$ :}\qquad \min_{ \theta^{[2]}}  \quad\mathcal{L}_\cD[\NN^{[2]}(\cdot;\theta_*^{[0]}, \theta_*^{[1]},\theta^{[2]})].
\end{equation*}
Due to the \mwrevise{ResNet} structure, upon defining the data
$$
\cD_1: =  \{\big(\NN^{[1]}(\bw;\theta_*^{[0]},\theta_*^{[1]}) ,\bc -\NN^{[1]}(\bw;\theta_*^{[0]},\theta_*^{[1]})\big)\ |\ (\bw,\bc)\in \cD\},
$$
this is equivalent to optimizing 
\begin{equation*}
    \min_{\theta^{[2]}} \quad \mathcal{L}_{\cD_1}[\NN^{[2]}(\cdot;\theta^{[2]})]
\end{equation*}
 with a $\NN^{[1]}-$ projected input and a residual from \textbf{OP-1} as output. 
 
 We can thus  expect $\NN^{[2]}$ to perform better than $\NN^{[1]}$ in terms of learning the map between the observation $\bw$ and state-labels $\bc$ because adding an approximate residual presumably
 increases the overall accuracy of the estimator. 
 
 Once a new block has been added, the parameters in preceding blocks presumably
 could be further adjusted. This can be done by repeating a sweep over
 all blocks, i.e., successively optimizing each block while freezing the hyper-parameters in all remaining blocks.
 
 The retained ghost-samples in $\cU_{ghost}$ can now be used to assess the generalization error.
 If this is found unsatisfactory, the current \mwrevise{ResNet} can be expanded by a further block, leading to analogous optimization tasks \textbf{OP-k}. Moreover, for $k>2$, the advantages of training shallow networks can be exploited by performing analogous block-optimization steps freezing the hyper-parameters in all but one block. 
 In addition, we are free to design various more elaborate sweeping strategies \mwrevise{to cover all trainable parameters need to be updated}. For instance,
one could \mwrevise{add a round of updates toward} all trainable parameters at the end of the training process. We refer to \S  \ref{num:sweep} for related experiments.

Another way of viewing this process is to consider the infinite optimization problem 
\begin{equation*}\label{res_inf_prob}
\textbf{OP:}\qquad \min_{\theta\in   \mathbb{R}^{\infty}}  \mathcal{L}_D[\NN(\cdot; \theta)].
\end{equation*}
where $\NN(\cdot;\theta) : = \lim_{i}\NN^{[i]}(\cdot;\theta^{[0]}, \theta^{[1]},\cdots, \theta^{[i]})$ is an idealized  \mwrevise{ResNet} with an infinite number of blocks. Then the above block expanding scheme can be interpreted as a greedy algorithm for approximately solving the infinite problem   \textbf{OP} 
\begin{equation*}
\begin{split}
\textbf{OP-i }=\qquad &\min_{\theta^{[0]},  \theta^{[1]}, \cdots, \theta^{[i]}}  \mathcal{L}_D[\NN(\cdot; \theta^{[0]},  \theta^{[1]}, \cdots, \theta^{[i]}, \bzero)].
\end{split}
\end{equation*}\label{OPs}
successively increasing $i$ until the test on ghost-samples falls below a set tolerance.

In Section \ref{sec:RB_compare}, we will first present the accuracy of the recovered solutions with the \textbf{nonlinear} map learned with a \mwrevise{ResNet} compared to solutions obtained with the Reduced Basis Method (RBM). 
We showed that in piece-wise constant case, \mwrevise{ResNet} can reach similar accuracy as RBM. 
In log-normal case, while the RBM can not be applied to obtain reasonable solutions, we showed that \mwrevise{ResNet} can still be used to obtain a nonlinear solution with $L_2$ lifting (see Section \ref{sec:log_normal}).

\section{Numerical Results}\label{sec:num-res}

The above description of a state estimation algorithm is so far merely a skeleton.
Concrete realizations require fixing concrete algorithmic ingredients such as batch sizes, learning rates, and the total number of training steps (whose meaning will be precisely explained later in this section). \mwrevise{The numerical experiments
reported in this section have two major purposes: (I) shed some light on  the 
dependence of optimization success on specific algorithmic settings, in particular,
regarding two different principal training strategies.
The first one represents standard procedures and applies Stochastic Gradient Descent (SGD)
variants {\em to all trainable parameters}, defining a given DNN with ResNet architecture.
Schemes of this type differ only by various algorithmic specifications listed below and will
be referred to as \RNN ~and ``training'' then refers to the corresponding {\em global} optimization. The second strategy does not aim at optimizing a fixed network
but intertwines optimization with a {\em blockwise network expansion}.
This means that initially only a shallow network is trained which is subsequently expanded
in a stepwise manner by additional blocks in a ResNet architecture. Such schemes are
referred to by \eRNN. The corresponding training, referred to as {\em blockwise training},
applies SGD only to the currently newly
added block, freezing all parameters in preceding blocks. Comparisons between \RNN ~and \eRNN
~ concern achievable generalization error accuracies, stability, robustness with regard to algorithmic settings, and efficiency.}
(II) compare the performance of neural networks with estimators that are based on Reduced Basis concepts for Scenario (S1) of piecewise constant affine parameter representations of the diffusion coefficients, introduced in \S \ref{ssec:1.1}. 
 (S1) is known to be a very
favorable scenario for Reduced Basis methods
that have been well studied for this kind of problems exhibiting excellent performance, \cite{CDDFMN,CDMN}.

Furthermore, we explore the performance of the \mwrevise{ResNet}-based estimator for scenario (S2), involving log-normal random diffusion parameters. In this case the performance of Reduced Basis is much less understood and hard
to certify. Corresponding experiments are discussed in Section \ref{sec:log_normal}.
To ease the description of the experiment configurations, we list the notations in Table \ref{tab:abbreviation}.





\newcommand{\hyp}{{\text -}}
\begin{table}[!h]\footnotesize
\hspace*{-1.2cm}
	\centering
	\begin{tabular}{|c|c|p{6cm}| p{7cm}| }
		\hline
				&\textbf{Abbreviation}&	\textbf{Explanation}& \textbf{Comments} \\
				\hline
	 &\mwrevise{ResNet} & Residual neural network & \mwrevise{General architecture is defined as in \eqref{architecture}. Activation function $\sigma = \tanh$.} \\
		\hline
			\multirow{2}{2cm}{\mwrevise{Training Scheme}} &\mwrevise{\RNN}& \mwrevise{A deep network with ResNet architecture is trained by standard SGD methods applied in a global fashion to the full collection
of trainable parameters.} & \\
		 \cline{2-4}
		  &\mwrevise{\eRNN}& \mwrevise{Deep networks with ResNet architecture are generated and trained by successively expanding the current network by a new block, confining SGD updates to a single block at a time} & \\
		 \hline
\multirow{3}{2cm}{\mwrevise{ResNet Configurations}}		 	&$\mathsf	{B}$&	Number of blocks  &\\
				\cline{2-4}
&$\mathsf	{W}$ &	Width of hidden layers 
&Random parameter initialization drawn from a  Gaussian distribution
\\
		\cline{2-4}
&	$\mathsf	{O}$ & Output dimension of $\bc$ & Number of orthogonal basis taken in $\Wp$\\
		\hline
	\multirow{4}{2cm}{\mwrevise{Training Hyper-parameters}}	&	$\mathsf	{l}$ &	$l^1$ regularization weight of trainable learning parameters & $0.00001$ if not specified\\
		\cline{2-4}
			&$\mathsf	{b}$ &	Batch size & $100$ if not specified\\
 		\cline{2-4}
 		&$\mathsf	{lr}$ &	Learning rate &\\
 		 		\cline{2-4}

 		&$\mwrevise{\mathsf{T}}$ &	\mwrevise{Number of total training steps}  &\\
		\hline
			\multirow{6}{2cm}{\mwrevise{Data Type}}&$\mathtt{Train}$& Training data type&\\
		\cline{2-4}
		&\textit{pwc} & Piecewise constant case training set (S1)&  $10,000$ samples subject to $16$ uniformly distributed sensors
		\\
	\cline{2-4}
		&\textit{log-normal}& Log-normal case training data (S2) & $1,000$ samples when $16$ sensors; $5,000$ samples when $49$ sensors\\
		\cline{2-4}
		&$\mathtt{sen}$& Number of uniformly distributed sensors & \\
	\cline{2-4}
		&POD-$H^1$& Dimension reduction of $\Wp$ with POD in $H^1$ sense&\\
		\cline{2-4}
		&POD-$L_2$& Dimension reduction of $\Wp$ with POD in $L_2$ sense&\\
		\hline
	\end{tabular}
	\caption{Abbreviations}
	\label{tab:abbreviation}
\end{table}


\subsection{Numerical Set Up}
\mwrevise{ResNets defined as in \eqref{architecture} with activation functions $\tanh$ are used for all experiments.  These networks will \mwrevise{be   optimized} with the \mwrevise{aid of the} Proximal Adagrad \mwrevise{algorithm. This  is} essentially a variant of   SGD which is capable of 
\mwrevise{adapting learning rates} per parameter. 
We further provide \mwrevise{in Section \ref{optimizer} a comparison between Adagrad with the Adam algorithm  to justify this choice.}
}

\mwrevise{
To ensure all experiments results are fair comparisons, for each set of experiments, we will fix the {number of} total training steps ($\mathsf{T}$). \mwrevise{``Training step''  means one iteration in Adagrad algorithm, which we sometimes also refer to  as {\it{one update step}}} of the trainable parameters. In particular, if the ResNet is trained in a block-wise sense (\eRNN), 
then each block will be attributed an equal number of $\mathsf{T}/\mathsf{B}$ steps to update the trainable parameters in this block. The choice of $\mathsf{T}$ will be specified in the loss history figures as well as in the error tables for each set of experiments.}

\mwrevise{With \eRNN~ it is, in principle, possible to apply in the course of the 
training process SGD repeatedly to blocks that had
been added at an earlier stage of the expansion process.
In most numerical examples, only one round of parameter updating will be applied to the
newly added last block. 
These blocks will not be revisited at a later time.
The only exception takes place in Section \ref{num:sweep}, where global updates {are carried out} in addition to the block-wise training for the ResNet. 
Here we wish to see whether \eRNN~ provides favorable initial guesses for a subsequent \RNN.
We often refer to any arrangement of the order of blockwise updates and revisiting blocks as \emph{training schedule}}. 
 
\mwrevise{For measurement of the errors, w}e wish to estimate the relative analogue to the ideal regression risk \eqref{obj}:

$${\mathcal{E}:= }\left(\frac{\cE_{\by\in\cY_{test}} ||u(\by)-u_{pred}(\by)||_{\U}^2}{\cE_{\by\in\cY_{test}}||u(\by)||_{\U}^2}\right)^{\frac{1}{2}},
$$
where $u_{pred} := \bPhi_h \bw + \bPsi_h \bc_{{pred}} $.
In particular, if the norm is taken to be the problem compliant norm $||\cdot||_{\U}$, due to Remark \ref{rem:accurate}, evaluating its empirical counter part:
$$ 
{\hat{\mathcal{E}} =}\left( \frac{\sum_{s }||\bc^s - \bc^s_{{pred}}||^2_2}{\sum_s||\bw||_2^2+||\bc^s||^2_2}\right)^{\frac{1}{2}}.$$
just requires computing Euclidean norms for the predicted coefficients $\{\bc_{pred}^{s}\}$. Of course, we expect sufficient large sample sizes provide accurate estimates
$$\hat{\mathcal{E}} \approx \mathcal{E}.$$

If one would like to evaluate  the accuracy of the recovered solution $u_{pred}(\by)$ in a norm that is different form the natural norm ($||\cdot||_{\U}$), quadrature in the truth space $\U_h$ is required  to approximately evaluate the respective norm of the functions.

\subsection{The Piecewise Constant Case (S1)}\label{sec:pc}

We first consider the aforementioned diffusion problem \eqref{ellip}  with $f=1$ and a piecewise constant diffusive parameter within $\Omega = [0, 1]^2$. More specifically, we consider a non-overlapping, $4\times 4$ uniform decomposition $\{\Omega_j\}_{j=1}^{16}$ of $\Omega$. The 
diffusion coefficient is a constant on each $\Omega_j$ as defined in \eqref{eq:model_a} .

\subsubsection{\mwrevise{\RNN} vs. \mwrevise{\eRNN}}\label{sec:PWC_eRNN}

\begin{figure}[h!]
\centering
	\centering
	\includegraphics[width=0.5\textwidth]{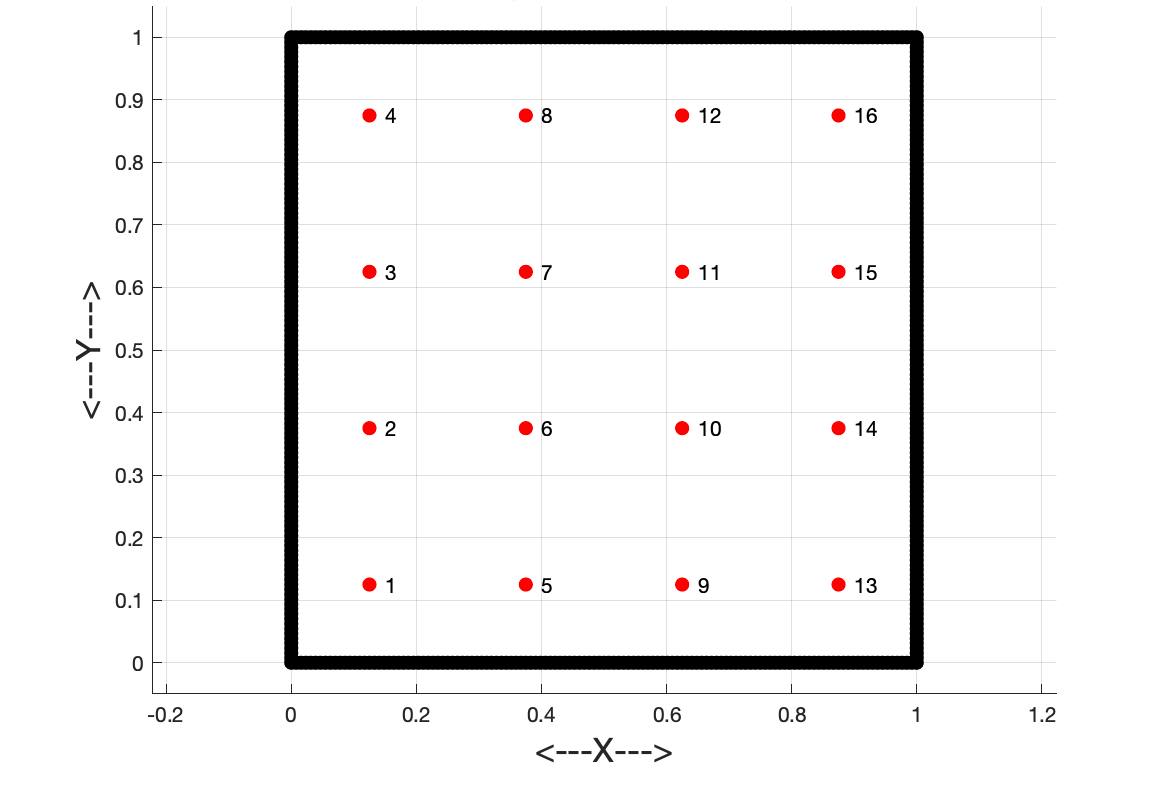}
	\caption{16 uniformly distributed sensors.}
\label{fig:sensors}
\end{figure}

The first group of experiments concerns a general performance comparison between 
the ResNet expansion strategy - in short \mwrevise{\eRNN}, outlined above and \mwrevise{a global update strategy which updates} the whole network with the same architecture, termed \mwrevise{\RNN} in what follows.  By ``performance'' we mean
training efficiency as well as corresponding achieved training and generalization losses.
Specifically, we consider first the case where $m=16$ sensors are placed uniformly in $\Omega$ (see Figure \ref{fig:sensors}) and the corresponding measurements are evaluated by averaging the solution at the four vertices of a square of side $0.001$ centered at the sensor location, see \eqref{four}. Thus, the number of sensors equals in this case the parametric dimension so that there is a chance that the measurements determine the state uniquely.

We use in total $10,000$ snapshots represented in the truth-space that serve as synthetic data, $500$ of which are used for testing purposes. To draw 
theoretical conclusions a larger amount of test data would be necessary. However,
intense testing has revealed that larger test sizes have no significant 
effect on the results in the scenarios under consideration. Based on computations, using these data we find that $k= 28$  basis functions suffice to reserve $99.5\%$ of the $H^1$-energy in $\Wp$, represented by the Hilbert-Schmidt norm of a full orthonormal basis in $\Wp$. The case studies documented by subsequent figures are referenced as follows:
  ``$pwc$'' refers to ``piecewise constant diffusion coefficients'' in scenario (S1);
  ``POD-$H^1$'' indicates that the sensors have been Riesz-lifted to $H^1(\Omega)$ which accommodates the measurement space $\W$.  {The SVD truncation threshold   is chosen to ensure accuracy in $H^1(\Omega)$; recall also that} ``sen16'' means that the recovery is based on data from $16$ sensors. 

The numerical results,  shown in Figure \ref{fig:pwc-ernn_s16}, indicate that the training loss resulting from \mwrevise{\eRNN} decays faster than the standard \mwrevise{\RNN} training in this case. The  relative generalization errors for both approaches after $12\times 10^5$ steps of training are displayed in  Table \ref{tab:w200_pwc},  where the generalization errors are evaluated on the test set of $500$ samples. The result shows that \mwrevise{\eRNN} outperforms \mwrevise{\RNN} in  a sense detailed later below. 
In particular, the training loss for
\mwrevise{\eRNN} drops   faster to a saturation level which can be achieved by \mwrevise{\RNN} 
only at the expense of a significantly larger training effort.
We also compare the estimated generalization error in $H^1$ with 
the achieved (expectedly smaller)  {error} in the weaker $L_2$-norm. 

One observes though that increasing network depth, i.e., employing a higher number of ResNet blocks
does not increase accuracy significantly in either model. This indicates that a moderate level of nonlinearity suffices in this scenario. This is not surprising considering the moderate number
of POD basis functions needed to accurately capture complement information. 
\newcommand{\fige}[7]{$\mathtt{Train}$-\textit{#1}-$\mathsf{sen}#2$-POD-$#3$, #4-$\mathsf{W}#5$-$\mathsf{O}#6$-$\mathsf{lr}#7$}
\newcommand{\figB}[7]{$\mathtt{Train}$-\textit{#1}-$\mathsf{sen}#2$-POD-$#3$, $\mathsf{B}#4$-$\mathsf{W}#5$-$\mathsf{O}#6$-$\mathsf{lr}#7$}
\newcommand{\figshort}[6]{$\mathtt{Train}$-\textit{#1}-$\mathsf{sen}#2$-POD-$#3$, $\mathsf{W}#4$-$\mathsf{O}#5$-$\mathsf{lr}#6$}

\begin{figure}[h!]
    \centering
    \includegraphics[width = 0.5\textwidth]{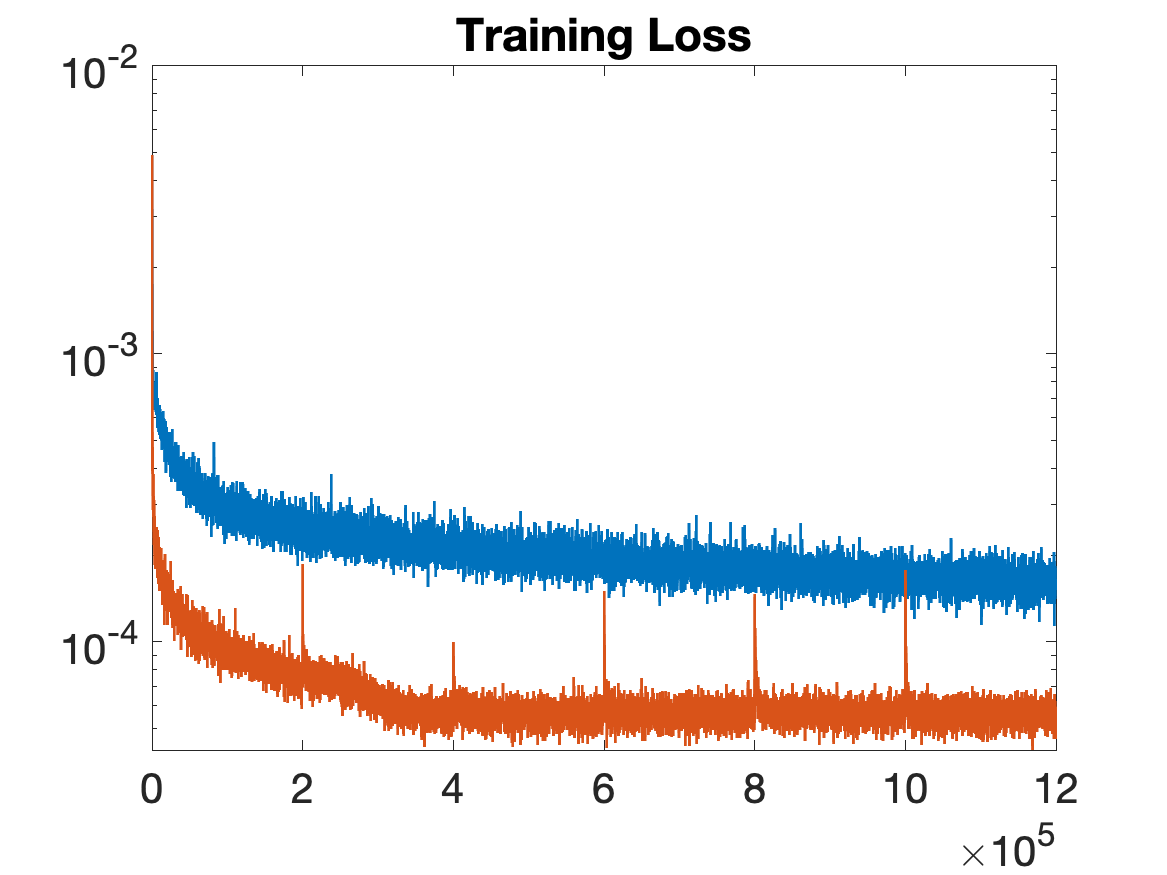}
    \caption{\textbf{\mwrevise{\eRNN(red)}/\mwrevise{\RNN(blue)} training loss comparison}. \figB{pwc}{16}{H^1}{6}{200}{28}{0.03}. 
    }
    \label{fig:pwc-ernn_s16}
\end{figure}

\begin{figure}[!htp]
	\centering
	\begin{subfigure}[b]{0.47\textwidth}
		\centering
		\includegraphics[width=\textwidth]{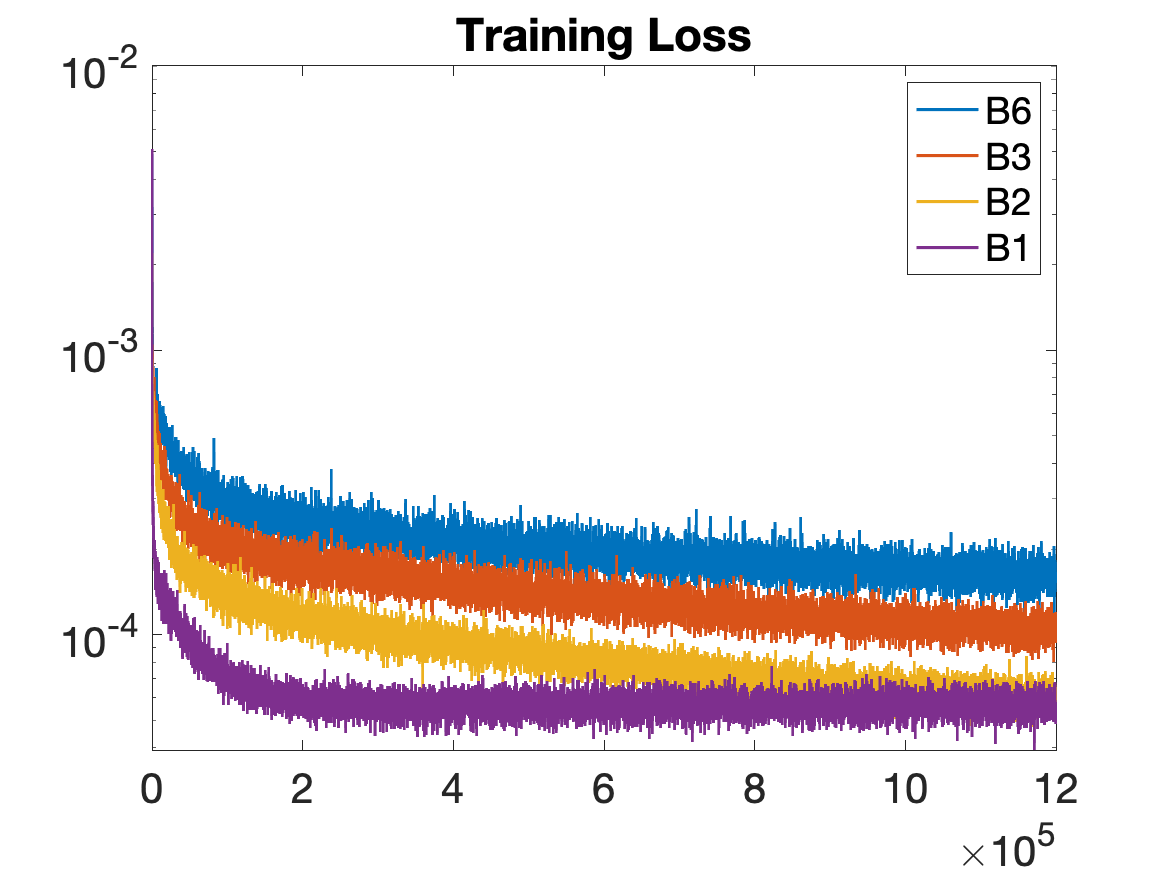}
		\caption{\mwrevise{\textbf{\RNN}} training loss}
		\label{RNN_train_s16_a}
	\end{subfigure}
	\hfill
	\begin{subfigure}[b]{0.47\textwidth}
		\centering
		\includegraphics[width=\textwidth]{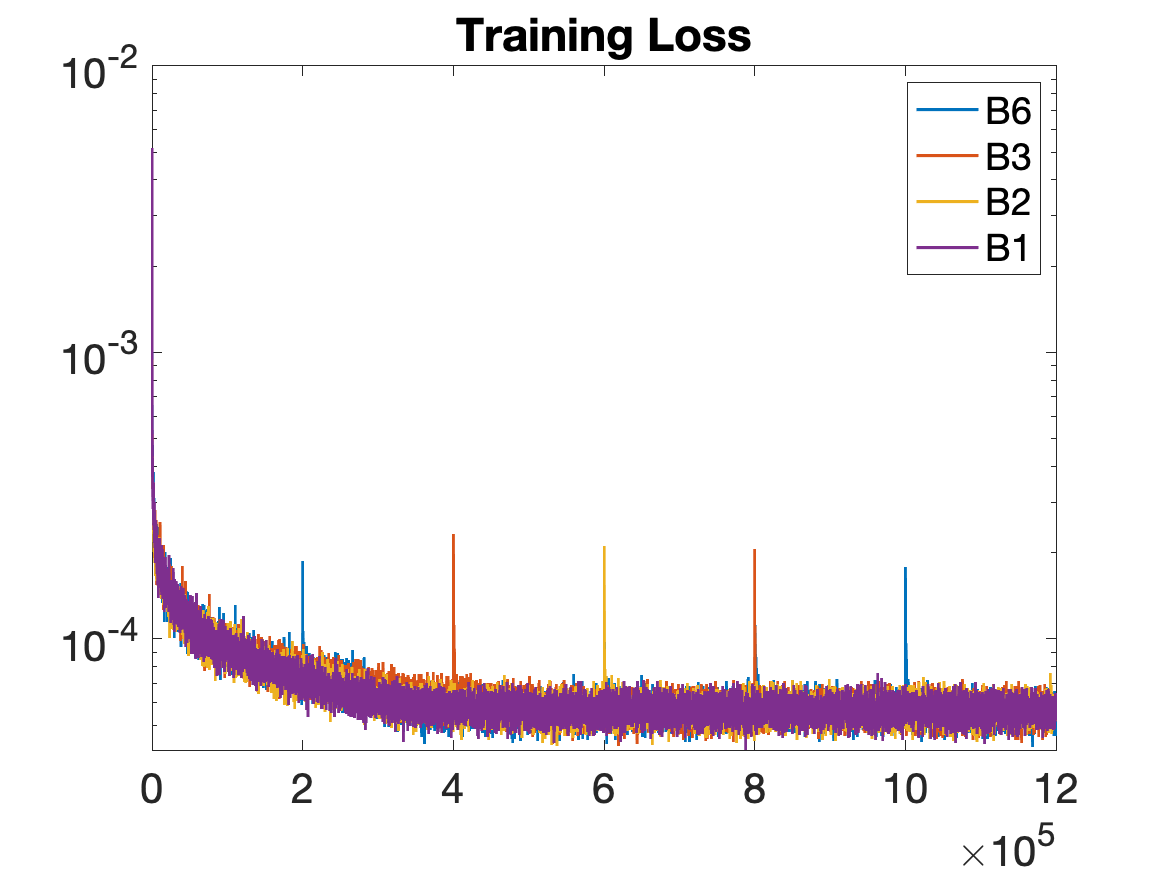}
		\caption{\mwrevise{\textbf{\eRNN}} training loss}
		\label{eRNN_train_s16_b}
	\end{subfigure}
	\caption{\textbf{Training loss of \mwrevise{\RNN}/\mwrevise{\eRNN} v.s. different number of Blocks}. \figshort{pwc}{16}{H^1}{200}{28}{0.03}.}
\end{figure}

\begin{table}[!h]
	\centering
	\begin{tabular}{c|c|c c |cc }
	 	&&\multicolumn{2}{|c}{$\hat{\mE}$}&\multicolumn{2}{|c}{Relative $L_2$ Error of $u_{pred}$}\\
	 	\hline
		\# of blocks&\#of trainable &\mwrevise{\eRNN} & \mwrevise{\RNN}&\mwrevise{\eRNN} & \mwrevise{\RNN}\\
 		\hline
         1 & 49,648 &9.76\%&9.76\% &3.01\%&3.01\%\\
         2 & 101,248 & 9.77\%& 10.13\%&3.01\% & 3.13\%\\
 	     3 &152,848& 9.76\%&13.53\% &3.01\%&4.09\% \\
	    6&307,648 &9.78\%&16.85\%&3.01\%&5.27\%\\
		\hline
	\end{tabular}
\hfill
    \caption{\textbf{ Generalization error v.s. different number of ResNet blocks (fixed total training steps $12\times 10^5$).} \fige{pwc}{16}{H^1}{\mwrevise{\eRNN}/\mwrevise{\RNN}}{200}{28}{0.03}. }
    \label{tab:w200_pwc}
\end{table}

\begin{figure}[h!]
\centering
\hspace*{-2.5cm}
\includegraphics[width=1.3\textwidth]{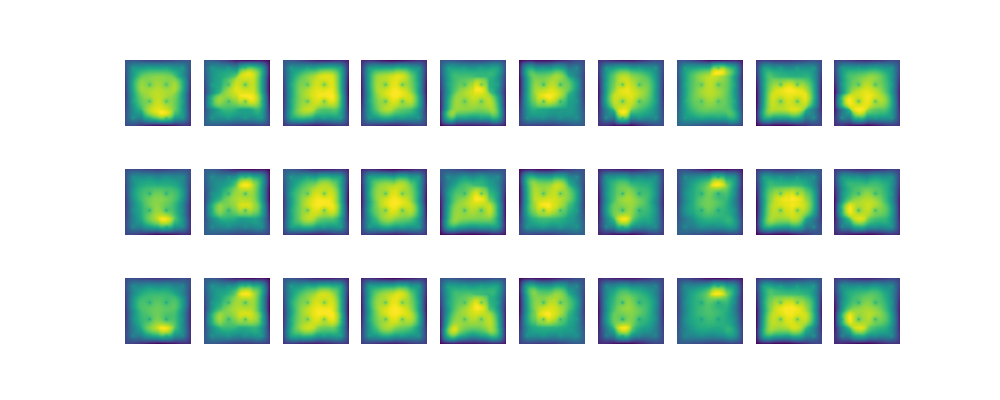}
\caption{\figB{pwc}{16}{H^1}{6}{200}{28}{0.03}. \textbf{Upper row}: reference projected solution $z_h^s = P_{\Wp} u_h^s$; 
\textbf{Middle row}: Prediction of $z_h^s$ with \textbf{\mwrevise{\RNN}};
\textbf{Lower row}: Prediction of $z_h^s$ with \textbf{\mwrevise{\eRNN}}.}
\label{fig:7}
\end{figure}

\newpage
To conclude, our findings can be summarized as follows: When using larger widths, e.g. $200$, \mwrevise{\eRNN} leads to a faster convergence and a somewhat smaller generalization error in comparison with plain \mwrevise{\RNN} training. In fact, while in \mwrevise{\eRNN} the generalization error at least does not increase when increasing network complexity, plain \mwrevise{\RNN} shows a degrading performance reflecting increasing difficulties in realizing expressive potential.
On the other hand, for smaller widths, such as $20$, the performance of both variants is comparable. It should be noted that already a single block and width $=20$ achieves an ``empirical accuracy level" that is improved only slightly by more complex networks (compare Table \ref{tab:w200_pwc} with Table \ref{tab:pwc_w20}). A significant increase in the number of trainable
parameters has not resulted in a significant decrease of training and generalization losses, indicating that \mwrevise{the global update strategy} does not exhaust the
expressive power of the underlying networks. This may rather indicate  
that larger neural network complexity widens  a   plateau  of local minima of about the same magnitude in the loss landscape while parameter choices realizing 
higher accuracy remain isolated and very hard to find. Of course,
this could be affected by different (more expensive) modelities in running \mwrevise{parameter updates} which incidentally would change the implicit regularization mechanism.

\begin{table}[!h]
	\centering
	\begin{tabular}{c|c|c c |cc}
	 	&&\multicolumn{2}{|c}{$\hat{\mE}$}&\multicolumn{2}{|c}{Relative $L_2$ Error of $u_{pred}$}\\
	 	\hline
		\# of blocks&\#of trainable &\mwrevise{\eRNN} & \mwrevise{\RNN}&\mwrevise{\eRNN} & \mwrevise{\RNN}\\
 		\hline
         1 & 1,768 &9.76\%&9.76\% &3.01\%&3.01\% \\
         2 & 3,328 & 9.76\%& 9.76\%&3.01\%&3.01\% \\
 	     3 &4,888& 9.76\%&9.77\%&3.01\% &3.01\%\\
	 6&9,568 &9.76\%&9.76\%&3.01\%&3.01\%\\
		\hline
	\end{tabular}
\hfill
    \caption{\textbf{ Generalization error v.s. different number of ResNet blocks (fixed total training steps $12\times 10^5$).} \fige{pwc}{16}{H^1}{\mwrevise{\eRNN}/\mwrevise{\RNN}}{20}{28}{0.03}. }
    \label{tab:pwc_w20}
\end{table}

\subsection{Log-normal case}\label{sec:log_normal}
As indicated earlier, uniform ellipticity in conjunction with affine parameter dependence of the diffusion coefficients offers very favorable conditions for the type of affine space recovery schemes
described in Section \ref{ssec:1-space}. In particular, affine parameter dependence as well as rapidly decaying Kolmogorov $n$-widths
\eqref{widths} are quite important for methods, resorting to Reduced Bases, to work well, while neural networks are far less dependent
on these preconditions. Therefore, we turn to scenario (S2) which is more challenging in both regards.
%


Specifically, we consider the  case with $a_0= 0$ and $a_1= 1$ in \eqref{eq:model_a_lognormal}. 
Recall that we aim to learn the map form $\bw\to\bc$ where $\bc$ is the POD coefficient vector of the solution in $\mathbb{W}^{\perp}_h$.

First, since the diffusion coefficients no longer depend affinely on the parameters $\by$, rigorously founded error surrogates are no longer computable in an efficient way. This impedes the theoretical foundation as well as the efficiency of methods using Reduced Bases and therefore provides a 
strong motivation exploring alternate methods.
Second, it is not clear whether the solution manifold $\cM$ still have rapidly decaying 
$n$-widths, so that affine spaces of moderate dimension will not give rise to accurate estimators.
In fact, computing the SVD of the snapshot projections $\{z_h^1, \ldots, z_h^{\widehat{N}}\}$, based on the Riesz  representations of the measurement functionals in $\U = H^1$,
shows only very slowly decaying singular values. This indicates that   $P_{\W_h^\perp}(\cM)$
cannot be well approximated by a linear space of moderate dimension.
Thus, when following the above lines, we would have to seek coefficient vectors $\bc$ 
in a space of dimension comparable to  $\widehat{N}$, which renders training \mwrevise{\RNN} prohibitive.
Finally, the diffusion coefficients may near-degenerate degrading uniform ellipticity. 
This raises the question whether $H^1$ is still an appropriate space to accommodate
a reasonable measurement space $\W$ which is at the heart of the choice of sensor coordinates \eqref{deco1}. These adverse effects are reflected by the numerical experiments discussed below.

Therefore, we choose in scenario (S2)
$\U = L_2(\Omega)$ which means 
we are content with a weaker metric for measuring accuracy.
As a consequence,  the representation of the functionals $\bfell(u)$ in $L_2(\Omega)$ 
is the $L_2$-orthogonal projection of these functionals to the truth-space $\U_h$ which then,
as before, span the measurement space $\mathbb{W}_h\subset \U_h$. 
As indicated earlier, we give up on quantifiable  gradient information but facilitate
a more effective approximation of $P_{\W_h^\perp}(\cM)$ where the projection is now
understood in the $L_2$-sense.
In fact, for $\widehat{N} = 1000$ snapshot samples, with the original $H^1$-Riesz-lifting, one needs the $1000$ dominant POD modes to sustain $99.5\%$ of the $H^1$ energy. Instead,
only $21$ dominating modes are required to realize the same accuracy in $L_2(\Omega)$. 
All results presented in this section will be subject to this change. Nevertheless, for comparison we record below in each experiment also the relative $H^1$-error which, as expected, is larger by an order of magnitude.

\subsubsection{\mwrevise{\eRNN} vs. \mwrevise{\RNN} (16 sensors)} \label{sec:log_normal_16}

We are interested to see whether, or under which circumstances, the  advantages of \mwrevise{\eRNN} over plain
\mwrevise{\RNN} training persists also in scenario (S2) where several problem characteristics are 
quite different. We consider similar test conditions as before, namely $16$ uniformly distributed sensors. In total, $1000$ snapshots are collected providing $1000$ synthetic data points, of which $950$  are used for training while $50$ are reserved for evaluation in this case. 
The dimension of observational data $\bw$ is then $m = 16$, while according to
the preceding remarks, the effective complement space dimension, accommodating the coefficients $\bc$,
is $k = 21$. 

\begin{figure}[h!]
    \centering
    \includegraphics[width = 0.5\textwidth]{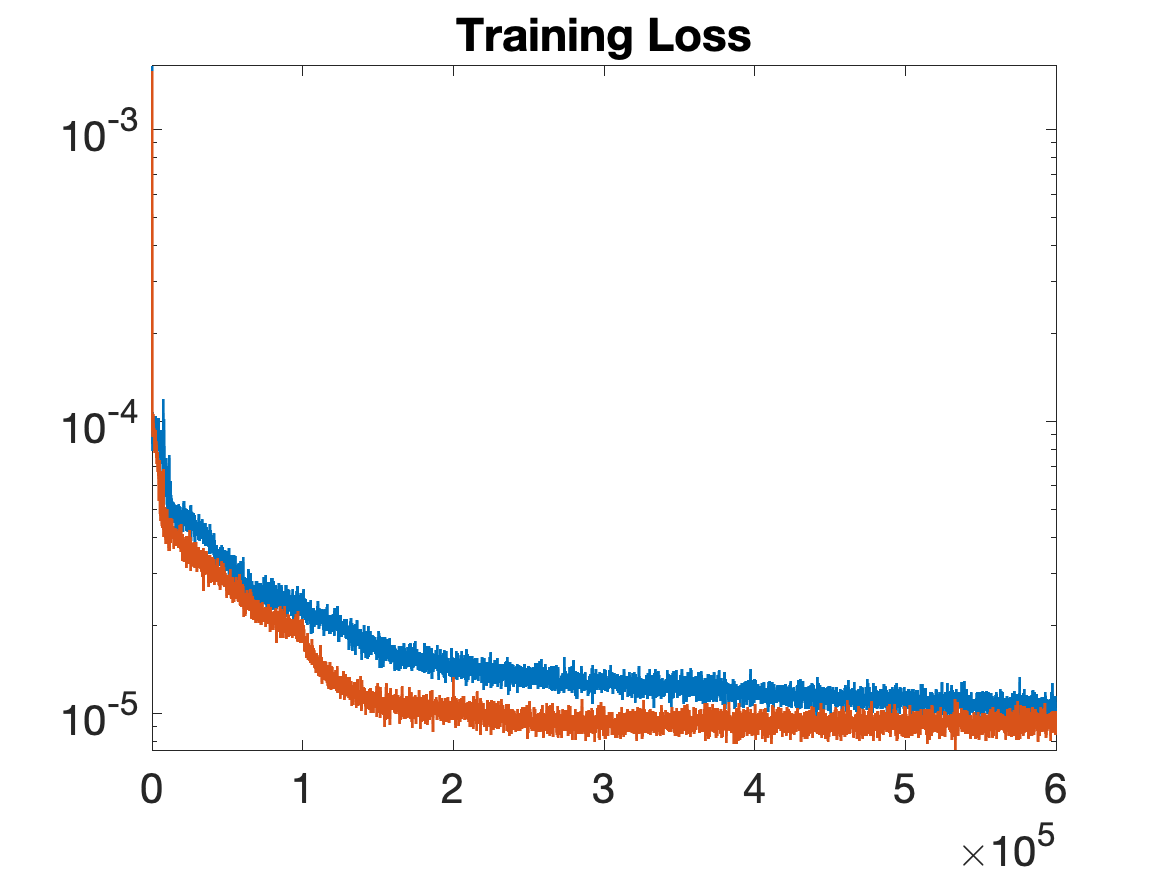}
    \caption{\textbf{\mwrevise{\eRNN (red)}/\mwrevise{\RNN (blue)} training loss comparison}. \figB{log-normal}{16}{L_2}{6}{20}{21}{0.02}.}
    \label{fig:lognormal-ernn_s16}
\end{figure}

The history of training losses is shown in Figure \ref{fig:lognormal-ernn_s16}, where we observe again that the training of \mwrevise{\eRNN} is more efficient than that of  \mwrevise{\RNN}. Specifically, when a new block is introduced in \mwrevise{\eRNN}, the training loss  decays more rapidly (see the corner of the loss curve in Figure \ref{fig:lognormal-ernn_s16} at step $1\times 10^5$). Thus, the expansion strategy is clearly 
beneficial in this case. Moreover, perhaps not surprisingly, we observe that \mwrevise{\RNN} suffers from a slow down in loss decay when training a larger number of parameters simultaneously (Figure \ref{RNN_train_s16_a}). By contrast, for
the block by block optimization in \mwrevise{\eRNN} and fixed width, the number of 
simultaneously updated parameters stays constant so that even for larger (deeper) networks, one  reaches a similar level of training loss at a smaller number of updates  (see Figure \ref{eRNN_train_s16_b}).  Correspondingly, a slightly better overall accuracy of \mwrevise{\eRNN} can be observed in terms of the generalization error (Table \ref{tab:sen16_lognomal_MSE}). One should keep in mind though that
we have allotted a fixed budget of training steps to all variants in this experiment. Thus, increasing depth, reduces the training effort spent on each block which may explain
  the relatively large generalization error obtained for \textbf{$\mathsf{B}12$}. 
Hence, when favoring accuracy improvements at the expense of more training steps, \mwrevise{\eRNN} offers
a clearly better potential while a global training seems to be rather limited. 
One the other hand, the results in Table \ref{tab:sen16_lognomal_MSE} also
indicate that not much gain in accuracy should be expected in this test case by
using more than one or two blocks.

In summary, \mwrevise{\eRNN} appears to offer 
advantages  {in training neural networks with }larger depth  {compared with plain \mwrevise{\RNN}}.
The relatively coarse information provided by the $16$ measurement data 
seems to leave more  room for {an enhanced nonlinearity of deeper networks to capture the manifold component
in $\W^\perp$. This will be seen in Section \ref{section4.2.2} to change somewhat when a larger number of sensors increases the accuracy of the ``zero-order approximation'' provided by the projection $P_{\W^\perp_h}$.}

     

\begin{figure}[!htp]
	\centering
	\begin{subfigure}[b]{0.47\textwidth}
		\centering
		\includegraphics[width=\textwidth]{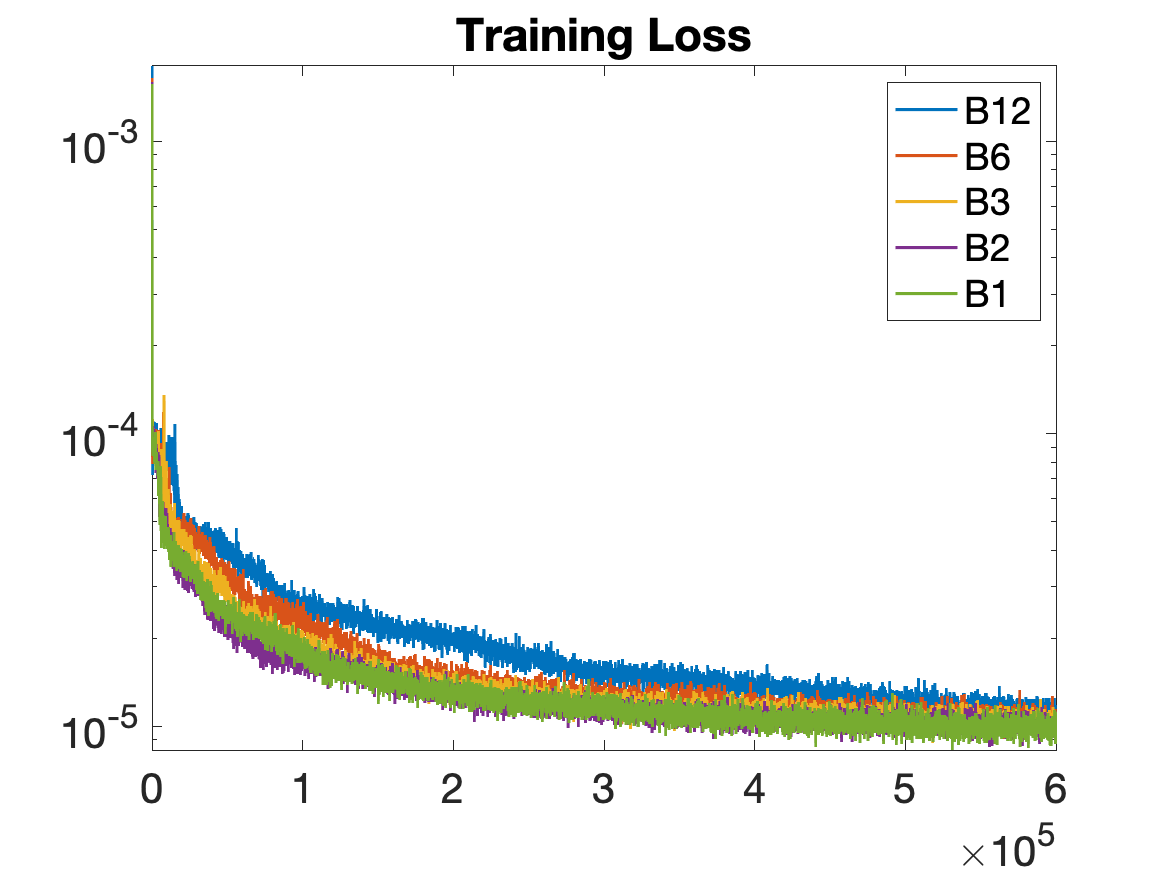}
		\caption{\mwrevise{\textbf{\RNN}} training loss}
		\label{}
	\end{subfigure}
	\hfill
	\begin{subfigure}[b]{0.47\textwidth}
		\centering
		\includegraphics[width=\textwidth]{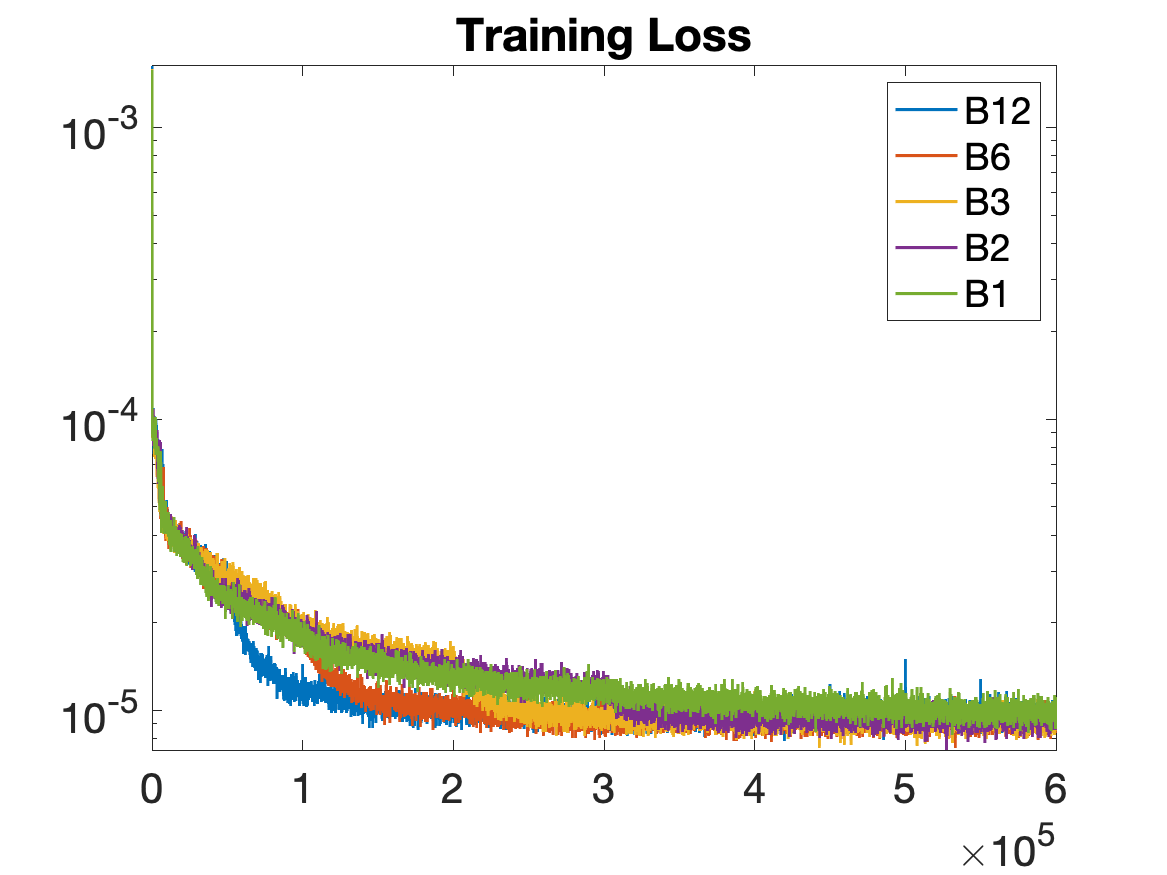}
		\caption{\mwrevise{\textbf{\eRNN}} training loss}
		\label{}
	\end{subfigure}
	\caption{\textbf{Training loss of \mwrevise{\RNN}/\mwrevise{\eRNN} v.s. different number of Blocks}. \figshort{log-normal}{16}{L_2}{20}{21}{0.02}.}
	\label{fig:sen16_different_block}
\end{figure}

\begin{table}[!h]
	\centering
	\begin{tabular}{c|c|c c |cc}
	 	&&\multicolumn{2}{|c}{$\hat{\mE}$}&\multicolumn{2}{|c}{relative $H^1$ Error of $u_{pred}$}\\
	 	\hline
		\# of blocks&\#of trainable &\mwrevise{\eRNN} & \mwrevise{\RNN}&\mwrevise{\eRNN} & \mwrevise{\RNN}\\
 		\hline
         1 & 1,516 &7.89\%&7.89\% & 47.61\% &47.61\% \\
         2 & 2,796 & 7.65\%& 7.94\%&47.33\% &48.26\%\\
 	     3 &4,076& 7.60\%&8.00\% &47.08\%& 48.31\%\\
	     6&7,916 &7.60\%&8.12\%&47.08\%& 48.65\%\\
 	 12&15,596 &7.76\%&8.42\%&47.42\%& 49.28\%\\
		\hline
	\end{tabular}
\hfill
    \caption{\textbf{ Generalization error v.s. different number of ResNet blocks (fixed total training steps $6\times 10^5$).} \fige{log-normal}{16}{L_2}{\mwrevise{\eRNN}/\mwrevise{\RNN}}{20}{21}{0.02}. }
    \label{tab:sen16_lognomal_MSE}
\end{table}

\begin{figure}[h!]
\centering
\hspace*{-2.7cm}
\includegraphics[width=1.3\textwidth]{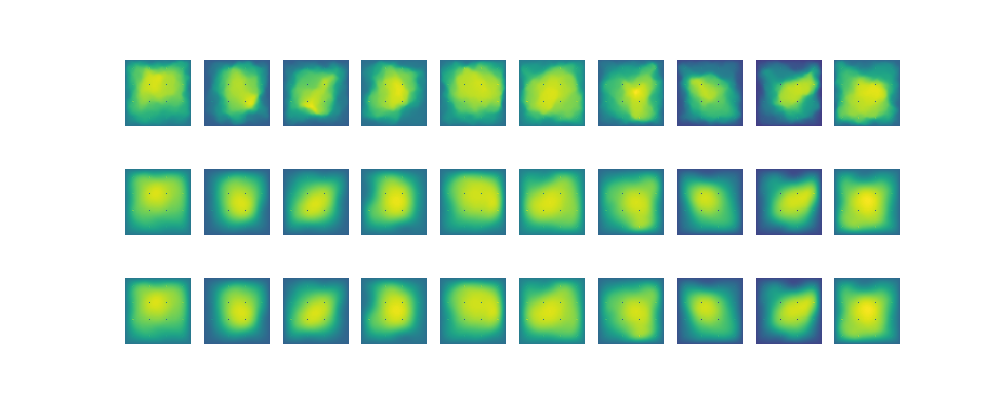}
\caption{\figB{log-normal}{16}{L_2}{6}{20}{21}{0.02}. \textbf{Upper row}: reference projected solution $z_h^s = P_{\mathbb{W}^{\perp}_h} u_h^s$; \textbf{Middle row}: Prediction of $z_h^s$ with \textbf{\mwrevise{\RNN}};
\textbf{Lower row}: Prediction of $z_h^s$ with \textbf{\mwrevise{\eRNN}}.}
\label{fig:universe}
\end{figure}


\subsubsection{\mwrevise{\eRNN} vs. \mwrevise{\RNN} (49 sensors)}\label{section4.2.2} 
In this subsection, we consider $49$ uniformly distributed sensors for measurements. In total, $6000$ snapshots are collected. Among them, $5950$ samples are used for training while $50$ are reserved for evaluation. The dimension of the latent space accommodating $\bc^s$ is now $k = 22$ after applying SVD and keeping $99.5\%$ energy in the $L_2$ sense. 

\begin{figure}[h!]
    \centering
    \includegraphics[width = 0.5\textwidth]{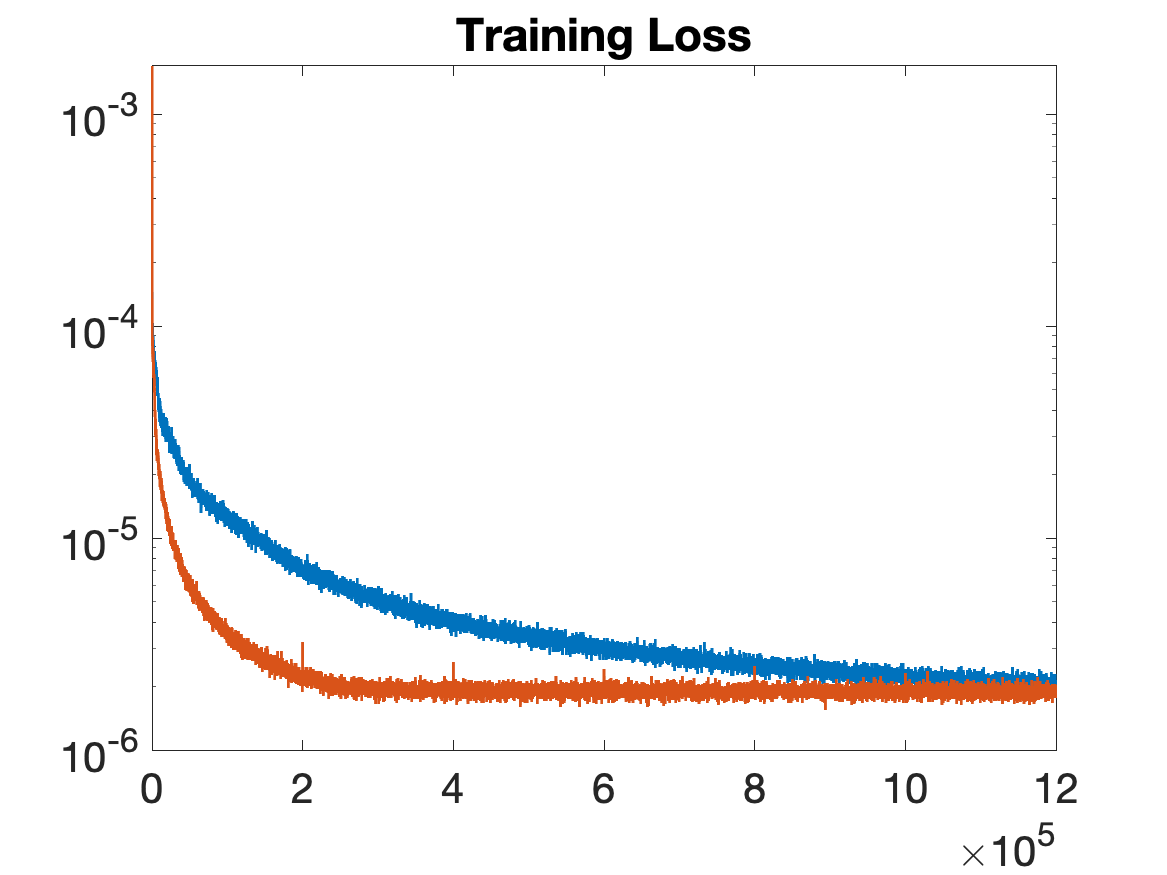}
    \caption{\textbf{\mwrevise{\eRNN (red)}/\mwrevise{\RNN (blue)} training loss comparison}. \figB{log-normal}{49}{L_2}{6}{20}{22}{0.03}.}
    \label{fig:lognormal-ernn}
\end{figure}

However, in this case, from Figure \ref{fig:lognormal-ernn}, we can see that although the training loss of \mwrevise{\eRNN} decays faster compared to \mwrevise{\RNN} at the beginning, both ended up at a similar level. We also do not observe significant changes in decay rates of the loss when employing the expansion strategy
\mwrevise{\eRNN} (see Figure \ref{fig:lognormal-ernn} at step $6\times 10^5$). In fact, the generalization error of \mwrevise{\eRNN} is only slightly smaller than that of \mwrevise{\RNN} (Table \ref{tab:sen_49_lognormal_generalization}). In addition, it is seen that, within the \mwrevise{fixed total number of  training steps}, the best performance is already achieved using a shallow \mwrevise{\RNN}.
This indicates that, within the achievable accuracy range, the map of interest is close to a linear one, given that the ``zero-order'' approximation 
 $P_{\W_h} u$ is now already {rather} accurate. 
However, we do notice that when applying a block by block training strategy in \mwrevise{\eRNN},  {while the} difference in generalization error  {is small}, the savings in training are huge because only one block is trained at a  time, and thus the number of parameters under training is fixed. Thus, for deep neural networks, such a sequential training {scheme} is expected to be beneficial compared to updating all parameters simultaneously.


\begin{figure}[!htp]
	\centering
	\begin{subfigure}[b]{0.47\textwidth}
		\centering
		\includegraphics[width=\textwidth]{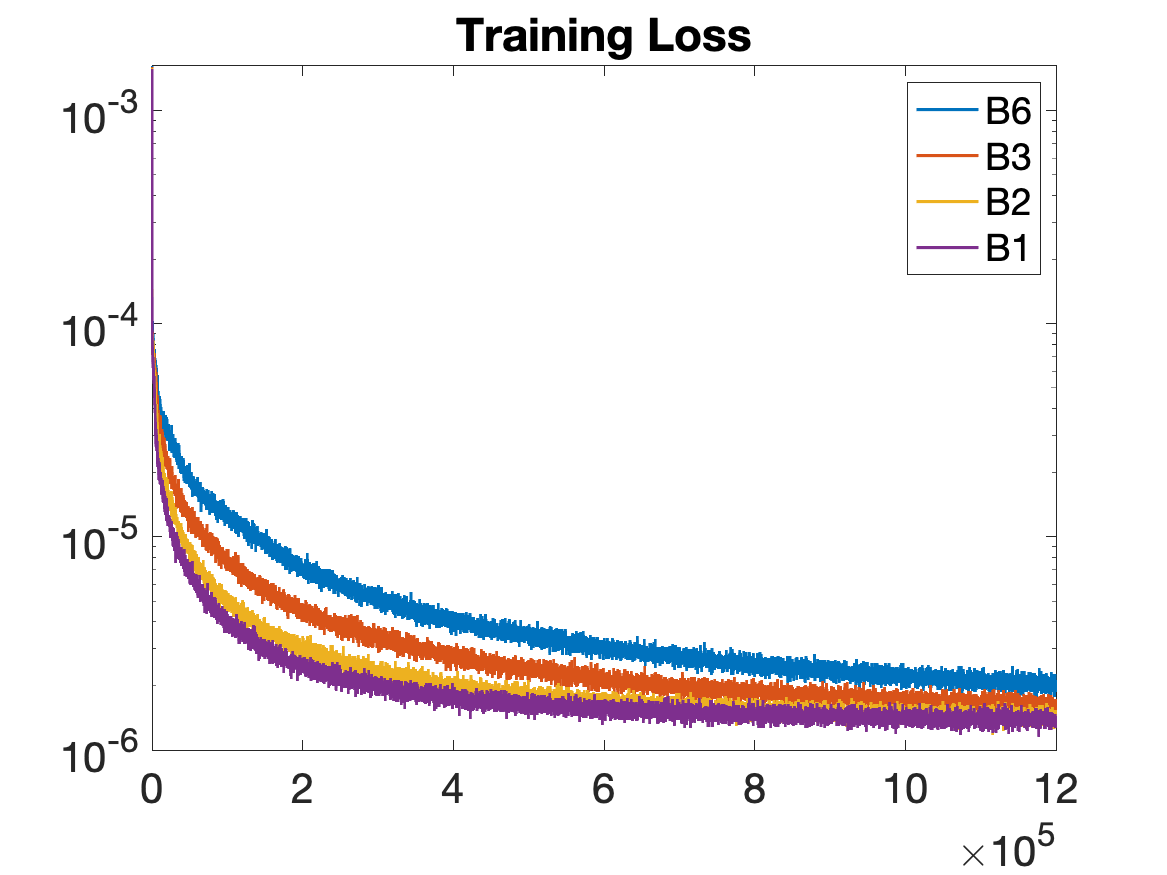}
		\caption{\mwrevise{\textbf{\RNN}} training loss}
		\label{RNN_train}
	\end{subfigure}
	\hfill
	\begin{subfigure}[b]{0.47\textwidth}
		\centering
		\includegraphics[width=\textwidth]{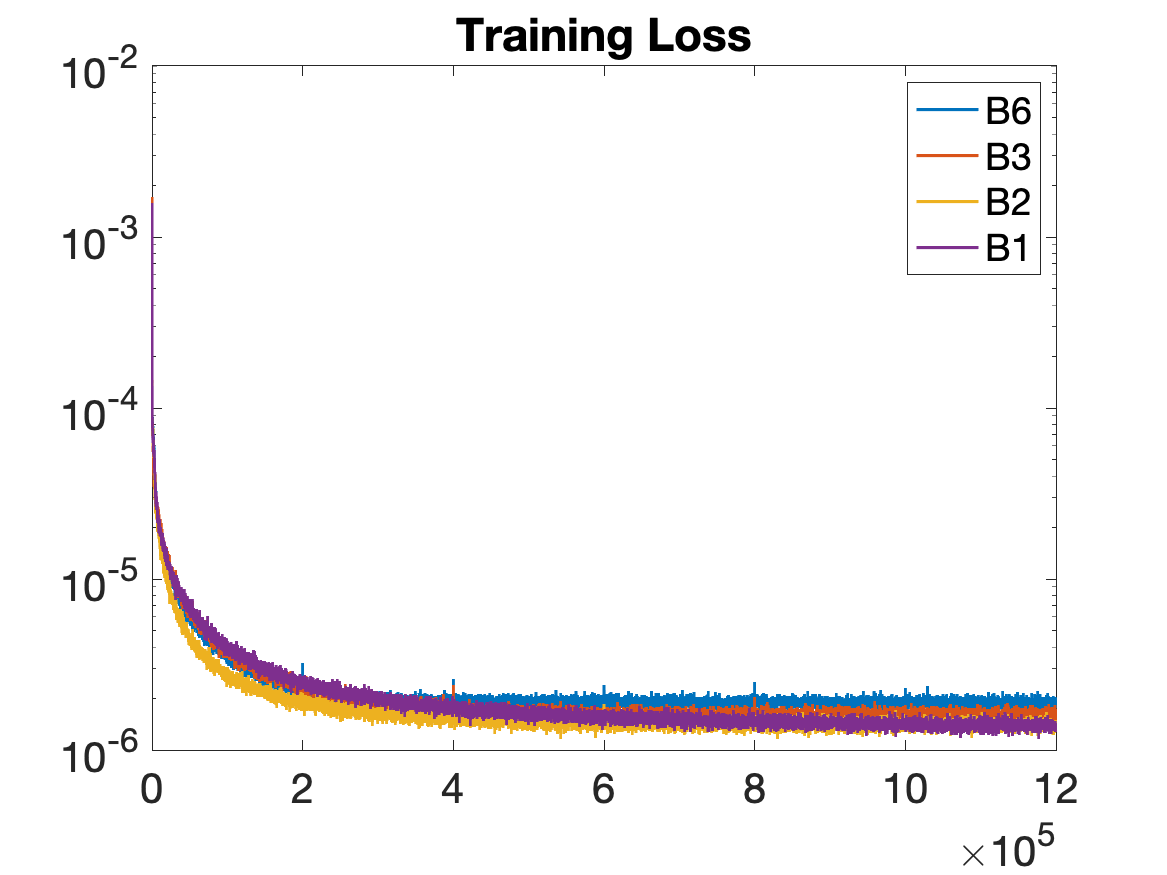}
		\caption{\mwrevise{\textbf{\eRNN}} training loss}
		\label{eRNN_train}
	\end{subfigure}
	\caption{\textbf{Training loss of \mwrevise{\RNN}/\mwrevise{\eRNN} v.s. different number of Blocks}. \figshort{log-normal}{49}{L_2}{20}{22}{0.03}.}
		\label{fig:B2_sen49_compairson}
\end{figure}

\begin{table}[!h]
	\centering
	\begin{tabular}{c|c|c c |cc }
	 	&&\multicolumn{2}{|c}{$\hat{\mE}$}&\multicolumn{2}{|c}{$H^1$ Error of $u_{pred}$}\\
	 	\hline
		\# of blocks&\#of trainable &\mwrevise{\eRNN} & \mwrevise{\RNN}&\mwrevise{\eRNN} & \mwrevise{\RNN}\\
 		\hline
         1 & 2,983 &3.08\%&3.08\% & 52.66\%&52.66\% \\
         2 &4,258 & 3.10\%& 3.13\%&52.64\% &52.71\%\\
 	     3 &5,578& 3.31\%&3.34\% &52.82\%& 53.06\%\\
	 6&9,538 &3.60\%&3.69\%& 53.27\%& 53.90\%\\
		\hline
	\end{tabular}
\hfill
    \caption{\textbf{ Generalization error v.s. different number of ResNet blocks (fixed total training steps $12\times 10^5$).} \fige{log-normal}{49}{L_2}{\mwrevise{\eRNN}/\mwrevise{\RNN}}{20}{22}{0.03}. }
    \label{tab:sen_49_lognormal_generalization}
\end{table}

\begin{figure}[h!]
\centering
\hspace*{-2.7cm}
\includegraphics[width=1.3\textwidth]{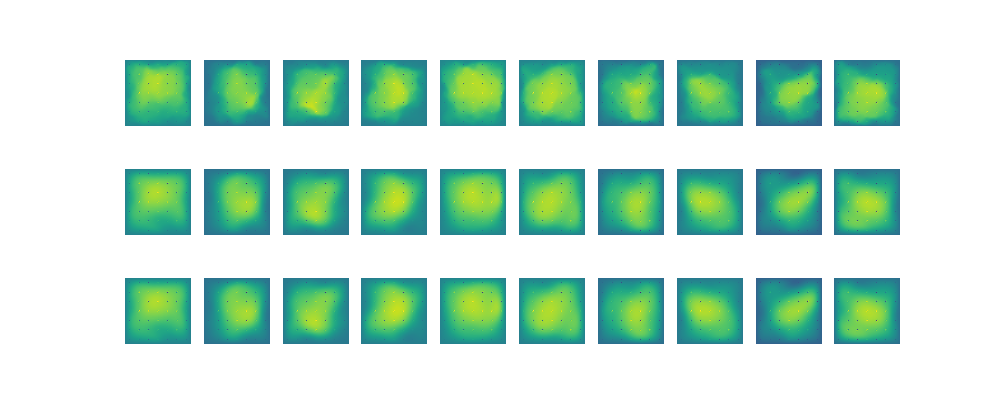}
\caption{\textbf{Upper row}: reference projected solution $z_h^s = P_{\mathbb{W}^{\perp}_h} u_h^s$; \textbf{Middle row}: Prediction of the $z_h^s$ with \textbf{\mwrevise{\RNN}}; \textbf{Lower row}: Prediction  of $z_h^s$ with \textbf{\mwrevise{\eRNN}}.}
\end{figure}


\section{Robustness of \mwrevise{\eRNN} With Regard to Algorithmic Settings}
In this section we further examine how sensitively \mwrevise{\eRNN} depends on various algorithmic settings such as different learning rates or neural network architecture. 
We also wish to see whether extra steps of optimization over all trainable parameters can improve the overall performance, especially the generalization accuracy of \mwrevise{\eRNN}. 
To be consistent, we stick to the example in Section \ref{sec:log_normal_16}, where a log-normal case is considered and $16$ uniformly placed sensors 
give rise to the space $\mathbb{W}_h$ and the corresponding complement $\mathbb{W}_h^T$ in the truth space $\U_h$. 
\subsection{Sensitivity to Learning Rate}
\begin{figure}[h!]
    \centering
     \begin{subfigure}[b]{0.47\textwidth}
        \includegraphics[width =\textwidth]{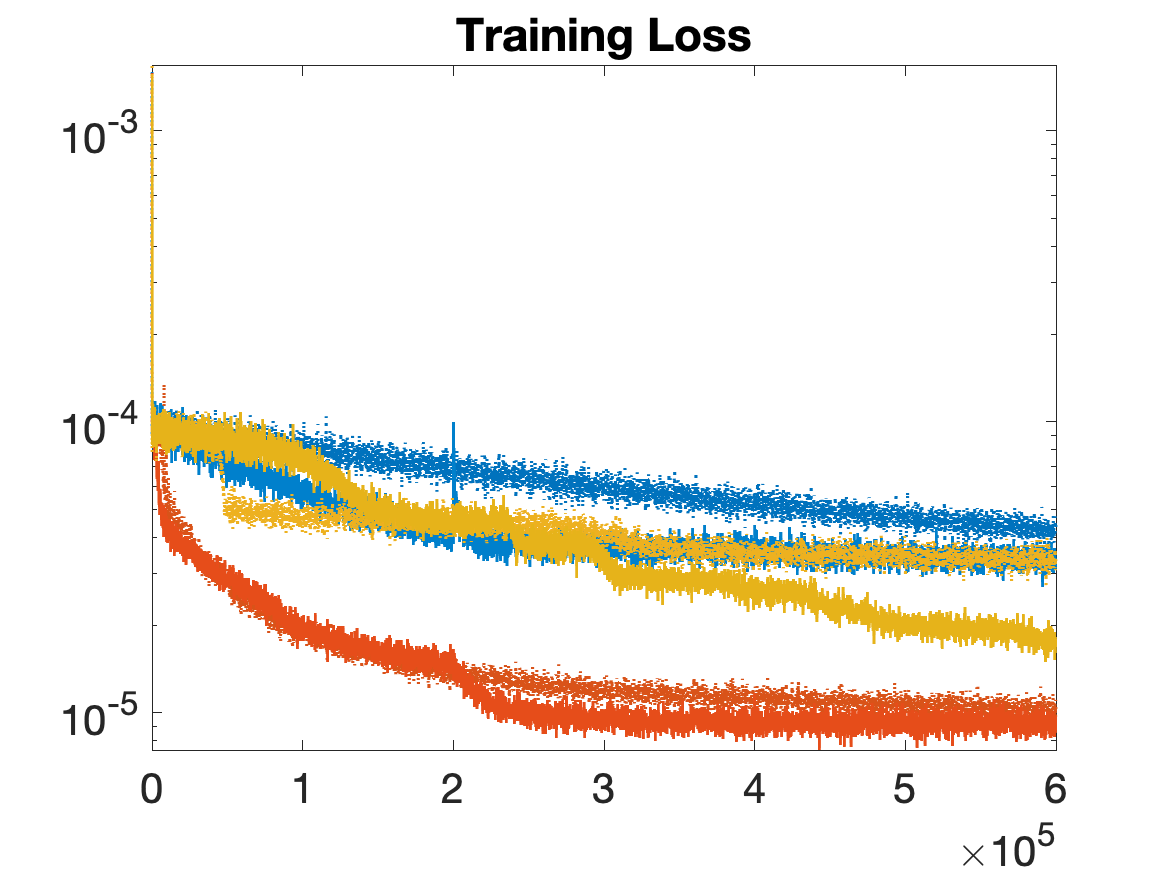}
        \caption{$\mathsf{B}3$}
    \end{subfigure}
    \hfill
    \begin{subfigure}[b]{0.47\textwidth}
        \includegraphics[width =\textwidth]{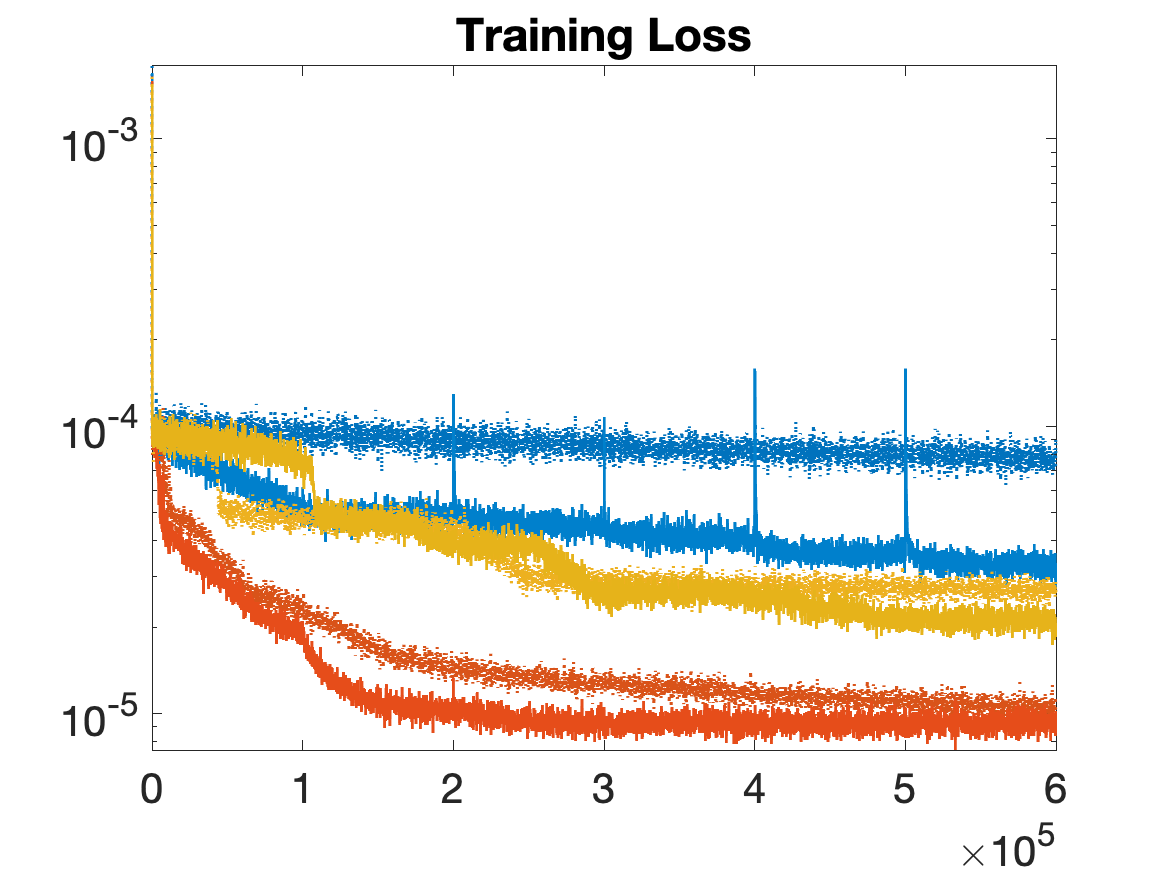}
        \caption{$\mathsf{B}6$}
    \end{subfigure}
    \caption{\textbf{Training loss of \mwrevise{\RNN (doted)}/\mwrevise{\eRNN (solid)} v.s. different learning rate \mwrevise{(blue: $\mathsf{lr} = 0.2$, red: $\mathsf{lr}= 0.02$, yellow: $\mathsf{lr} = 0.2$)}}. $\mathtt{Train}$-\textit{log-normal}-$\mathsf{sen}16$-POD-$L_2$, $\mathsf{W}20$-$\mathsf{O}21$. }
    \label{fig:learning_rate_comparison}
\end{figure}

We first check how \mwrevise{different} learning rates affect the performance of \mwrevise{\eRNN} in comparison with \mwrevise{\RNN}. \mwrevise{In particular, the different learning rates considered are merely an initialization of the learning rates that are applied during the training process. The optimizer of our choice (AdaGrad, the adaptive gradient algorithm) is a modified stochastic gradient descent algorithm which will automatically adjust the learning rate per parameter as the training proceeds.}
Corresponding findings can be summarized as follows (see Figure \ref{fig:learning_rate_comparison} and
Table \ref{tab:lr_comparison}):\\
\begin{itemize}
    \item The training loss of \mwrevise{\eRNN} converges faster for a wide range of {\em constant initial} learning rates.
    \item The generalization errors for \mwrevise{\eRNN} are smaller.
\end{itemize}

\begin{table}[!h]
	\centering
	\begin{subtable}[t]{\textwidth}
	\centering
	\begin{tabular}{c|c c |cc }
	 	&\multicolumn{2}{|c}{$\hat{\mE}$}&\multicolumn{2}{|c}{Relative $H^1$ Error of $u_{pred}$}\\
	 	\hline
		lr&\mwrevise{\eRNN} & \mwrevise{\RNN}&\mwrevise{\eRNN} & \mwrevise{\RNN}\\
 		\hline
         0.002 &10.75\%&14.50\%&54.25\%&61.29\% \\
         0.02 & 7.60\%& 8.00\%&47.08\% &48.31\%\\
 	     0.2 & 14.78\%&16.66\% &60.62\%& 62.96\%\\
		\hline
	\end{tabular}
	\caption{$\mathsf{B}3$}
	\end{subtable}
	\begin{subtable}[t]{\textwidth}
	\centering
		\begin{tabular}{c|cc |cc }
	 	&\multicolumn{2}{|c}{$\hat{\mE}$}&\multicolumn{2}{|c}{Relative $H^1$ Error of $u_{pred}$}\\
	 	\hline
		lr&\mwrevise{\eRNN} & \mwrevise{\RNN}&\mwrevise{\eRNN} & \mwrevise{\RNN}\\
 		\hline
         0.002 &11.76\%&13.38\% &55.66\%&60.13\% \\
         0.02 & 7.60\%& 8.12\%&47.08\% &48.65\%\\
 	     0.2 & 14.53\%&23.44\%&60.80\%& 76.68\%\\
		\hline
	\end{tabular}
	\caption{$\mathsf{B}6$}
	\end{subtable}
    \caption{\textbf{ Generalization error v.s. different learning rates (fixed total training steps $6\times 10^5$).} $\mathtt{Train}$-\textit{log-normal}-$\mathsf{sen}16$-POD-$L_2$,$\mathsf{W}20$-$\mathsf{O}21$. }
    \label{tab:lr_comparison}
\end{table}

\subsection{Dependence on Neural Network Architecture}

In this Subsection we explore the effect of varying the architecture again for the application in scenario (S2).
\subsubsection{Width}

\begin{figure}[h!]
    \centering
    \begin{subfigure}[b]{0.47\textwidth}
        \includegraphics[width =\textwidth]{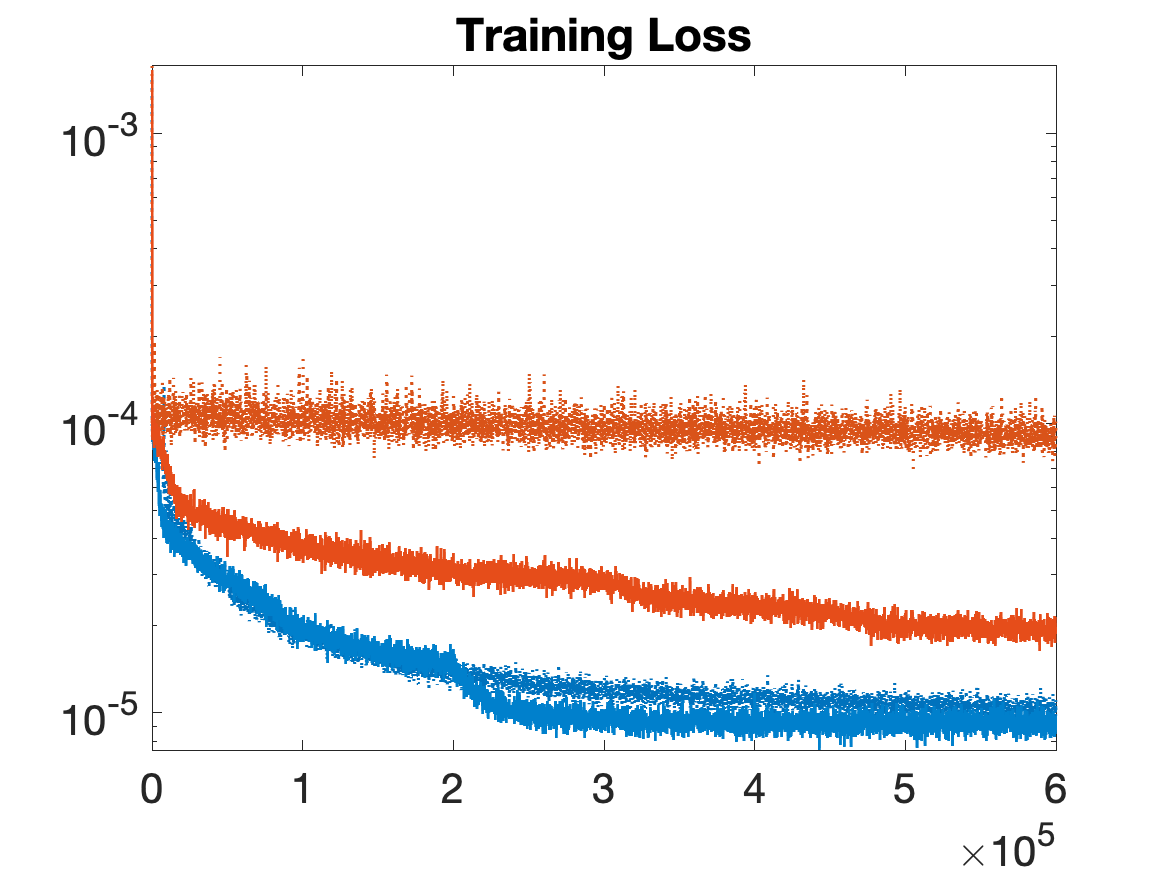}
        \caption{$\mathtt{B3}$}
    \end{subfigure}
    \begin{subfigure}[b]{0.47\textwidth}
        \includegraphics[width =\textwidth]{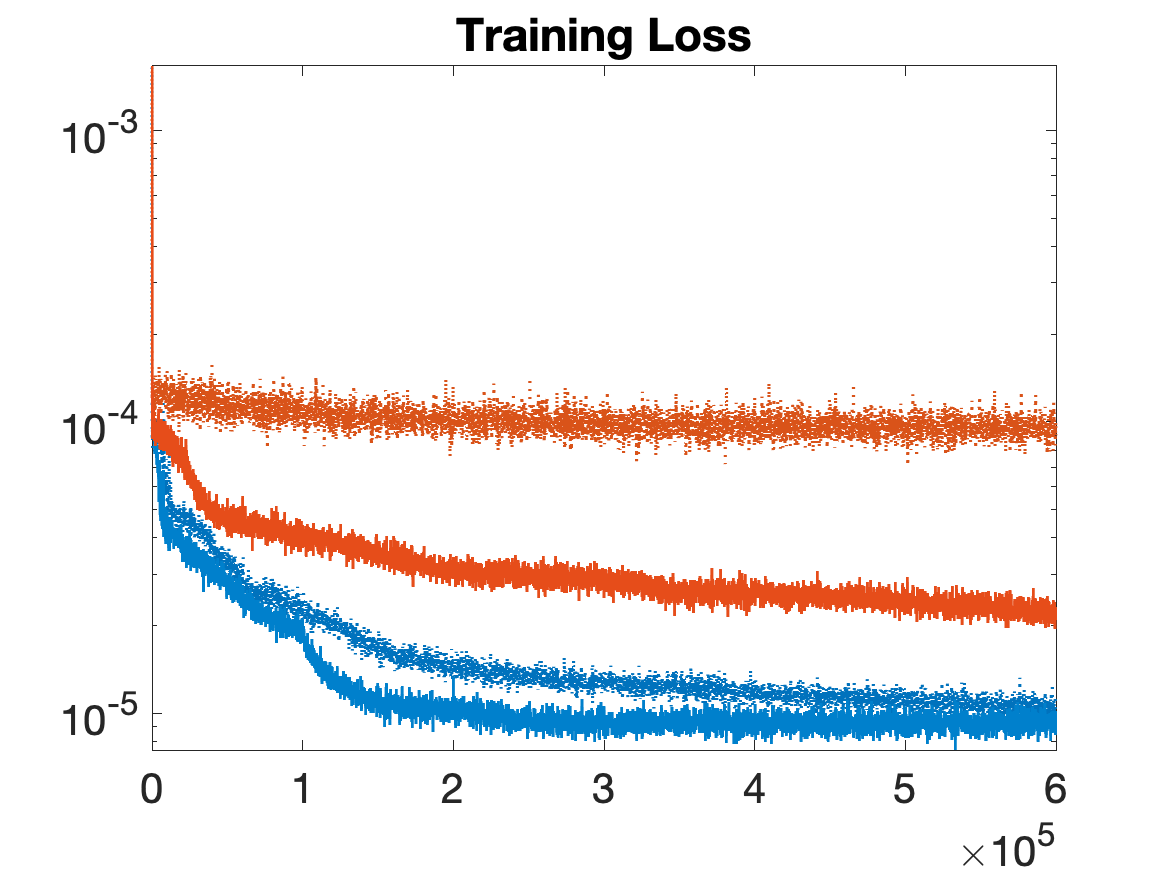}
        \caption{$\mathtt{B6}$}
    \end{subfigure}
    \caption{\textbf{Training loss of \mwrevise{\RNN (doted)}/\mwrevise{\eRNN (solid)} v.s. different width \mwrevise{(blue: $\mathsf{W}20$, red: $\mathsf{W}200$)}}. $\mathtt{Train}$-\textit{log-normal}-$\mathsf{sen}16$-POD-$L_2$, $\mathsf{O}21$-$\mathsf{lr}0.02$}
    \label{fig:width_comparison}
\end{figure}

As indicated by Figure \ref{fig:width_comparison} and as expected, neural networks of equal depth but larger widths require more training iterations  to converge and each update is computationally  more intense than for narrow ones. For both narrow and wide neural networks, \mwrevise{\eRNN} appears to entail a faster decay of the loss during the training process as well as a more accurate prediction.

\begin{table}[!h]
	\centering
	\centering
	\begin{subtable}[t]{\textwidth}
	\centering
	\begin{tabular}{c|c c |cc }
	 	&\multicolumn{2}{|c}{$\hat{\mE}$}&\multicolumn{2}{|c}{Relative $H^1$ Error of $u_{pred}$}\\
	 	\hline
		Width &\mwrevise{\eRNN} & \mwrevise{\RNN}&\mwrevise{\eRNN} & \mwrevise{\RNN}\\
 		\hline
         $\mathsf{W}20$ & 7.60\%& 8.00\%&47.08\% &48.31\%\\
         $\mathsf{W}200$ & 11.41\%& 25.10\%&55.17\% &77.85\%\\
		\hline
	\end{tabular}
	\caption{$\mathsf{B}3$}
	\end{subtable}
	\begin{subtable}[t]{\textwidth}
		\centering
	\begin{tabular}{c|c c |cc}
	 	&\multicolumn{2}{|c}{$\hat{\mE}$}&\multicolumn{2}{|c}{Relative $H^1$ Error of $u_{pred}$}\\
	 	\hline
		Width &\mwrevise{\eRNN} & \mwrevise{\RNN}&\mwrevise{\eRNN} & \mwrevise{\RNN}\\
 		\hline
         $\mathsf{W}20$  & 7.60\%& 8.12\%&47.08\% &48.65\%\\
         $\mathsf{W}200$ &12.06\%& 24.94\%& 56.24\%&75.99\% \\
		\hline
	\end{tabular}
		\caption{$\mathsf{B}6$}
	\end{subtable}
    \caption{\textbf{ Generalization error v.s. different neural network width (fixed total training steps $6\times 10^5$).} $\mathtt{Train}$-\textit{log-normal}-$\mathsf{sen}16$-POD-$L_2$, $\mathsf{O}21$-$\mathsf{lr}0.02$. }
\end{table}

\subsubsection{Dependence on Depth}
\begin{figure}[h!]
    \centering
    \includegraphics[width = 0.6\textwidth]{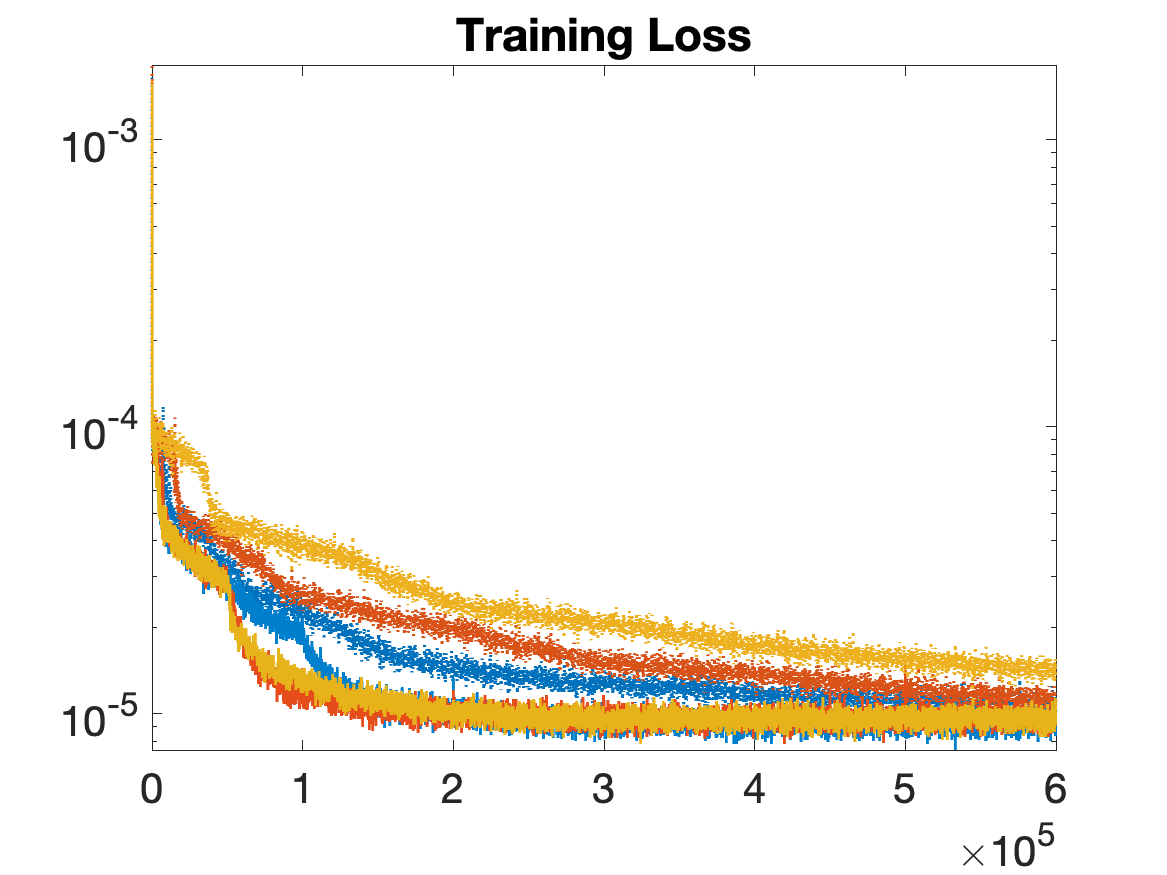}
   \caption{\textbf{Training loss of \mwrevise{\RNN (doted)}/\mwrevise{\eRNN (solid)} v.s. different depth (\# of blocks, \mwrevise{blue: $\mathsf{B}6$, red: $\mathsf{B}12$, yellow: $\mathsf{B}14$})}. $\mathtt{Train}$-\textit{log-normal}-$\mathsf{sen}16$-POD-$L_2$, $\mathsf{W}20$-$\mathsf{O}21$-$\mathsf{lr}0.02$}
    \label{fig:depth_comparison}
\end{figure}

Similar to what can be observed from Figure \ref{fig:sen16_different_block}, we observe that,
for the given test problems, deeper \mwrevise{\eRNN}s can converge at a similar rate as shallower ones.
By contrast, perhaps not surprisingly, the convergence of plain \mwrevise{\RNN}s will   generally slow down as depth increases. In other words, the deeper, the slower is the convergence of \mwrevise{\RNN}. see Figure \ref{fig:depth_comparison}.

\begin{table}[!h]
	\centering
	\centering
	\begin{tabular}{c|c c |cc}
	 	&\multicolumn{2}{|c}{$\hat{\mE}$}&\multicolumn{2}{|c}{Relative $H^1$ Error of $u_{pred}$}\\
	 	\hline
		\# of Blocks &\mwrevise{\eRNN} & \mwrevise{\RNN}&\mwrevise{\eRNN} & \mwrevise{\RNN}\\
 		\hline
         $\mathsf{B}6$ & 7.60\%& 8.12\%&47.08\% &48.65\%\\
         $\mathsf{B}12$ & 7.76\%& 8.42\%&47.42\% &49.28\%\\
         $\mathsf{B}14$ & 8.17\%& 9.40\%&47.78\% &51.00\%\\
		\hline
	\end{tabular}
    \caption{\textbf{ Generalization error v.s. different number of ResNet blocks (fixed total training steps $6\times 10^5$).} $\mathtt{Train}$-\textit{log-normal}-$\mathsf{sen}16$-POD-$L_2$, $\mathsf{W}20$-$\mathsf{O}21$-$\mathsf{lr}0.02$. }
\end{table}

\subsection{Duration of Training}
The above discussion of the influence of width and depth on the learning outcome is based on an assumption that a fixed \mwrevise{number of training steps is applied.}
It is, in general, completely unclear whether such a budget \mwrevise{of training steps} suffices to exploit the expressive power of the network.
We are therefore interested to see how \mwrevise{\eRNN} and \mwrevise{\RNN} compare when lifting the complexity constraints.
This is the more interesting as an optimization step on a single block is not quite comparable 
with a descent step over the whole network. This discrepancy increases of course with
increasing depth. As a first step in this direction, we inspect the effect of quadrupling 
\mwrevise{the total number of training \mwrevise{steps}  to} $2.4\times 10^{6}$.

It turns out that earlier findings are confirmed. Eventually, given enough training
time and effort, \mwrevise{\RNN} can achieve about the same accuracy as \mwrevise{\eRNN} which indicates that the ``achievable'' expressivity offered by the \mwrevise{\RNN} architecture has been exploited by both optimization strategies. 
There is a slight gain of accuracy in comparison with the previous cap of $6\times 10^5$,
namely $0.02\%\sim 0.03\%$.\\
\begin{table}[!h]
	\centering
	\centering
	\begin{tabular}{c|cc|cc}
	 	&\multicolumn{2}{|c}{$\hat{\mE}$}&\multicolumn{2}{|c}{Relative $H^1$ Error of $u_{pred}$}\\
	 	\hline
		Total training steps &\mwrevise{\eRNN} & \mwrevise{\RNN}&\mwrevise{\eRNN} & \mwrevise{\RNN}\\
 		\hline
         $6\times 10^5$ & 7.60\%& 8.12\%&47.08\% &48.65\%\\
         $24\times 10^5$ &7.75\%& 7.74\%&47.10\%&47.09\%\\
		\hline
	\end{tabular}
    \caption{\textbf{ Generalization error v.s. different number of training steps.} $\mathtt{Train}$-\textit{log-normal}-$\mathsf{sen}16$-POD-$L_2$, $\mathsf{B}6$-$\mathsf{W}20$-$\mathsf{O}21$-$\mathsf{lr}0.02$. }
\end{table}

\begin{figure}[h!]
    \centering
    \includegraphics[width = 0.6\textwidth]{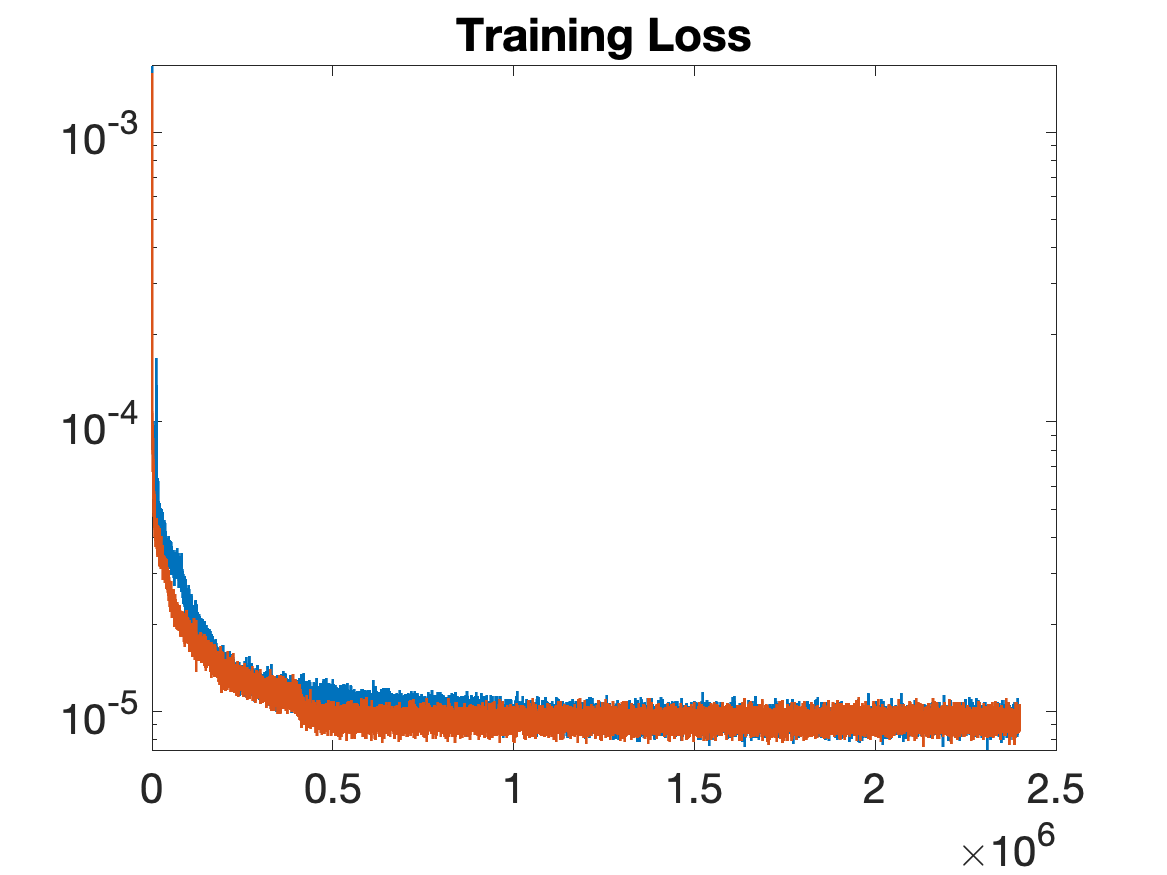}
     \caption{\textbf{Training loss of \mwrevise{\RNN (red)}/\mwrevise{\eRNN (blue)} for large number of training steps ($2.4\times 10^6$)}. $\mathtt{Train}$-\textit{log-normal}-$\mathsf{sen}16$-POD-$L_2$, $\mathsf{B}6$-$\mathsf{W}20$-$\mathsf{O}21$-$\mathsf{lr}0.02$}
    \label{fig:my_label}
\end{figure}

\subsection{Training Schedules}\label{num:sweep}
\begin{figure}
    \centering
    \hspace{-0.5cm}
    \includegraphics[width = 0.6\textwidth]{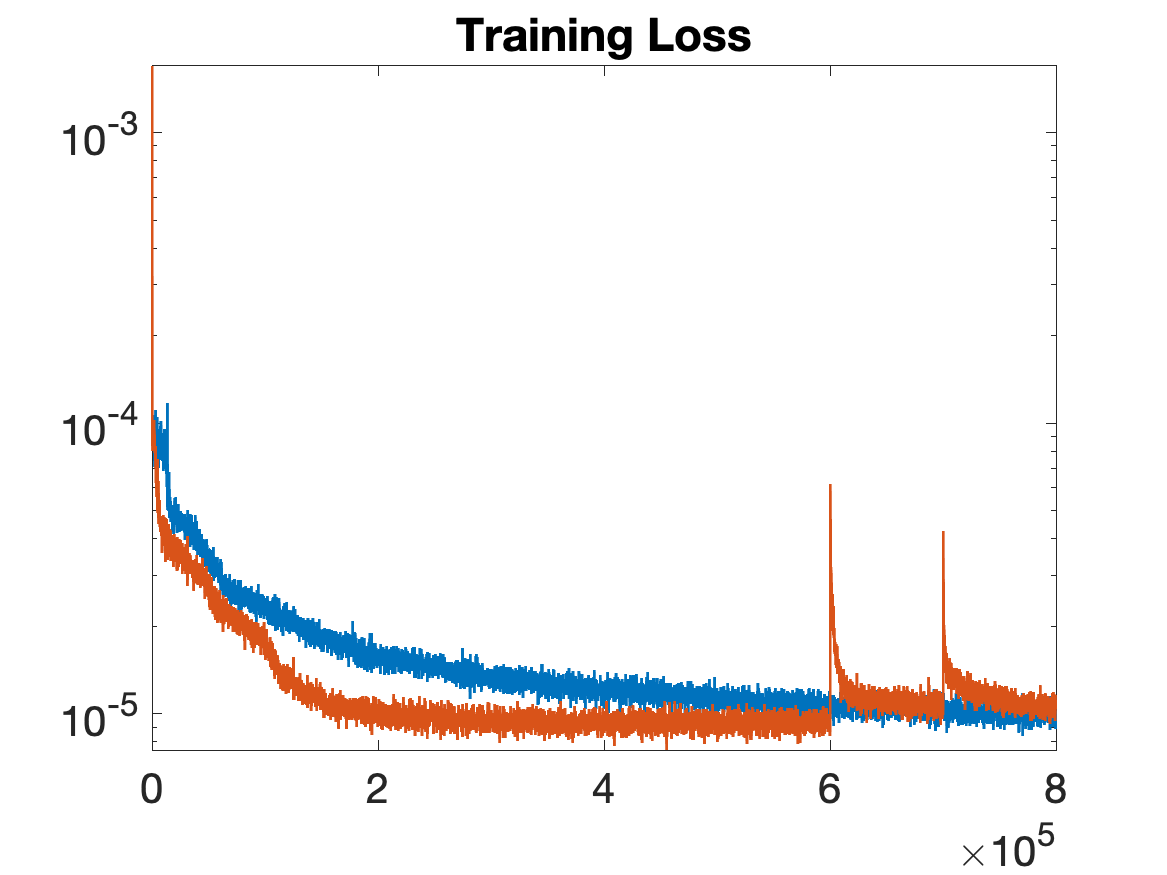}
     \caption{\textbf{Training loss of \mwrevise{\eRNN (red)}/\mwrevise{\RNN (blue)} with $2\times10^5$ steps of global updates. $\mathtt{Train}$-\textit{log-normal}-$\mathsf{sen}16$-POD-$L_2$, $\mathsf{B}6$-$\mathsf{W}20$-$\mathsf{O}21$-$\mathsf{lr}0.02$}}.
    \label{fig:sweep}
\end{figure}

In this example, we add an additional $2\times 10^5$ training steps to update all trainable parameter simultaneously \mwrevise{in addition to a block-wise training}.  Compared with a pure block by block training schedule (see Figure \ref{fig:lognormal-ernn_s16}), the additional global optimization effort
 does not seem to improve on the training success but rather worsens it (see Figure \ref{fig:sweep}).

\begin{table}[!h]
	\centering
	\begin{tabular}{p{5cm}|cc|cc }
	 	&\multicolumn{2}{|c}{$\hat{\mE}$}&\multicolumn{2}{|c}{Relative $H^1$ Error of $u_{pred}$}\\
	 	\hline
		Total training steps &\mwrevise{\eRNN} & \mwrevise{\RNN}&\mwrevise{\eRNN} & \mwrevise{\RNN}\\
 		\hline
         {$6\times 10^5$ block-wise updates}  & 7.60\%& 8.12\%&47.08\% &48.65\%\\
         {$6\times 10^5$ block-wise updates + $2\times 10^5$ global updates }  &7.75\%& 7.99\% &47.10\%&47.91\%\\
		\hline
	\end{tabular}
    \caption{\textbf{ Generalization error with/without $2\times 10^5$ global updates.} $\mathtt{Train}$-\textit{log-normal}-$\mathsf{sen}16$-POD-$L_2$, $\mathsf{B}6$-$\mathsf{W}20$-$\mathsf{O}21$-$\mathsf{lr}0.02$.}
\end{table}

\subsection{Wide Neural Network Subject to Long-term Training}
We  check at last whether the neural networks can do better than in previous experiments in terms of approximating the map between $\mathbb{W}_h$ and $\mathbb{W}^{\perp}_h$ when significantly increasing training time.
\begin{figure}[h!]
    \centering
    \includegraphics[width =0.6\textwidth]{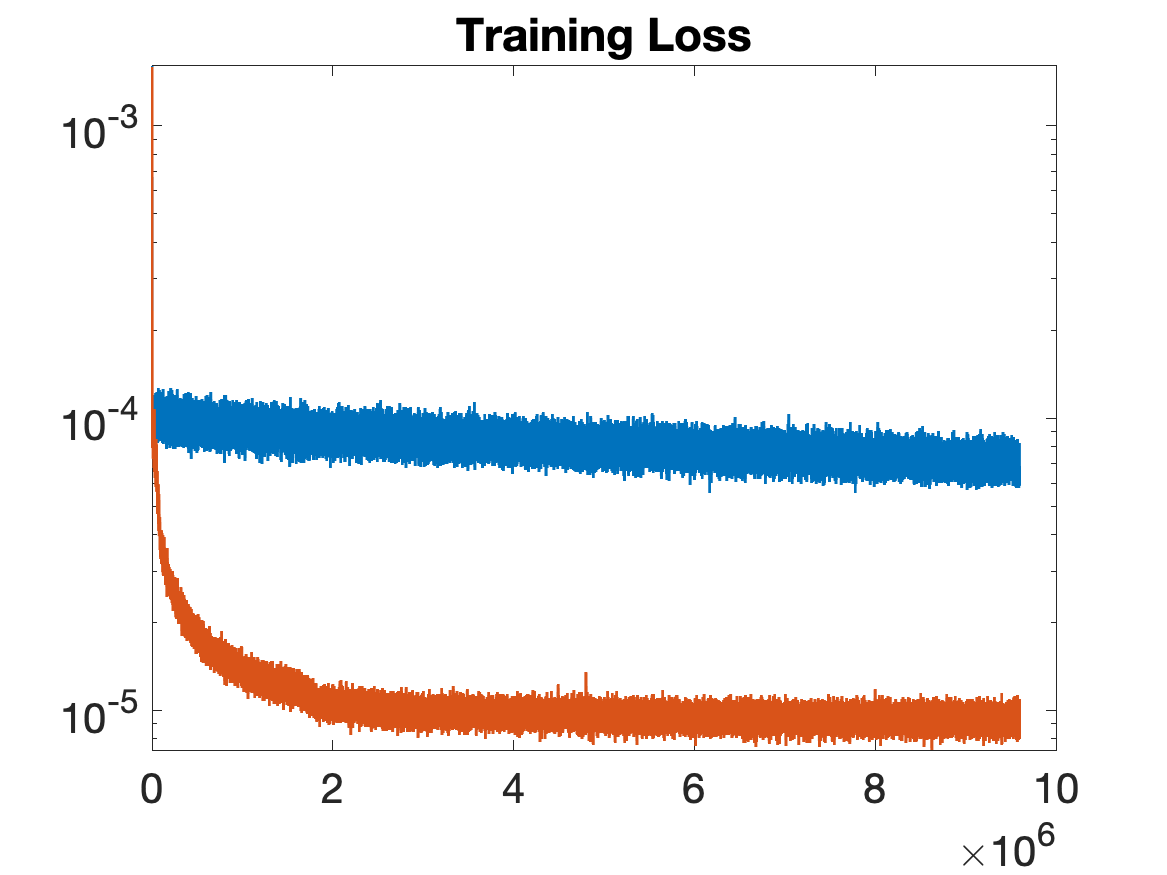}
    \caption{\textbf{Training loss of \mwrevise{\eRNN(red)}/\mwrevise{\RNN(blue)} for large number of training steps ($9.6\times 10^6$) of wide neural network ($\mathsf{W}200$)}. $\mathtt{Train}$-\textit{log-normal}-$\mathsf{sen}16$-POD-$L_2$, $\mathsf{B}6$-$\mathsf{W}200$-$\mathsf{O}21$-$\mathsf{lr}0.02$}
\end{figure}
\begin{table}[!h]
	\centering
	{
	\centering
	\begin{tabular}{c|cc|cc}
	 	&\multicolumn{2}{|c}{$\hat{\mE}$}&\multicolumn{2}{|c}{Relative $H^1$ Error of $u_{pred}$}\\
	 	\hline
		Total training steps &\mwrevise{\eRNN} & \mwrevise{\RNN}&\mwrevise{\eRNN} & \mwrevise{\RNN}\\
 		\hline
         short ($6\times 10^5$) \& narrow ($\mathsf{W}20$)& 7.60\%& 8.12\%&47.08\% &48.65\%\\
         long ($24\times 10^5$) \& narrow ($\mathsf{W}20$) &7.75\%& 7.74\%&47.10\%&47.09\%\\
         short ($6\times 10^5$) \& wide ($\mathsf{W}200$) &12.06\%& 24.94\%& 56.24\%&75.99\% \\
         long ($24\times 10^5$) \& wide ($\mathsf{W}200$) &8.82\%& 23.09\%&49.66\%&78.19\% \\
         super long ($96\times 10^5$) \& wide ($\mathsf{W}200$) &7.77\%& 21.83\%& 47.17\%&71.95\% \\
		\hline
	\end{tabular}
	}
    \caption{\textbf{ Generalization error v.s. different number of training steps.} $\mathtt{Train}$-\textit{log-normal}-$\mathsf{sen}16$-POD-$L_2$, $\mathsf{B}6$-$\mathsf{O}21$-$\mathsf{lr}0.02$. }
    \label{tab:super_long_time}
\end{table}
The results show that earlier findings persist. Table \ref{tab:super_long_time} shows that, while with \mwrevise{\eRNN} there is no real benefit of larger width, at least estimation quality does not apper to degrade over long training periods.  Instead, \mwrevise{\RNN} seems to be adversely affected by  larger network complexity.  

More specifically, Table \ref{tab:super_long_time} shows that the smallest generalization error is achieved by the relatively low training effort for 
narrow networks. Much larger networks, instead, seem to achieve about that same 
accuracy level only at the expense of a significantly larger training effort,
leaving little hope for substantial further accuracy improvements by continued
training. In contrast, globally updating corresponding networks achieves the best result again for narrow networks but, in agreement with earlier tests, at the expense of four times as many training steps than for block-wise training.
For wide networks even extensive training effort does not seem to reproduce accuracy levels, achieved earlier for smaller networks.

\subsection{Optimizer}\label{optimizer}
\begin{figure}[h!]
    \centering
    \includegraphics[width=0.7\linewidth]{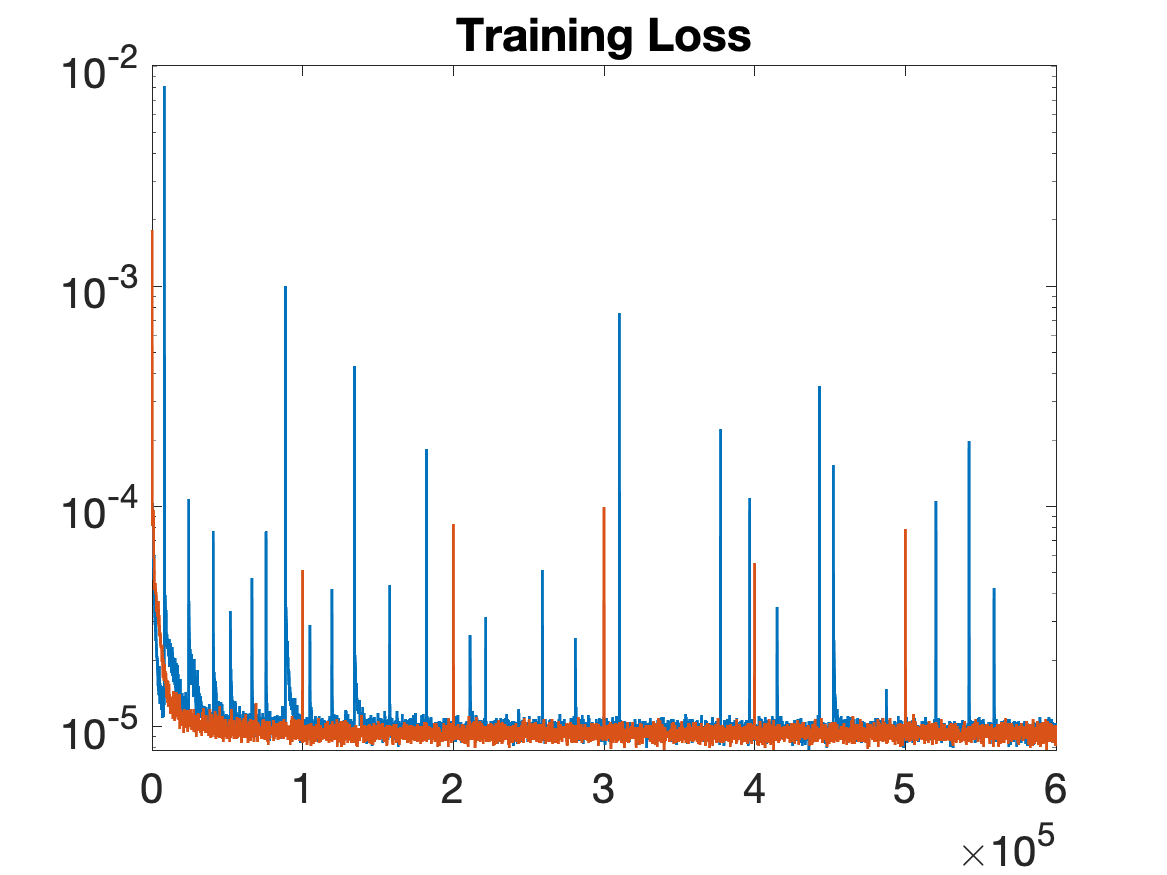}
    \caption{\mwrevise{\textbf{Training history using Adam.} \fige{pwc}{16}{H^1}{{\eRNN (red)}/{\RNN (blue)}}{20}{21}{0.02/0.002}. }}
    \label{fig:adma}
\end{figure}

\begin{figure}[h!]
\centering
\begin{subfigure}{.45\textwidth}
  \centering
  \includegraphics[width=\linewidth]{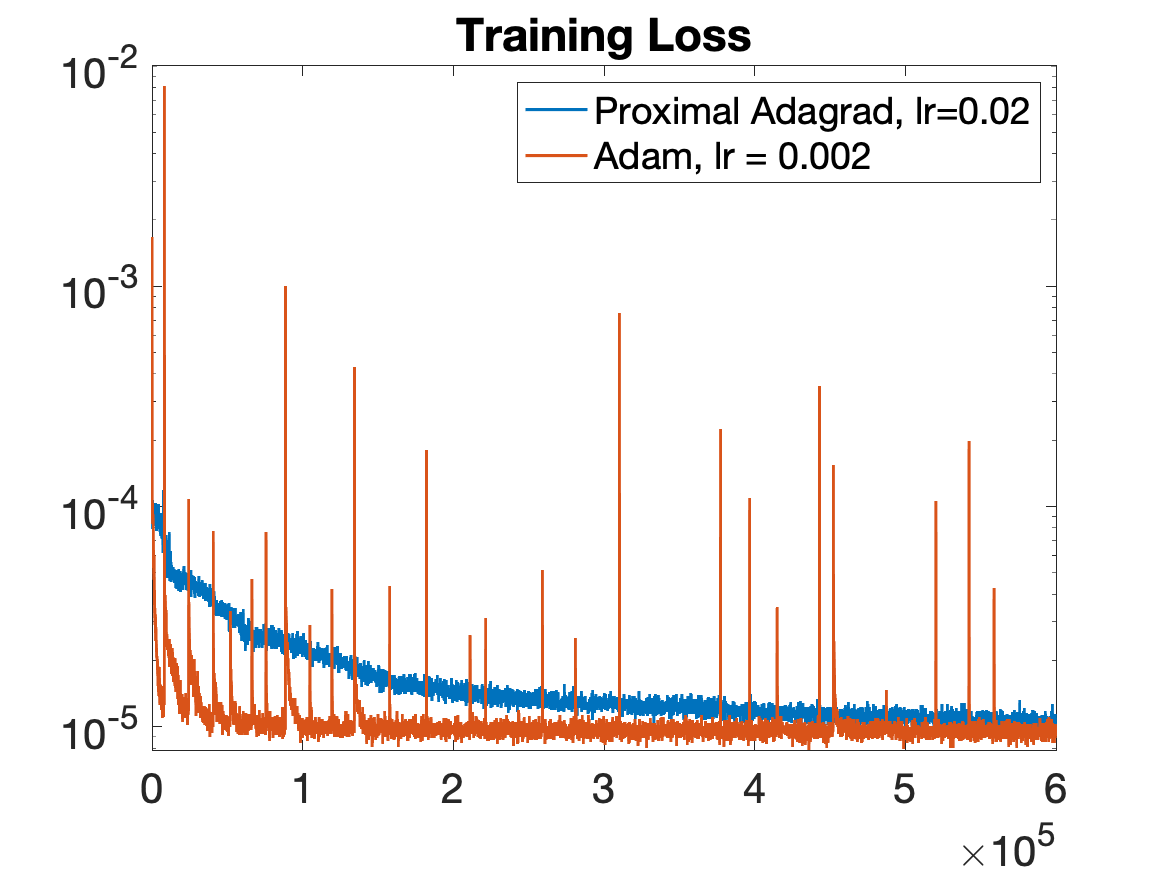}
    \caption{\RNN}
\end{subfigure}
\hfil
\begin{subfigure}{.45\textwidth}
  \centering
  \includegraphics[width=\linewidth]{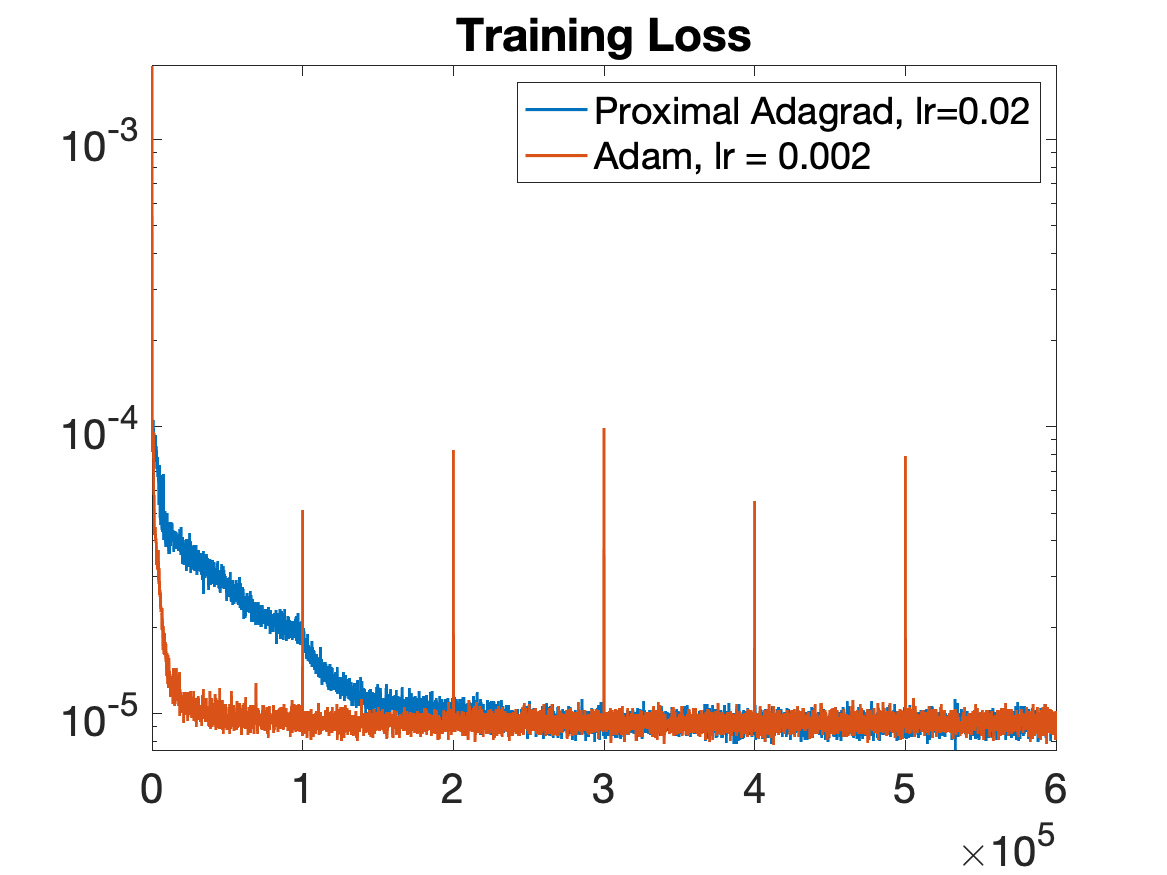}
   \caption{\eRNN}
\end{subfigure}
\caption{\mwrevise{\textbf{Training history using Adam v.s. Proximal Adagrad.} \fige{pwc}{16}{H^1}{{\RNN}/{\eRNN}}{20}{21}{0.02/0.002}. }}
    \label{fig:optimizer_comparison}
\end{figure}

\mwrevise{Here we present a comparison of the  Proximal Adagrad with Adam. From Figure \ref{fig:adma} we also observe that the training process of \eRNN ~is more stable than \RNN. From Figure \ref{fig:optimizer_comparison} and Table \ref{tab:adam}, we can further see that both optimizers seem to provide similar results by the end of the training process. But
Adam \mwrevise{seems to be overall a little less stable}.  The reason for taking the learning rate for Adam to be $0.001$ is because {larger rates like $0.01$ appear to produce} meaningless results.}


\begin{table}[!h]
	\centering
	\begin{tabular}{c|c|c c|cc }
	 	&&\multicolumn{2}{|c}{$\hat{\mE}$}&\multicolumn{2}{|c}{Relative $L_2$ Error of $u_{pred}$}\\
	 	\hline
		\# of blocks& learning rate& {eResNet}&{ResNet} & {eResNet}&{ResNet}\\
 		\hline
	    Proximal Adagrad&0.02 &7.76\%&8.42\%&47.42\%&49.28\%\\
	   Adam &0.002&7.47\%&7.82\%&47.14\%&  47.23\%\\
		\hline
	\end{tabular}
\hfill
    \caption{\textbf{ Generalization error v.s. different optimizer(fixed total training steps $6\times 10^5$).} \fige{pwc}{16}{H^1}{{eResNet}/{ResNet}}{20}{21}{0.02/0.002}. }
    \label{tab:adam}
\end{table}


\section{Comparison with Affine Reduced Basis Estimators}\label{sec:RB_compare}

\begin{figure}[h!]
	\centering
	\includegraphics[width = 0.6\textwidth]{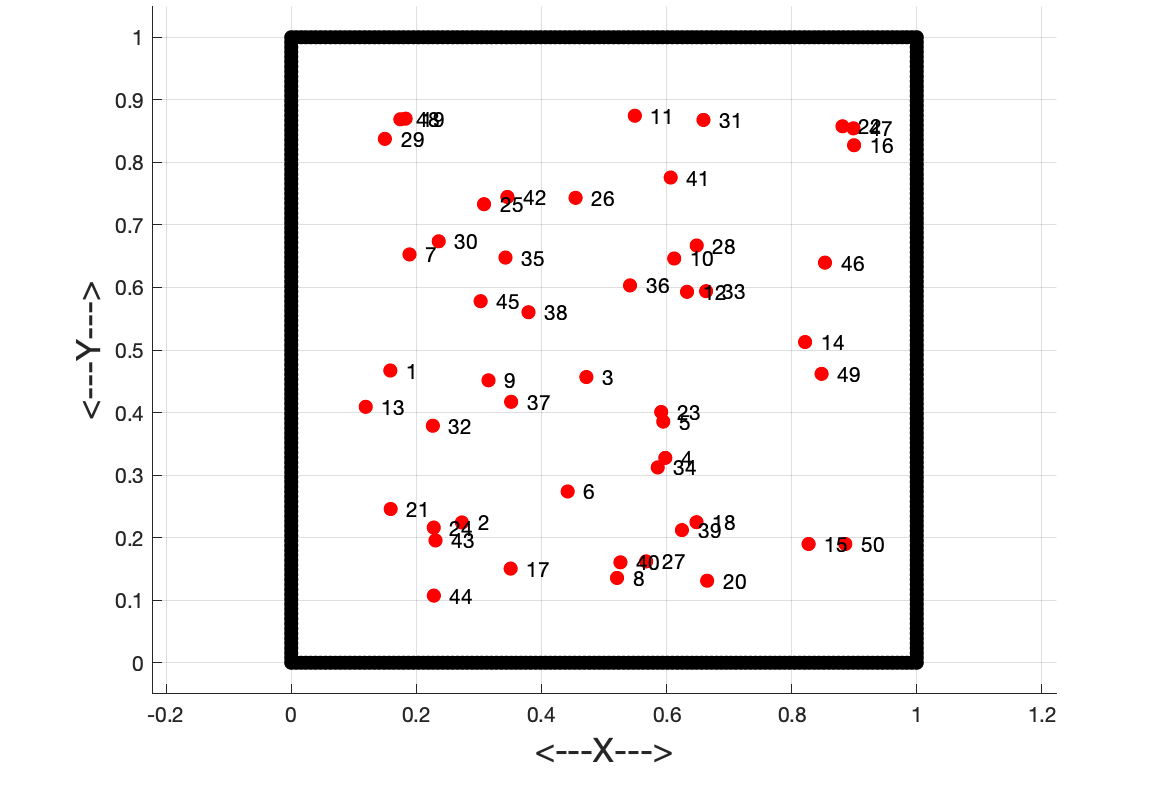}
	\caption{An illustration of $50$ random sensors.}
	\label{fig:random_sen}
\end{figure}

\begin{figure}[h!]
     \centering
         \centering
         \includegraphics[width=0.5\textwidth]{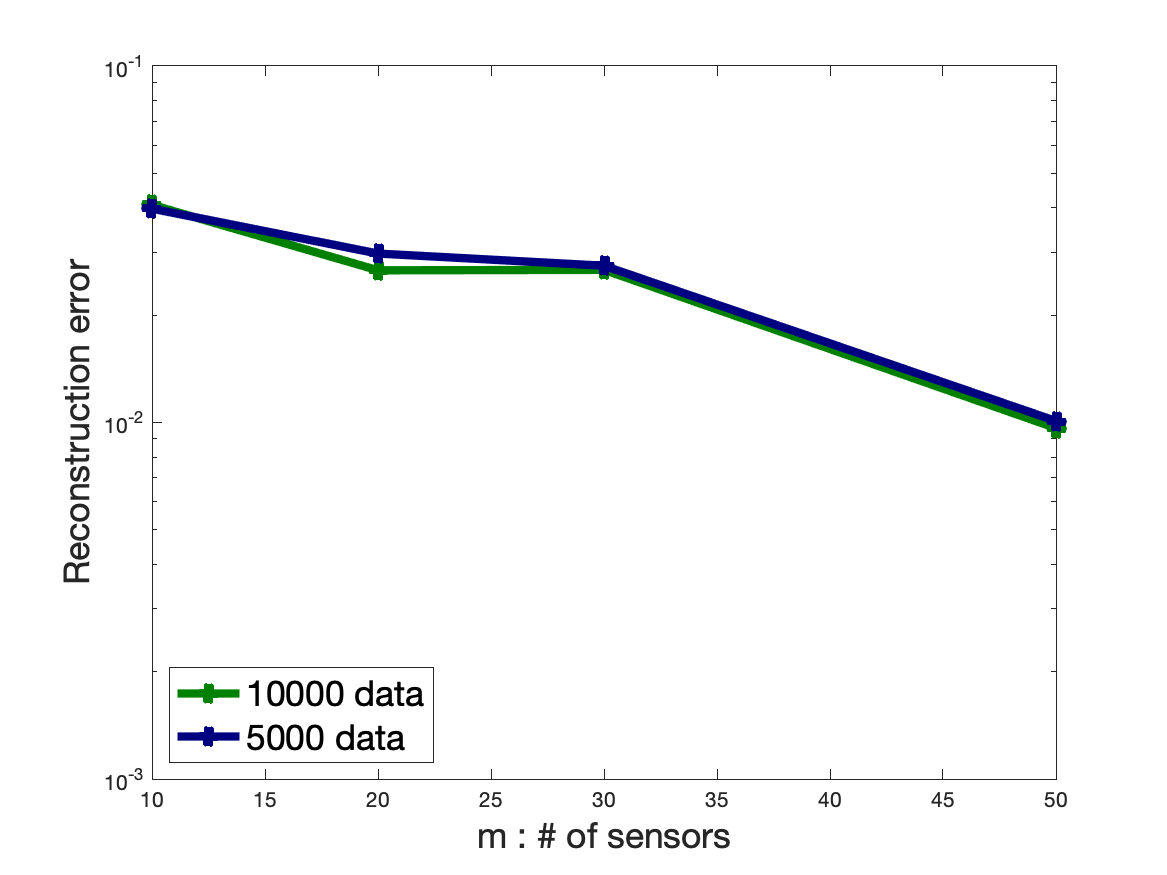}
        \caption{Max reconstruction $H^1$ errors among $1000$ testing samples with \mwrevise{\eRNN}.}
        \label{fig:compare_RB}
\end{figure}

As discussed  {earlier} in more detail,   we have run these experiments with varying choices of learning rates and architecture modifications consistently obtaining essentially the same magnitude of training loss and generalization error. This indicates a certain saturation effect as well as reliability in generating consistent results. Nevertheless, one wonders to what extent the actual expressive power of the networks is at least nearly exhausted and how to gauge the results in comparison with alternate approaches.
We have therefore compared the accuracy achieved by \mwrevise{\eRNN} with results for Affine Space estimators
from \mw{\cite{CDDN}}, mentioned earlier in Section \ref{ssec:1-space}. For the sake of such a comparison we consider the same type of up to $50$ randomly 
distributed 
sensors depicted in  Figure \ref{fig:random_sen}).
We generate again $10,000$ data points   $9000$ of which are are used for training while the rest 
is used for evaluation. 
Based on the experiences gained in previous experiments we have used a $\mathsf{B}2\hyp\mathsf{W}100$ \mwrevise{\eRNN} as a model for 
learning the observation-to-state mapping. That is, the \mwrevise{\RNN} of $2$ blocks is trained in the expansion manner. We confine the training of the \mwrevise{\eRNN} to a fixed number of $100,000$ steps. 
The dependence of the tested generalization error on the number of sensors   is shown in Figure  \ref{fig:compare_RB}. As expected, an increasing number of sensors provides more detailed
information on $P_{\Wp}(\cM)$. The green and blue curve show that further increasing the number of training data has little
effect on the achieved accuracy. For the given fixed budget of $10,000$ steps the evaluation
shows an $H^1$-error of about $10^{-2}$.

The performance of several versions of Affine Space estimators under the same test conditions 
has been reported in \cite{CDDN}. The computationally most expensive but also most accurate
Optimal Affine Space estimator achieves in this experiment roughly an $H^1$-error of size  {$3 \times 10^{-3}$} which is slightly better than the accuracy $8\times10^{-3}$,   
observed for \mwrevise{\eRNN}. However, when lifting the cap of at most $10,000$ training steps, 
the observed maximal $H^1$-error for \mwrevise{\eRNN} drops also to $5\times 10^{-3}$ in the $50$-sensor case,
which is at the same level as what the best affine estimator achieves.
In summary, it seems that for this type of problems both types of estimators achieve about the same
level of estimation accuracy and nonlinearity of the neural network lifting map does not seem to offer
substantial advantages for scenario (S1). There is instead a noteworthy difference regarding computational cost in relation to predictable training success. The above example shows that
neural networks may have  significant disadvantages with regard to optimization success and incurred computational cost in comparison with affine-space recovery schemes where, however, \mwrevise{\eRNN} shows a consistent level of reliability
that avoids degrading accuracy under over-parametrization.


\begin{figure}[hbt!]
     \centering
     \begin{subfigure}[b]{0.3\textwidth}
         \centering
         \includegraphics[width=1.2\textwidth]{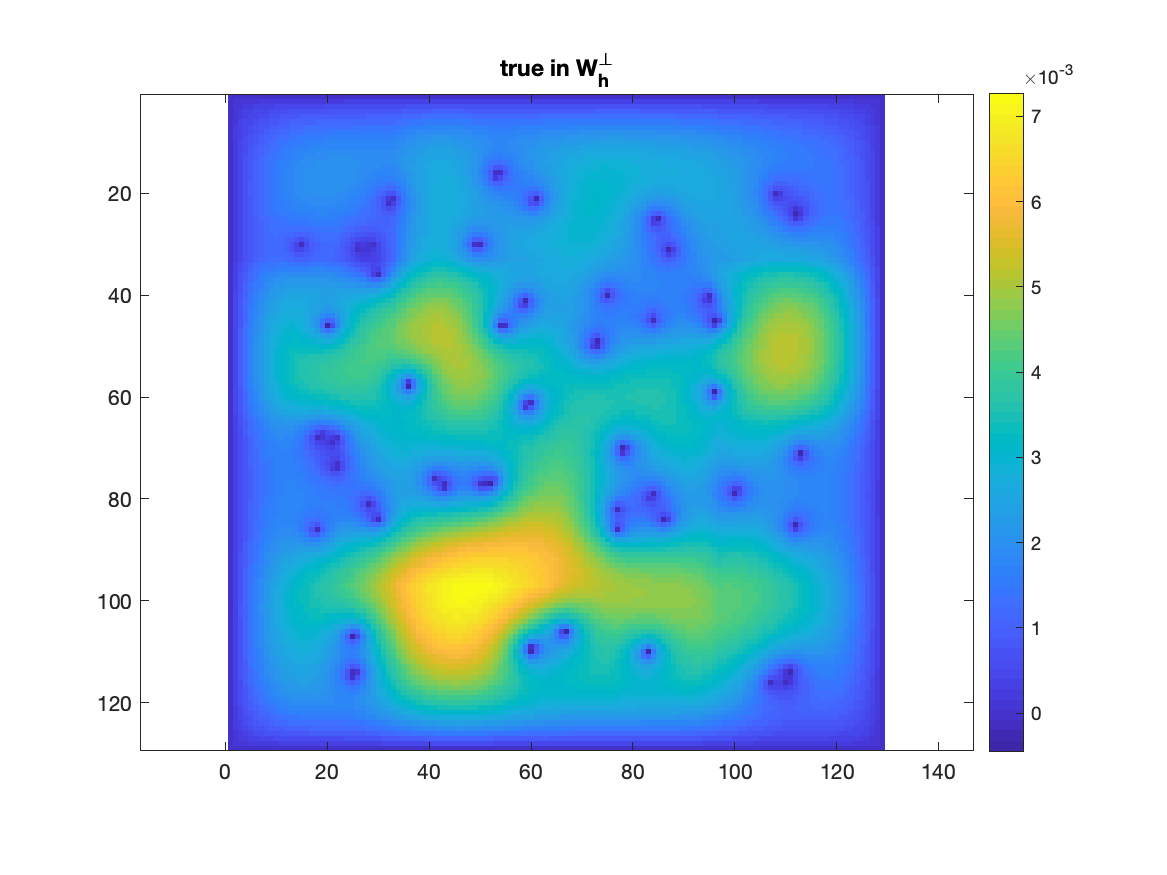}
         \caption{Exact solution}
     \end{subfigure}
     \hfill
     \begin{subfigure}[b]{0.3\textwidth}
         \centering
\includegraphics[width=1.2\textwidth]{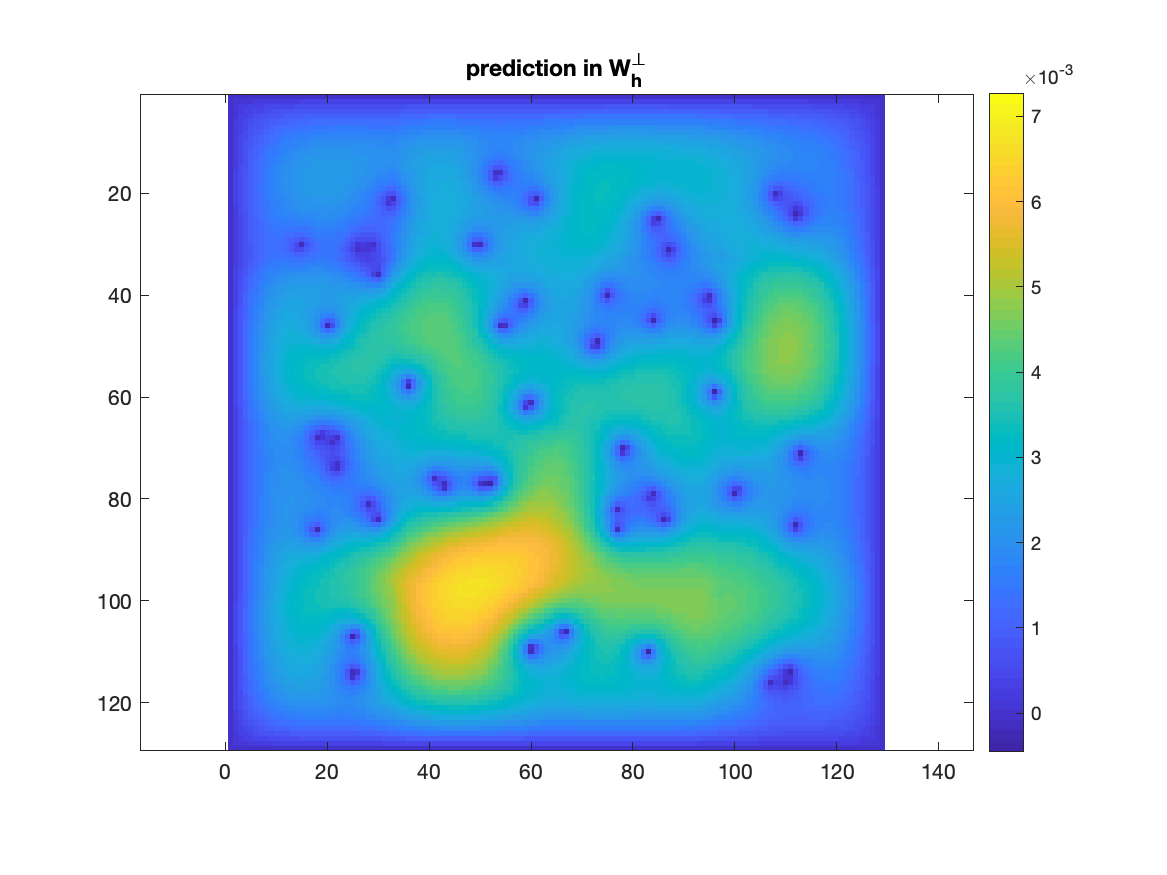}
         \caption{Prediction}
     \end{subfigure}
     \hfill
     \begin{subfigure}[b]{0.3\textwidth}
         \centering
\includegraphics[width=1.2\textwidth]{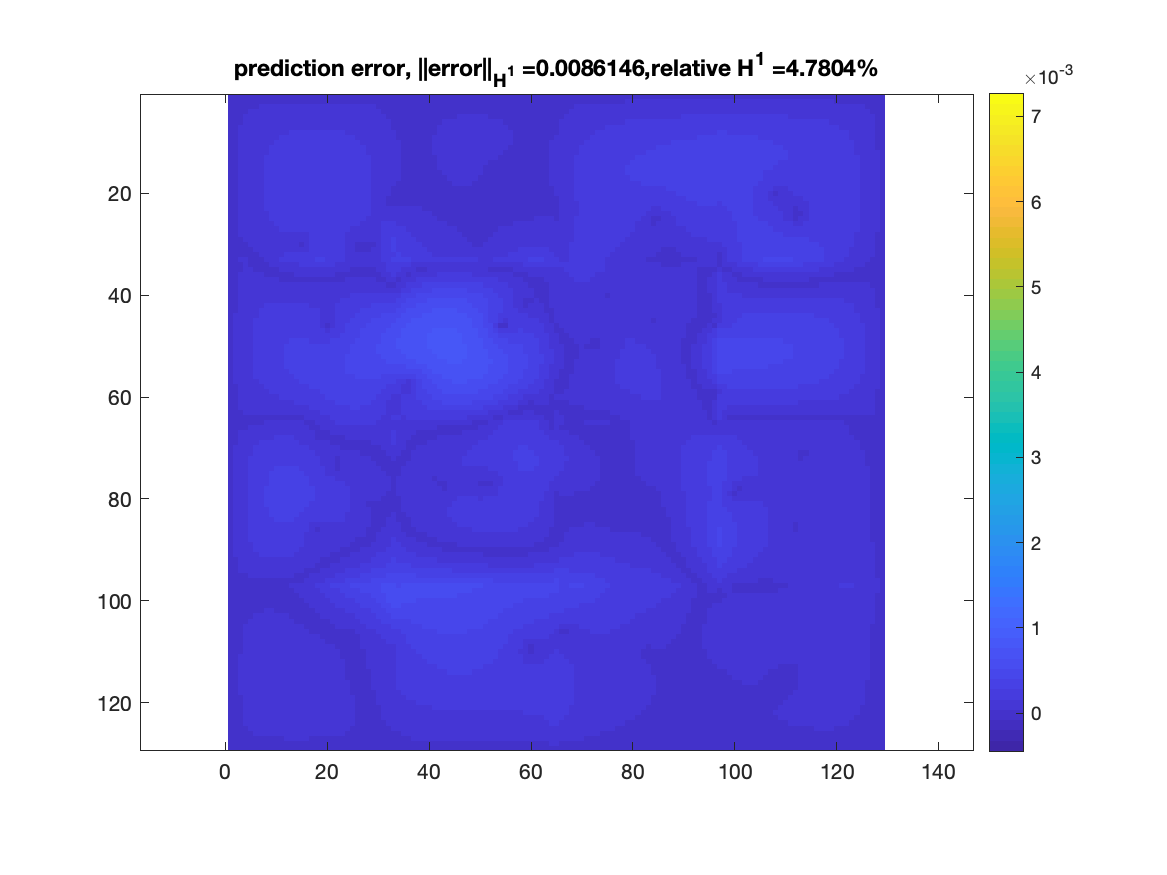}
         \caption{Error}
     \end{subfigure}
        \caption{An example of of exact vs. \mwrevise{\eRNN} predicted solution in $\mathbb{W}_h^T$ for 50 random sensors case. }
\end{figure}

\pagebreak

\section{Summary}\label{sec:conclusion}

\begin{figure}
    \centering
\begin{subfigure}[t]{\textwidth}
    \centering
    \includegraphics[width = 0.8\textwidth]{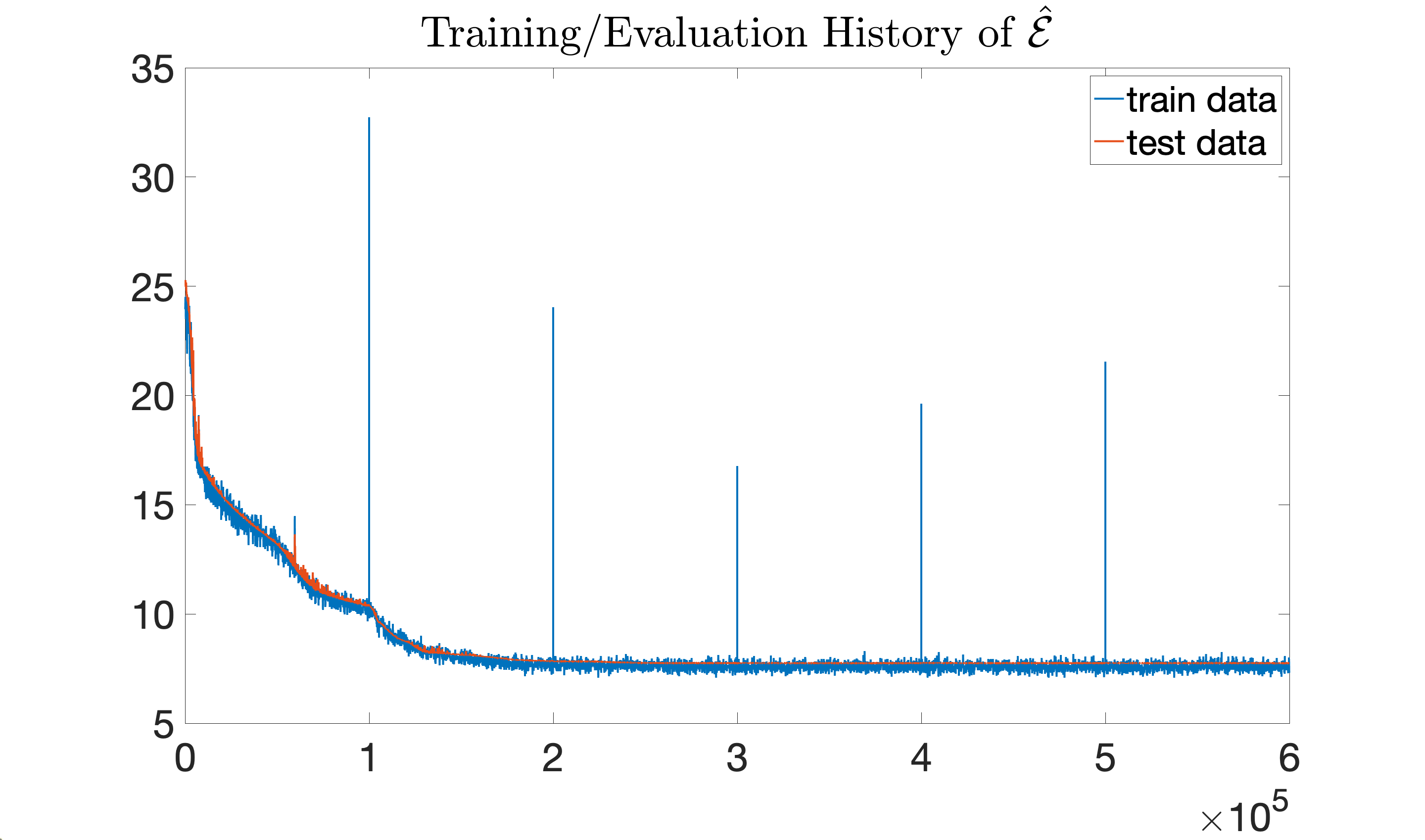}
    \caption{$\mathsf{B} 6$: Train/Evaluation Error history}
\end{subfigure}
\hfill
\begin{subfigure}[t]{\textwidth}
    \centering
     \includegraphics[width = 0.8\textwidth]{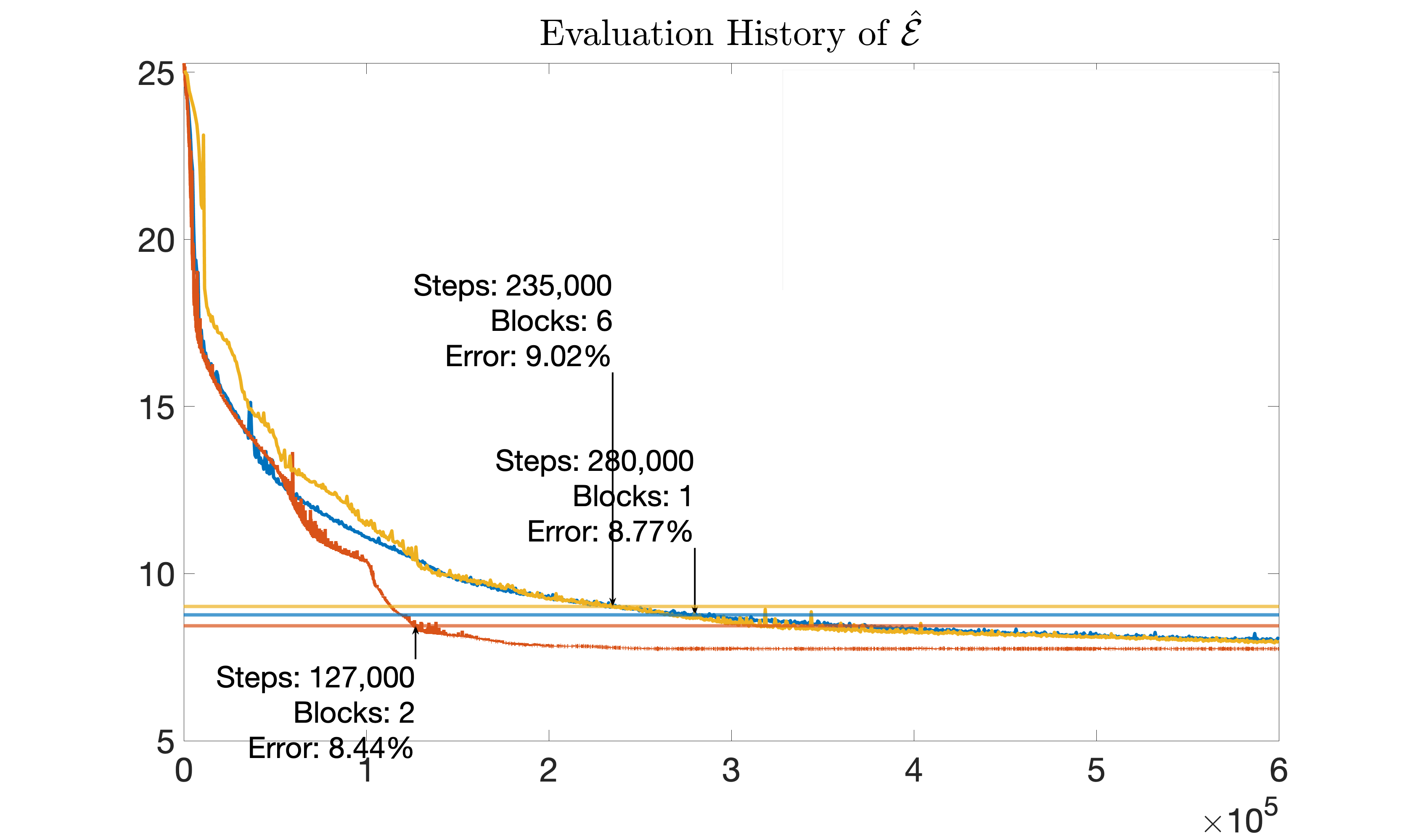}
     \caption{$\mathsf{B}1/\mathsf{B}6$: Evaluation Error History.  \mwrevise{Horizontal lines: $10/9$ level of the terminal evaluation error (blue: \RNN-$\mathsf{B}1$, red: \eRNN-$\mathsf{B}6$, yellow: \RNN-$\mathsf{B}6$)}}
\end{subfigure}
    \caption{$\mathtt{Train}$-\textit{log-normal}-$\mathsf{sen}16$-POD-$L_2$, $\mathsf{W}20$-$\mathsf{O}21$-$\mathsf{lr}0.02$.}
    \label{fig:history}
\end{figure}

Overall, the experiments for the two application scenarios (S1), (S2)
reflect the following general picture. Regardless of a specific training mode
accuracy improves rapidly at the beginning and essentially saturates at a moderate
network complexity. Beyond that point further 
improvements require a relatively substantial training effort which may not even be rewarded in the \mwrevise{\RNN} mode. Instead, \mwrevise{\eRNN} usually responds with slight improvements and essentially never with an accuracy degradation. While in a number of
cases \mwrevise{\eRNN} achieves a smaller generalization error than plain \mwrevise{\RNN}, by and large,
the differences in accuracy are not overly significant. Once the generalization error curve starts flattening, additional increases of network complexity
seem to just increase over-parametrization and widen a flat plateau fluctuating
around {``achievable''} local minima. {Instead a realization of theoretically 
possible expressive power seems to remain highly improbable.} 
Aside from an increased robustness with respect to algorithmic settings, the main advantage of \mwrevise{\eRNN} over \mwrevise{\RNN} seems to lie in substantial savings of computational work needed to nearly realize an apparently  achievable generalization accuracy. This is illustrated by Figure \ref{fig:history}, (b),
recording the work needed to achieve $10/9$ of the smallest generalization error
achieved by the respective training modality. It is also interesting to note that
the generalization errors at various optimization stages are not much larger 
than the corresponding relative loss-size, reflecting reliability of the schemes,
see Figure \ref{fig:history}, (a).

\section*{Acknowledgments}
This work was supported by National Science Foundation under grant DMS-2012469. We thank the reviewers for their valuable comments regarding the presentation of the material. 

\end{document}